\newif \ifSIAM

\ifSIAM
\documentclass[review,onefignum,onetabnum]{siamart190516}
\else
\documentclass[10pt,a4paper]{article}
\usepackage[a4paper,left=2.5cm,right=2.5cm,top=2.5cm,bottom=2.5cm]{geometry}
\fi

\usepackage{bm}
\usepackage{amssymb}
\usepackage{amsmath}
\usepackage{mathtools}
\usepackage{enumitem}
\usepackage[normalem]{ulem} 
\usepackage{cancel}
\usepackage[textsize=footnotesize,color=blue!30!white]{todonotes}
\usepackage{pgfplots}
\newcommand{\SlopeTriangle}[6]
{

    \pgfplotsextra
    {
        \pgfkeysgetvalue{/pgfplots/xmin}{\xmin}
        \pgfkeysgetvalue{/pgfplots/xmax}{\xmax}
        \pgfkeysgetvalue{/pgfplots/ymin}{\ymin}
        \pgfkeysgetvalue{/pgfplots/ymax}{\ymax}

        \pgfmathsetmacro{\xArel}{#1}
        \pgfmathsetmacro{\yArel}{#3}
        \pgfmathsetmacro{\xBrel}{#1-#2}
        \pgfmathsetmacro{\yBrel}{\yArel}
        \pgfmathsetmacro{\xCrel}{\xArel}

        \pgfmathsetmacro{\lnxB}{\xmin*(1-(#1-#2))+\xmax*(#1-#2)} 
        \pgfmathsetmacro{\lnxA}{\xmin*(1-#1)+\xmax*#1} 
        \pgfmathsetmacro{\lnyA}{\ymin*(1-#3)+\ymax*#3} 
        \pgfmathsetmacro{\lnyC}{\lnyA+#4*(\lnxA-\lnxB)}
        \pgfmathsetmacro{\yCrel}{\lnyC-\ymin)/(\ymax-\ymin)} 

        \coordinate (A) at (rel axis cs:\xArel,\yArel);
        \coordinate (B) at (rel axis cs:\xBrel,\yBrel);
        \coordinate (C) at (rel axis cs:\xCrel,\yCrel);

        \draw[#6]   (A)-- node[anchor=north] {#5}
                    (B)--
                    (C)--
                    cycle;
    }
}

\ifSIAM
\usepackage{lipsum}
\usepackage{amsfonts}
\ifpdf
  \DeclareGraphicsExtensions{.eps,.pdf,.png,.jpg}
\else
  \DeclareGraphicsExtensions{.eps}
\fi

\newsiamremark{remark}{Remark}
\newsiamremark{hypothesis}{Hypothesis}
\crefname{hypothesis}{Hypothesis}{Hypotheses}
\newsiamthm{claim}{Claim}
\else
\usepackage{amsthm}
\usepackage[breaklinks,bookmarks=false]{hyperref}
\fi

\hypersetup{colorlinks, linkcolor=blue, citecolor=blue,
urlcolor=blue, plainpages=false, pdfwindowui=false, pdfstartview={FitH}, pdftitle={$p$-robust equilibrated flux reconstruction in H(curl) based on local minimizations; application to a posteriori analysis of the curl--curl problem}, pdfauthor={Martin Vohral\'ik} }

\numberwithin{equation}{section}

\newif \ifbbm

\newif \iffract

\newcommand\cred[1]{%
  \protect\leavevmode
  \begingroup
    #1%
  \endgroup
}

\newcommand\tcf[1]{%
  \protect\leavevmode
  \begingroup
    #1%
  \endgroup
}

\newcommand\cye[1]{%
  \protect\leavevmode
  \begingroup
    #1%
  \endgroup
}

\newcommand\cgr[1]{%
  \protect\leavevmode
  \begingroup
    #1%
  \endgroup
}

\newcommand{\be}{\begin{equation}}
\newcommand{\ee}{\end{equation}}
\newcommand{\bea}{\begin{eqnarray}}
\newcommand{\eea}{\end{eqnarray}}
\newcommand{\bean}{\begin{eqnarray*}}
\newcommand{\eean}{\end{eqnarray*}}
\def\ba#1\ea{\begin{align}#1\end{align}}
\def\ban#1\ean{\begin{align*}#1\end{align*}}
\def\bat#1\eat{\begin{alignat}#1\end{alignat}}
\def\batn#1\eatn{\begin{alignat*}#1\end{alignat*}}
\def\bs#1\es{\begin{split}#1\end{split}}
\newcommand{\bse}{\begin{subequations}}
\newcommand{\ese}{\end{subequations}}
\newcommand{\bt}{\begin{theorem}}
\newcommand{\et}{\end{theorem}}
\newcommand{\bpr}{\begin{proposition}}
\newcommand{\epr}{\end{proposition}}
\newcommand{\bcj}{\begin{conjecture}}
\newcommand{\ecj}{\end{conjecture}}
\newcommand{\bl}{\begin{lemma}}
\newcommand{\el}{\end{lemma}}
\newcommand{\bc}{\begin{corollary}}
\newcommand{\ec}{\end{corollary}}
\newcommand{\bp}{\begin{proof}}
\newcommand{\ep}{\end{proof}}
\newcommand{\bd}{\begin{definition}}
\newcommand{\ed}{\end{definition}}
\newcommand{\br}{\begin{remark}}
\newcommand{\er}{\end{remark}}
\newcommand{\bas}{\begin{assumption}}
\newcommand{\eas}{\end{assumption}}
\newcommand{\bex}{\begin{example}}
\newcommand{\eex}{\end{example}}
\newcommand{\bqo}{\begin{quote}}
\newcommand{\eqo}{\end{quote}}
\newcommand{\bdc}{\begin{description}}
\newcommand{\edc}{\end{description}}
\newcommand{\bi}{\begin{itemize}}
\newcommand{\ei}{\end{itemize}}
\newcommand{\ben}{\begin{enumerate}}
\newcommand{\een}{\end{enumerate}}

\ifSIAM
\newtheorem{assumption}[theorem]{Assumption}

\newtheorem{conjecture}[theorem]{Conjecture}

\else

\newtheorem{lemma}{Lemma}[section]
\newtheorem{corollary}[lemma]{Corollary}
\newtheorem{proposition}[lemma]{Proposition}
\newtheorem{assumption}[lemma]{Assumption}
\newtheorem{theorem}[lemma]{Theorem}
\newtheorem{definition}[lemma]{Definition}
\newtheorem{remark}[lemma]{Remark}
\newtheorem{example}[lemma]{Example}

\newtheorem{conjecture}[lemma]{Conjecture}

\fi

\newcommand\Om{\Omega}
\newcommand\om{\omega}
\newcommand\oma{{\omega_\ver}} 
\newcommand\ome{{\omega_\edg}} 
\newcommand\toma{{\widetilde\omega_\ver}}
\newcommand\tome{{\widetilde\omega_\edg}}
\newcommand\omat{{\omega_\vertt}}
\newcommand\GD{{\Gamma_{\mathrm D}}} 
\newcommand\GN{{\Gamma_{\mathrm N}}} 
\newcommand\GT{{\Gamma_{\mathrm T}}} 
\newcommand\gD{{\gamma_{\mathrm D}}} 

\newcommand\Gr{\nabla} 
\newcommand\Dv{\nabla {\cdot}} 
\newcommand\Crl{\nabla {\times}} 

\newcommand\dv{\mathrm{div}}

\newcommand\crl{\mathrm{curl}}
\newcommand\scp{{\cdot}} 
\newcommand\vp{{\times}} 

\DeclarePairedDelimiter\norm{\|}{\|} 

\newcommand\pp{{p^\prime}}
\newcommand\qq{{q^\prime}}

\ifbbm

\else

\fi

\newcommand\HO{H^1}
\newcommand\Ho{H^1(\Om)}
\newcommand\tHo{{\bm H}^1(\Om)}
\newcommand\Hoi[1]{H^1(#1)}
\newcommand\tHoi[1]{{\bm H}^1(#1)}
\newcommand\Hois[1]{H^1_*(#1)}

\newcommand\HoD{H^1_{0,\mathrm{D}}(\Om)}

\newcommand\Lt{L^2(\Om)}
\newcommand\LT{L^2}
\newcommand\tLt{\bm{L}^2(\Om)}
\newcommand\tLti[1]{\bm{L}^2(#1)}

\newcommand\Lti[1]{L^2(#1)}

\newcommand\HDV{\bm{H}(\dv)}
\newcommand\Hdv{\bm{H}(\dv,\Om)}
\newcommand\HdvN{\bm{H}_{0,\mathrm{N}}(\dv,\Om)}

\newcommand\Hdvi[1]{\bm{H}(\dv,#1)}

\newcommand\Hdva{\bm{H}_{0}(\dv, \oma)}
\newcommand\Hdvat{\bm{H}_{0}(\dv, \toma)}

\newcommand\HC{\bm{H}(\crl)}
\newcommand\Hc{\bm{H}(\crl, \Om)}
\newcommand\Hci[1]{\bm{H}(\crl, #1)}
\newcommand\HcT{\bm{H}_{0,\mathrm{T}}(\crl, \Om)}
\newcommand\HcD{\bm{H}_{0,\mathrm{D}}(\crl, \Om)}
\newcommand\HcN{\bm{H}_{0,\mathrm{N}}(\crl, \Om)}
\newcommand\Hca{\bm{H}_{0}(\crl, \oma)}

\newcommand\Hcdual[1]{\bm{H}^\dagger(\crl, #1)}

\newcommand\RT{\bm{\mathcal{R\hspace{-0.1em}T}}\hspace{-0.25em}}
\newcommand\RTproj[1]{\bm{I}^{\bm{\mathcal{R\hspace{-0.1em}T}}}_{\cye{\elm,} #1}}
\newcommand\ND{\bm{\mathcal{N}}\hspace{-0.2em}}
\newcommand\NDproj[1]{\bm{I}^{\bm{\mathcal{N}}}_{#1}}

\newcommand{\cohom}{{\bm{\mathcal H}}(\Om,\GD)}

\ifbbm

\else

\fi

\newcommand\ie{i.e.}
\newcommand\cf{cf.}
\newcommand\eg{e.g.}
\newcommand\eal{{\em et al.}}
\newcommand\eq{:=}

\newcommand\ls{\lesssim}

\newcommand\nn{\nonumber}
\newcommand\pt{\partial}
\newcommand{\<}{\langle}
\renewcommand{\>}{\rangle}

\newcommand{\suma}{\sum_{\ver \in \Vh}}

\newcommand\nv{\bm 0}
\newcommand\reff[2]{\stackrel{\eqref{#1}}{#2}}

\ifbbm

\else

\fi

\newcommand{\elm}{{K}} 

\newcommand{\elmi}{{K_i}} 
\newcommand{\elmj}{K_j} 
\newcommand{\elmjo}{{K_{j_1}}} 
\newcommand{\elmjt}{{K_{j_2}}} 
\newcommand{\elml}{{K_{|\Ta|}}} 
\newcommand{\elmo}{{K_1}} 
\newcommand{\sd}{{F}} 
\newcommand{\sdt}{{F'}}
\newcommand{\edg}{{e}} 

\newcommand{\ver}{{\bm{a}}} 
\newcommand{\vertt}{{\bm{b}}}

\newcommand\FK{\mathcal{F}_\elm}

\newcommand\EhD{\mathcal{E}_h^{\mathrm{D}}} 

\newcommand\Th{{\mathcal{T}_h}} 
\newcommand\Ta{{\mathcal{T}_{\ver}}} 
\newcommand\tTa{{\widetilde{\mathcal{T}}_{\ver}}} 

\newcommand\Vh{\mathcal{V}_h} 

\newcommand\VK{\mathcal{V}_\elm}

\newcommand\psia{\psi^\ver}
\newcommand\psib{\psi^\vertt}
\newcommand\tAh{\tA_h}

\newcommand\fr{\thh} 
\newcommand\frh{\thh_h} 

\newcommand\tdah{\hat\bdl_h^\ver}
\newcommand\tda{\bdl_h^\ver}
\newcommand\tha{\bm{\theta}_h^\ver}
\newcommand\thatt{\bm{\theta}_h^\vertt}
\newcommand\thab{\bar{\bm{\theta}}_h^\ver}
\newcommand\tza{\bm{\epsilon}_h^\ver}
\newcommand\tzo{\bm{\epsilon}_h^1}
\newcommand\tzi{\bm{\epsilon}_h^i}
\newcommand\tzj{\bm{\epsilon}_h^j}
\newcommand\tzl{\bm{\epsilon}_h^{|\Ta|}}
\newcommand\tzjo{\bm{\epsilon}_h^{j_1}}
\newcommand\tzjt{\bm{\epsilon}_h^{j_2}}

\newcommand\tzab{\bar{\bm{\epsilon}}_h^\ver}
\newcommand\lhi{\lambda_h^i}
\newcommand\lho{\lambda_h^1}
\newcommand\lhj{\lambda_h^j}
\newcommand\lhij{\lambda_h^{i,j}}

\newcommand\pphi{p_h^i}
\newcommand\ppho{p_h^1}

\ifbbm

\else

\fi

\newcommand\vf{\varphi}

\newcommand\tx{\bm{x}}

\newcommand\thh{\bm{h}}

\newcommand\tj{\bm{j}}

\newcommand\tn{\bm{n}}

\newcommand\tq{\bm{q}}
\newcommand\tr{\bm{r}}

\newcommand\tv{\bm{v}}
\newcommand\tw{\bm{w}}

\newcommand\tA{\bm{A}}

\newcommand\tH{\bm{H}}

\newcommand\tK{\bm{K}}
\newcommand\tL{\bm{L}}

\newcommand\tV{\bm{V}}

\newcommand{\bPi}{\bm{\Pi}}

\newcommand{\bdl}{\bm{\delta}}
\newcommand\bvf{\bm{\varphi}}

\newcommand\bpsi{\bm{\psi}}
\newcommand\btau{\bm{\tau}}
\newcommand\biot{\bm{\iota}}

\newcommand\bvphi{\bm{\varphi}}

\newcommand\F{\mathcal{F}}

\newcommand\T{\mathcal{T}}

\newcommand\RR{{\mathbb R}}

\ifbbm

\newcommand\PP{{\mathbb P}}

\else

\newcommand\PP{{\mathcal P}}

\fi

\iffract

\newcommand\ft{{\frac 1 2}}
\newcommand\mft{{-\frac 1 2}}

\else

\newcommand\ft{{1/2}}
\newcommand\mft{{-1/2}}

\fi

\newcommand\Hsa{H^1_*(\oma)}
\newcommand\Hsaw{H^1_*(\toma)}

\newcommand\Clift{C_{\mathrm{L}}}

\newcommand{\htzab}{\widehat{\bm \epsilon}}

\newcommand{\mean}{\bm{\mu}^\ver}
\newcommand{\tmean}{\widetilde {\bm{\mu}}^\ver}
\newcommand{\JAC}{\mathbb J}

\ifSIAM

\headers{$\MakeLowercase{p}$-robust equilibrated flux reconstruction in $\bm{H}(\MakeLowercase{\mathrm{curl})}$}{T. Chaumont-Frelet and M. Vohral\'ik}

\title{$\MakeLowercase{p}$-robust equilibrated flux reconstruction in $\HC$ based on local minimizations. Application to a~posteriori analysis of the curl--curl problem\thanks{Submitted to the editors DATE.
\funding{This project has received funding from the European Research Council (ERC) under the European Union's Horizon 2020 research and innovation program (grant agreement No 647134 GATIPOR).}}}

\author{%
Th\'eophile Chaumont-Frelet\thanks{Inria, 2004 Route des Lucioles, 06902 Valbonne, France \&
Laboratoire J.A. Dieudonn\'e, Parc Valrose, 28 Avenue Valrose, 06108 Nice, France
(\email{theophile.chaumont@inria.fr}, \url{https://tchaumont.github.io}).}
\and
Martin Vohral\'ik\thanks{Inria, 2 rue Simone Iff, 75589 Paris, France \&
CERMICS, Ecole des Ponts, 77455 Marne-la-Vall\'ee, France (\email{martin.vohralik@inria.fr}, \url{https://who.rocq.inria.fr/Martin.Vohralik}).}
}

\else

\title{$p$-robust equilibrated flux reconstruction in $\HC$ based on local minimizations. Application to a~posteriori analysis of the curl--curl problem\thanks{This project has received funding from the European Research Council (ERC) under the European Union's Horizon 2020 research and innovation program (grant agreement No 647134 GATIPOR).}}

\author{Th\'eophile Chaumont-Frelet\footnotemark[2]\,\,\footnotemark[3] \and Martin Vohral\'ik\footnotemark[4]\,\,\footnotemark[5]}

\fi

\def\KW{Sobolev space $\HC$,
Sobolev space $\HDV$,
equilibrated flux reconstruction,
$p$-robustness,
a posteriori error estimate,
divergence-free decomposition,
broken polynomial extension}

\begin{document}
\maketitle

\ifSIAM

\newtheorem{theoremclr}[theorem]{\cye{Theorem}}

\else

\newtheorem{theoremclr}[lemma]{\cye{Theorem}}
\renewcommand{\thefootnote}{\fnsymbol{footnote}}

\footnotetext[2]{Inria, 2004 Route des Lucioles, 06902 Valbonne, France
(\href{mailto:theophile.chaumont@inria.fr}{\texttt{theophile.chaumont@inria.fr}}).}
\footnotetext[3]{Laboratoire J.A. Dieudonn\'e, Parc Valrose, 28 Avenue Valrose, 06108 Nice, France}
\footnotetext[4]{Inria, 2 rue Simone Iff, 75589 Paris, France
(\href{mailto:martin.vohralik@inria.fr}{\texttt{martin.vohralik@inria.fr}}).}
\footnotetext[5]{CERMICS, Ecole des Ponts, 77455 Marne-la-Vall\'ee, France}

\renewcommand{\thefootnote}{\arabic{footnote}}

\fi

\begin{abstract}
\ifSIAM\else\noindent\fi
We present a local construction of $\HC$-conforming piecewise polynomials
satisfying a prescribed curl constraint. We start from a piecewise polynomial
not contained in the $\HC$ space but satisfying a suitable orthogonality property.
The procedure employs minimizations in vertex patches and the outcome is,
up to a generic constant independent of the underlying polynomial degree,
as accurate as the best-approximations over the entire local versions of $\HC$.
This allows to design guaranteed, fully computable, constant-free, and
polynomial-degree-robust a posteriori error estimates of the Prager--Synge type for
N\'ed\'elec's finite element approximations of the curl--curl problem. A divergence-free decomposition of a divergence-free $\HDV$-conforming piecewise polynomial, relying on  over-constrained minimizations in Raviart--Thomas' spaces, is the key ingredient.
Numerical results illustrate the theoretical developments.
\end{abstract}

\ifSIAM 

\else
\fi

\ifSIAM
\else


\fi

\section{Introduction} \label{sec_intr}

A posteriori error estimation by equilibrated flux reconstruction has achieved a great
attention for \cye{elliptic model} problems like the Poisson problem. For an $\HO$-conforming
discretization whose flux is not in $\HDV$, one has to reconstruct a flux in $\HDV$ satisfying a
prescribed divergence constraint. To design high-performance algorithms, the procedure must
furthermore be localized and can not involve a solution of any supplementary global problem.
Then, a guaranteed, fully computable, and constant-free upper bound on the
unknown discretization error follows from the equality of Prager and Synge~\cite{Prag_Syng_47}.
There are several techniques of such an equilibrated flux reconstruction. Following Ladev{\`e}ze
and Leguillon~\cite{Lad_Leg_83} and Ainsworth and Oden~\cite{Ainsw_Oden_a_post_FE_00}, normal
fluxes on mesh faces can first be constructed and then lifted elementwise \cye{as in Nicaise~\eal\ \cite{Nic_Wit_Wohl_a_post_Lam_eq_fl_08}},
dual Vorono\"i-type grids can be employed for local non-overlapping minimizations in $\HDV$
as in Luce and Wohlmuth~\cite{Luce_Wohl_local_a_post_fluxes_04} or Hannukainen
\eal~\cite{Han_Sten_Voh_a_post_S_12}, or a localization by the partition of unity
via the finite element hat basis functions can be used for an overlapping combination of
best-possible vertex-patch fluxes as in Destuynder and M{\'e}tivet~\cite{Dest_Met_expl_err_CFE_99} or
Braess and Sch{\"o}berl~\cite{Braess_Scho_a_post_edge_08}. This last approach is conceptual and,
as established in Braess \eal~\cite{Brae_Pill_Sch_p_rob_09} and Ern and
Vohral\'\i k~\cite{Ern_Voh_p_rob_3D_20}, it gives estimates robust with respect to the polynomial
degree $p$ (henceforth termed $p$-robust).

In contrast, there is only a handful of results available for the curl--curl problem, where,
for an $\HC$-conforming discretization whose curl is not in $\HC$, one has to locally reconstruct
a flux in $\HC$ satisfying a prescribed curl constraint. An approach based on patchwise
minimizations for the lowest-order case $p=0$ has been designed
in~\cite{Braess_Scho_a_post_edge_08}. Its generalization for arbitrary $p \geq 1$,
however, turns surprisingly difficult and, to the best of our knowledge, has not been presented
yet. Several workarounds appeared in the literature recently, though. A conceptual discussion
appears in Licht~\cite{Licht_FEEC_a_post_H_curl_19}, whereas a construction following in
spirit~\cite{Lad_Leg_83, Ainsw_Oden_a_post_FE_00} has been proposed and analyzed in Gedicke
\eal~\cite{Ged_Gee_Per_a_post_Maxw_20}. This last approach has been recently modified in
Gedicke \eal~\cite{Ged_Gee_Per_Sch_post_Maxw_21} in order to achieve $p$-robustness.
A broken patchwise equilibration procedure that bypasses the Prager--Synge theorem
is proposed and proved $p$-robust in Chaumont-Frelet \eal~\cite{Chaum_Ern_Voh_Maxw_22};
it relies on smaller edge patches, but the arising estimates are not constant-free.

The purpose of this contribution is to design an equilibrated flux reconstruction in $\HC$
employing best-possible local fluxes. In doing so, we rely on localization by the partition
of unity via the hat functions and overlapping flux combinations, in generalization of
the concept of~\cite{Braess_Scho_a_post_edge_08} to arbitrary $p \geq 0$. Consequently,
we identify the equivalent in $\HC$ of the concept of equilibrated flux reconstruction in
$\HDV$ from~\cite{Dest_Met_expl_err_CFE_99, Braess_Scho_a_post_edge_08, Brae_Pill_Sch_p_rob_09,
Ern_Voh_p_rob_15, Ern_Voh_p_rob_3D_20}.
This is then used for a posteriori error estimation when the N\'ed\'elec (edge) finite elements
of arbitrary degree $p \geq 0$ are used for approximation of the curl--curl problem.
It leads to guaranteed, fully computable, and constant-free a posteriori error estimates that are locally efficient and robust with respect to the polynomial degree $p$; (higher-order) data oscillation terms \cye{are rigorously included in our analysis}.
Our $p$-robust efficiency proofs are based on the seminal volume and tangential trace
$p$-robust extensions on a single tetrahedron of
Costabel and McIntosh~\cite[Proposition~4.2]{Cost_McInt_Bog_Poinc_10} and Demkowicz
\eal~\cite[Theorem~7.2]{Demk_Gop_Sch_ext_II_09}.
These results were recently extended into a stable broken polynomial extension for
\cye{a single tetrahedron in Chaumont-Frelet~\eal\ \cite[Theorem~2]{Chaum_Ern_Voh_curl_elm_20},}
an edge patch of tetrahedra \cye{in Chaumont-Frelet~\eal\ \cite[Theorem~3.1]{Chaum_Ern_Voh_Maxw_22}}\cye{, and for a vertex patch of tetrahedra in Chaumont-Frelet and Vohral\'ik~\cite[Theorem~3.3, see also Corollary~4.3]{Chaum_Voh_p_rob_3D_H_curl_22}.}

An important step in the construction of our estimators is to decompose the given
divergence-free right-hand side into locally supported divergence-free contributions.
Starting from the available (lowest-order Galerkin) orthogonality property, we propose
a multi-stage procedure relying on two central technical results of independent interest:
over-constrained minimization in Raviart--Thomas' spaces leading to suitable elementwise orthogonality properties, and a decomposition of a divergence-free piecewise polynomial with the above elementwise orthogonality properties into local divergence-free contributions.
These issues are related to the developments on divergence-free decompositions in
Scheichl~\cite{Schei_MFE_dec_03}, Alonso Rodr\'{\i}guez
\eal~\cite{AlonsRodr_Cam_San_div_free_18,AlonsRodr_Cam_Ghi_Val_div_free_17},
and the references therein.

This contribution is organized as follows. Section~\ref{sec_not} fixes the \cye{notation}.
Section~\ref{sec_set} \cye{introduces the curl--curl problem\cye{, its N\'ed\'elec finite element discretization,} and identifies therefrom two abstract assumptions under which our analysis is performed.}
In Section~\ref{sec_mot}, we motivate our approach at the continuous level. \cye{Section~\ref{sec_main_res} then presents our main results:} in Section~\ref{sec_div_free_dec}, we develop a divergence-free decomposition of the given target curl; in Section~\ref{sec_equil}, we present the equilibrated flux reconstruction based on local minimization in $\HC$, as well as its $p$-robust stability; finally, these abstract results are applied in Section~\ref{sec_a_post} to the N\'ed\'elec finite element discretization of the curl--curl problem. Section~\ref{sec_num} is dedicated to a numerical illustration\cye{, whereas Section~\ref{sec_det_proofs} collects some technical details and proofs}. Finally, in Appendices~\ref{app_over_constr_min}
and~\ref{app_div_free_dec}, we present the two central technical results on over-constrained
minimization and divergence-free decomposition.

\section{Notation} \label{sec_not}

\cye{The purpose of this section is to set the necessary notation.}
Let $\om, \Om \subset \RR^3$ be open, Lipschitz polyhedra; $\Om$ will be
used to denote the computational domain, while we reserve the
notation $\om \subseteq \Om$ for its simply connected subsets.
Notice that we do not require $\Om$ to be simply connected.

\subsection{Sobolev spaces \texorpdfstring{$\HO$, $\HC$, and $\HDV$}{H1, H(curl), and H(div)}} \label{sec_bas_cont_sp}

We let $\Lti{\om}$ be the space of scalar-valued square-integrable functions defined on $\om$;
we use the notation $\tLti{\om} \eq [\Lti{\om}]^3$ for vector-valued functions with each
component in $\Lti{\om}$. We denote by $\norm{{\cdot}}_\om$ the $\Lti{\om}$ or $\tLti{\om}$
norm and by $({\cdot},{\cdot})_\om$ the corresponding scalar product; we drop the index when
$\om = \Om$. We will extensively work with the following three Sobolev spaces: 1) $\Hoi{\om}$,
the space of scalar-valued $\Lti{\om}$ functions with weak gradients in $\tLti{\om}$,
$\Hoi{\om} \eq \{v \in \Lti{\om}; \, \Gr v \in \tLti{\om}\}$; 2) $\Hci{\om}$,
the space of vector-valued $\tLti{\om}$ functions with weak curls in $\tLti{\om}$,
$\Hci{\om} \eq \{\tv \in \tLti{\om}; \, \Crl \tv \in \tLti{\om}\}$; and 3) $\Hdvi{\om}$,
the space of vector-valued $\tLti{\om}$ functions with weak divergences in $\Lti{\om}$,
$\Hdvi{\om} \eq \{\tv \in \tLti{\om}; \, \Dv \tv \in \Lti{\om}\}$. We refer the reader to
Adams~\cite{Adams_75} and Girault and Raviart~\cite{Gir_Rav_NS_86} for an in-depth description
of these spaces. Moreover, component-wise $\Hoi{\om}$ functions will be denoted by
$\tHoi{\om} \eq \{\tv \in \tLti{\om}; \, \tv_i \in \Hoi{\om}, \, i=1, \ldots, 3\}$.
We will employ the notation $\<{\cdot},{\cdot}\>_S$ for the integral product on boundary (sub)sets $S \subset \pt \om$.

\subsection{\cye{Sobolev spaces with partially vanishing traces on \texorpdfstring{$\pt\Om$}{}}} \label{sec_cont_sp_BC}

Let $\GD$, $\GN$ be two disjoint, relatively open, and possibly empty subsets of the
computational domain boundary $\pt \Om$ such that $\pt \Om = \overline \GD \cup \overline \GN$.
We assume in addition that each boundary face of the mesh $\Th$ \cye{defined below} lies entirely either
in $\overline \GD$ or in $\overline \GN$. Then $\HoD$ is the subspace of $\Ho$ formed by
functions vanishing on $\GD$ in the sense of traces, $\HoD \eq \{v \in \Ho; \, v=0$ on $\GD\}$.
Let $\tn_\Om$ be the unit normal vector on $\pt \Om$, outward to $\Om$.
Let $\mathrm{T} = \mathrm{D}$ or $\mathrm{N}$; then $\HcT$ is the subspace of $\Hc$ formed by
functions with vanishing tangential trace on $\GT$,
$\HcT \eq \{\tv \in \Hc; \, \tv \vp \tn_{\Om}=\cye{\nv}$ on $\GT\}$, where $\tv \vp \tn_{\Om} = \cye{\nv}$ on $\GT$ means that $(\Crl \tv, \bvf) - (\tv, \Crl \bvf) = 0$ for all functions $\bvf \in \tHo$ such that $\bvf \vp \tn_{\Om} = \nv$ on $\pt \Om \setminus \GT$.
Finally, $\HdvN$ is the subspace of $\Hdv$ formed by functions with vanishing normal
trace on $\GN$, $\HdvN \eq \{\tv \in \Hdv; \, \tv \scp \tn_{\Om}=0$ on $\GN\}$, where
$\tv \scp \tn_{\Om} = 0$ on $\GN$ means that $(\tv, \Gr \vf) + (\Dv \tv, \vf) = 0$ for
all functions $\vf \in \HoD$. Fernandes and Gilardi~\cite{Fer_Gil_Maxw_BC_97} present a
thorough characterization of tangential (resp. normal) traces of $\Hc$ (resp. $\Hdv$) on
a part of the boundary $\pt \Om$.

\subsection{Cohomology space} \label{sec_cohom}

The space $\cohom$ of functions $\tv \in \HcD \cap \HdvN$ such that $\Crl \tv = \nv$ and
$\Dv \tv = 0$ is the ``cohomology'' space associated with the domain $\Om$ and the partition
of its boundary $\pt \Om = \overline \GD \cup \overline \GN$. When $\Om$ is simply connected
and $\GD$ is connected, this space is trivial; then the conditions associated with it below
can be disregarded. In the general case, $\cohom$ is finite-dimensional, and its dimension
depends on the topology of $\Om$ and $\GD$, see~\cite{Fer_Gil_Maxw_BC_97, Gross_Kot_elm_04}.

\subsection{Tetrahedral mesh, patches of elements, and the hat functions} \label{sec_mesh}

Let $\Th$ be a simplicial mesh of the domain $\Om$, \ie,
$\cup_{\elm \in \Th} \elm = \overline \Om$, where any element $\elm \in \Th$ is a
closed tetrahedron with nonzero measure, and where the intersection of two
different tetrahedra is either empty or their common vertex, edge, or face. The shape-regularity
parameter of the mesh $\Th$ is the positive real number
$\kappa_\Th \eq \max_{\elm \in \Th} h_\elm / \rho_\elm$, where $h_\elm$ is the diameter of the
tetrahedron $\elm$ and $\rho_\elm$ is the diameter of the largest ball contained in $\elm$.
These assumptions are standard and allow for strongly graded meshes with local refinements.
We will use the notation $a \ls b$ when there exists a positive constant $C$ only depending on
$\kappa_\Th$ such that $a \leq C b$.

We denote the set of vertices of the mesh $\Th$ by $\Vh$; it is composed of interior vertices
lying in $\Om$ and of vertices lying on the boundary $\pt \Om$.
For an element $\elm \in \Th$, $\FK$ denotes the set of its faces and $\VK$ the set of its vertices.
Conversely, for a vertex $\ver \in \Vh$, $\Ta$ denotes the patch of the elements of $\Th$ that share $\ver$, and $\oma$ is the corresponding open subdomain with diameter $h_\oma$.
A particular role below will be played by the continuous, piecewise affine ``hat'' function $\psia$ which takes value $1$ at the vertex $\ver$ and zero at the other vertices.
We note that $\oma$ corresponds to the support of $\psia$ and that the functions $\psia$ form the partition of unity
\be \label{eq_PU}
    \suma \psia = 1.
\ee
We will also need the patch $\tTa$ \cye{extended by one layer of neighbors} and the associated subdomain $\toma$, corresponding to the supports of the hat functions $\psib$ for all vertices $\vertt$ contained in the patch $\Ta$.

\subsection{Piecewise polynomial spaces}
\label{sec_disc_sp}

Let $q \geq 0$ be an integer. For a single tetrahedron $\elm \in \Th$, denote
by $\PP_q(\elm)$ the space of scalar-valued polynomials on $\elm$ of total degree at most $q$,
and by $[\PP_q(\elm)]^3$ the space of vector-valued polynomials on $\elm$ with each
component in $\PP_q(\elm)$. The N\'ed\'elec~\cite{Bof_Brez_For_MFEs_13,Ned_mix_R_3_80} space of
degree $q$ on $\elm$ is then given by
\be \label{eq_N_K}
    \ND_q(\elm) \eq [\PP_q(\elm)]^3 + \tx \vp [\PP_q(\elm)]^3.
\ee
Similarly, the Raviart--Thomas~\cite{Bof_Brez_For_MFEs_13,Ra_Tho_MFE_77} space of degree $q$
on $\elm$ is given by
\be \label{eq_RT_K}
    \RT_q(\elm) \eq [\PP_q(\elm)]^3 + \PP_q(\elm) \tx.
\ee
We note that~\eqref{eq_N_K} and~\eqref{eq_RT_K} are equivalent to the writing with a direct
sum and only homogeneous polynomials in the second terms.
The second term in~\eqref{eq_N_K} can also equivalently be given by homogeneous $(q+1)$-degree polynomials $\tv_h$ such that $\tx \scp \tv_h(\tx) = 0$ for all $\tx \in \elm$.

We will below extensively use the broken, piecewise polynomial spaces formed from the above
element spaces
\ban
    \PP_q(\Th) & \eq \{  v_h \in  \Lt; \,  v_h|_\elm \in \PP_q(\elm)\quad \forall \elm \in \Th\},\\
    \ND_q(\Th) & \eq \{\tv_h \in \tLt; \, \tv_h|_\elm \in \ND_q(\elm)\quad \forall \elm \in \Th\},\\
    \RT_q(\Th) & \eq \{\tv_h \in \tLt; \, \tv_h|_\elm \in \RT_q(\elm)\quad \forall \elm \in \Th\}.
\ean
To form the usual finite-dimensional Sobolev subspaces, we will write $\PP_q (\Th) \cap \Ho$
(for $q \geq 1$), $\ND_q(\Th) \cap \Hc$, $\RT_q(\Th) \cap \Hdv$ (both for $q \geq 0$), and
similarly for the subspaces reflecting the different boundary conditions. The same notation will
also be used on the patches \cye{of mesh elements} $\Ta$.

\subsection{\texorpdfstring{$\LT$}{L2}-orthogonal projectors and the Raviart--Thomas interpolator}\label{sec_projs}

For $q \geq 0$, let $\Pi_q$ denote the $\Lti{\elm}$-orthogonal projector onto
$\PP_q(\elm)$\cye{. Since this is an elementwise procedure, we keep the same notation for} the $\Lt$-orthogonal projector onto $\PP_q(\Th)$ \cye{given
for $v \in \Lt$ as $\Pi_q(v) \in \PP_q(\Th)$ such that $(\Pi_q(v),w_h) = (v,w_h)$ for all $w_h \in \PP_q(\Th)$}. Then, $\bPi_q$ is given componentwise by $\Pi_q$.

Let $\elm \in \Th$ be a mesh tetrahedron and \cye{let} $\tv \in [C^1(\elm)]^3$ be given.
Following~\cite{Bof_Brez_For_MFEs_13, Ra_Tho_MFE_77}, the canonical $q$-degree Raviart--Thomas
interpola\cye{nt} $\RTproj{q}(\tv) \in \RT_q(\elm)$, $q \geq 0$, is given by
\bse \label{eq_RT_proj} \bat{2}
\<\RTproj{q}(\tv) \scp \tn_{\elm}, r_h\>_\sd & = \<\tv \scp \tn_{\elm}, r_h\>_\sd \qquad & & \forall r_h \in \PP_q(\sd), \quad \forall \sd \in \FK, \label{eq_RT_proj_face} \\
(\RTproj{q}(\tv), \tr_h)_\elm & = (\tv, \tr_h)_\elm & & \forall \tr_h \in [\PP_{q-1}(\elm)]^3. \label{eq_RT_proj_vol}
\eat \ese
Less regular functions can be used in~\eqref{eq_RT_proj}, but $\tv \in [C^1(\elm)]^3$
will be sufficient for our purposes; we will actually only employ polynomial $\tv$ \cye{as arguments of $\RTproj{q}$}.
This interpolator crucially satisfies, on the tetrahedron $\elm$, the commuting property
\be \label{eq_com_prop}
    \Dv \RTproj{q}(\tv) = \cye{\Pi}_q (\Dv \tv) \qquad \forall \tv \in [C^1(\elm)]^3.
\ee

\subsection{\cye{Sobolev spaces on the patch subdomains \texorpdfstring{$\oma$}{}}} \label{sec_cont_sp_oma}

Let $\ver \in \Vh$ be an interior vertex. Then we set 1) $\Hsa \eq \{v \in \Hoi{\oma}; \, (v,1)_\oma = 0\}$, so
that $\Hsa$ is the subspace of those $\Hoi{\oma}$ functions whose mean value vanishes;
2) $\Hca \eq \{\tv \in \Hci{\oma}; \, \tv \vp \tn_{\oma}=\cye{\nv}$ on $\pt \oma\}$,
where the tangential trace is understood \cye{by duality} as above \cye{in Section~\ref{sec_cont_sp_BC}}; and, similarly, 3)
$\Hdva \eq \{\tv \in \Hdvi{\oma}; \, \tv \scp \tn_{\oma}=0$ on $\pt \oma\}$.
We will also need 4) $\Hcdual{\oma} \eq \Hci{\oma}$ \cye{(the symbol $\dagger$ is used here for notational purposes)}.
The situation is more subtle for boundary vertices. As a first possibility,
if $\ver \in \GN$ (\ie, $\ver \in \Vh$ is a boundary vertex such that all the faces sharing
the vertex $\ver$ lie in $\GN$), then the spaces $\Hsa$, $\Hca$, $\Hdva$, and $\Hcdual{\oma}$
are defined as above.
Secondly, when $\ver \in \overline \GD$, then at least one of the faces
sharing the vertex $\ver$ lies in $\overline \GD$\cye{; we denote by $\gD$ the subset of $\GD$ formed by all mesh faces sharing the vertex $\ver$ and lying in $\overline \GD$.} In this situation, we let 1) $\Hsa \eq \{v \in \Hoi{\oma}; \, v = 0$ on $\gD\}$; 2)
$\Hca \eq \{\tv \in \Hci{\oma}; \, \tv \vp \tn_{\oma}=\cye{\nv}$ on $\pt \oma \setminus \gD\}$; 3)
$\Hdva \eq \{\tv \in \Hdvi{\oma}; \, \tv \scp \tn_{\oma}=0$ on $\pt \oma \setminus \gD\}$; and 4)
$\Hcdual{\oma} \eq \{\tv \in \Hci{\oma}; \, \tv \vp \tn_{\oma}=\cye{\nv}$ on $\gD\}$.

\subsection{\cye{Functional inequalities}} \label{sec_ineq}

\cye{To work with data oscillation terms, we will employ the three following functional inequalities. First,} from~\cite[Theorems 3.4 and 3.5]{Cost_Dau_Nic_sing_Maxw_99},
\cite[Theorem~2.1]{Hipt_Pechs_discr_reg_dec_19}, and the discussion
in~\cite[Section~3.2.1]{Chaum_Ern_Voh_Maxw_22}, it follows that there exists a constant
$\Clift$ such that for all $\tv \in \HcD$, there exists $\tw \in \tHo \cye{\cap \HcD}$) such that $\Crl \tw = \Crl \tv$,
and
\be
\label{eq_estimate_lift}
\norm{\Gr \tw} \leq \Clift \norm{\Crl \tv}.
\ee
\cye{When either $\GD$ or $\GN$ has zero measure and if $\Om$ is convex, one can take $\Clift = 1$,
see~\cite{Cost_Dau_Nic_sing_Maxw_99} together with~\cite[Theorem~3.7]{Gir_Rav_NS_86}
for Dirichlet boundary conditions and~\cite[Theorem~3.9]{Gir_Rav_NS_86} for Neumann boundary
conditions.}

\cye{Second, for any mesh element $\elm \in \Th$ and $\tv \in \tHoi{\elm}$,} there hold\cye{s} the \cye{Poincar\'e inequality}
\be \label{eq_PF}
    \norm{\tv - \bPi_0(\tv)}_{\cye{\elm}} \leq \cye{\frac{h_\elm}{\pi}} \norm{\Gr \tv}_{\cye{\elm}} \cye{,}
\ee
\cye{where $\bPi_0(\tv)$ denotes the componentwise mean value of $\tv$ on $\elm$.}


\cye{Third,} the Poincar\'e--Friedrichs--Weber inequality, see~\cite[Proposition 7.4]{Fer_Gil_Maxw_BC_97}
and more precisely~\cite[Theorem A.1]{Chaum_Ern_Voh_Maxw_22} for the form of the constant,
will be useful: for all vertices $\ver \in \Vh$ and all vector-valued functions
$\tv \in \Hcdual{\oma} \cap \Hdva$ with $\Dv \tv =0$, we have
\be
    \label{eq_weber_patch}
    \norm{\tv}_{\oma} \lesssim h_{\oma} \norm{\Crl \tv}_{\oma}.
\ee
Strictly speaking, the inequality is established in~\cite[Theorem A.1]{Chaum_Ern_Voh_Maxw_22}
for edge patches, but the proof can be easily extended to vertex patches.

\section{\cye{Setting}} \label{sec_set}

\cye{The purpose of this section is to introduce the curl--curl problem and its N\'ed\'elec finite element approximation. We also identify, in a form of two self-standing assumptions, the kernel properties solely needed for our analysis.}

\subsection{Current density}
\label{eq_curr}

\cye{The following assumption is central for us:}

\bas[Current density $\tj$] \label{ass_jh}
Let $\tj$ be $\HdvN$-conforming, di\-ver\-gen\-ce-free, and $\Lt$-orthogonal
to the cohomology space $\cohom$, \ie,
\bse \label{eq_j} \ba
    \tj & \in \HdvN, \label{eq_j_Hdv}\\
    \Dv \tj & = 0, \label{eq_j_div}\\
    (\tj,\bvphi) & = 0 \quad \forall \bvphi \in \cohom. \label{eq_j_cohom}
\ea\ese
\eas

\cye{Let us recall from Section~\ref{sec_cohom} that when $\Om$ is simply connected
and $\GD$ is connected, then~\eqref{eq_j_cohom} can be disregarded.}
Sometimes, to illustrate the main ideas, we will additionally suppose that $\tj$
is a piecewise $p$-degree Raviart--Thomas polynomial, $\tj \in \RT_p(\Th) \cap \HdvN$.
Assumption~\ref{ass_jh} equivalently means that $\tj$ belongs to the
range of the curl operator, \ie, there exists $\tv \in \HcN$
such that $\Crl \tv = \tj$.

\subsection{The curl--curl problem}
\label{eq_Maxw}

The curl--curl problem we study here reads: find the magnetic vector potential $\tA: \Om \to \RR^3$ such that
\bse \label{eq_maxwell_strong}
\bat{3}
    \Crl (\Crl \tA) & = \tj, & \quad \Dv \tA & = 0 \qquad & & \text{in } \Om, \\
    \tA \vp \tn_{\Om} & = \nv, & & & &   \text{on } \GD,\\
    (\Crl \tA) \vp \tn_{\Om} & = \nv, & \quad \tA \scp \tn_{\Om} & = 0, &  & \text{on } \GN,
\eat
\ese
with the additional requirement that $(\tA,\bvphi) = 0$ for all $\bvphi$ from the cohomology
space $\cohom$ introduced in Section~\ref{sec_cohom} to ensure uniqueness.
Introducing $\tK(\Om) \eq \left \{\tv \in \HcD; \; \Crl \tv = \nv
\right \}$, the weak formulation of problem~\eqref{eq_maxwell_strong}, \cf, \eg,
\cite{Bof_Brez_For_MFEs_13}, consists in finding a pair
$(\tA,\tq) \in \HcD \times \tK(\Om)$ such that
\bse \label{eq_maxwell_weak}
\bat{2}
(\tA, \bvphi) &= 0 \qquad \qquad \quad && \forall \bvphi \in \tK(\Om)\cye{,}
\label{eq_maxwell_weak_1}\\
(\Crl \tA, \Crl \tv) + (\tq, \tv) &= (\tj, \tv) && \forall \tv \in \HcD. \label{eq_maxwell_weak_2}
\eat
\ese
Picking the test function $\tv = \tq$ in~\eqref{eq_maxwell_weak_2},
we see that $\tq = \nv$. Thus $\tA$ is such that
\bse \label{eq_maxwell_weak_II} \ba
& \tA \in \HcD \\
& (\Crl \tA, \Crl \tv) = (\tj, \tv) \qquad \forall \tv \in \HcD.
\ea\ese

\cye{\br[Characterization~\eqref{eq_maxwell_weak_II} of the magnetic vector potential $\tA$] All the main developments below actually rely solely on~\eqref{eq_maxwell_weak_II}, so that in particular the vector field $\tA$ can in our setting only be defined up to a curl-free component.
Remark that the existence of $\tA$ satisfying~\eqref{eq_maxwell_weak_II} is a direct consequence of Assumption~\ref{ass_jh}, and that a direct consequence of~\eqref{eq_maxwell_weak_II} is that $\Crl \tA \in \HcN$ with $\Crl (\Crl \tA) = \tj$. \er}

\subsection{N\'ed\'elec finite element approximation} \label{eq_appr}

For an integer $p \geq 0$ \cye{that we consider fixed henceforth}, let the N\'ed\'elec finite element space be given by
$\tV_h \eq \ND_p(\Th) \cap \HcD$. The subspace
$\tK_h \eq \{ \tv_h \in \tV_h; \; \Crl \tv_h = \nv \}$ is simply $\Gr (\PP_{p+1}(\Th) \cap \HoD)$ when $\Om$ is simply connected and $\GD$ is connected, and can be readily identified by introducing ``cuts'' in the mesh mimicking the construction
of the cohomology space $\cohom$, see~\cite[Chapter 6]{Gross_Kot_elm_04}.
The finite element approximation of~\eqref{eq_maxwell_weak} is a pair
$(\tAh, \tq_h) \in \tV_h \times \tK_h$ such that
\bse \label{eq_maxwell_discrete}
\bat{2}
(\tAh, \bvphi_h) &= 0 \qquad \qquad\quad && \forall \bvphi_h \in \tK_h,\\
(\Crl \tAh, \Crl \tv_h) + (\tq_h, \tv_h) &= (\tj, \tv_h) && \forall \tv_{h} \in \tV_h.
\eat\ese
Observing that $\tK_h \subset \tK$, this actually leads to $\tAh \in \tV_h$ such that
\be \label{eq_maxwell_discrete_II}
(\Crl \tAh, \Crl \tv_h) = (\tj, \tv_h) \qquad \forall \tv_h \in \tV_h.
\ee
\cye{In the developments below, we can actually still weaken~\eqref{eq_maxwell_discrete_II} and rely solely on:}

\bas[Discrete magnetic vector potential $\tAh$]
\label{ass_Ah}
Let $\tAh$ be a piecewise $p$-degree N\'ed\'elec polynomial satisfying a
lowest-order N\'ed\'elec \cye{orthogonality property},
\bse
\label{eq_Ah}
\ba
    & \tAh \in \ND_p(\Th) \cap \HcD, \label{eq_Ah_conf} \\
    & (\Crl \tAh,\Crl \tv_h) = (\tj,\tv_h) \label{eq_Ah_orth} \qquad
\forall \tv_h \in \ND_0(\Th) \cap \HcD.
\ea \ese
\eas

\section{Motivation} \label{sec_mot}

Let $\tj$ satisfy Assumption~\ref{ass_jh}. We motivate here our approach by showing
how an equilibrated flux $\fr$ may be constructed locally from any $\tA$
satisfying~\eqref{eq_maxwell_weak_II} at the continuous level. These observations are the basis of the
actual flux equilibration procedure involving $\tAh$ satisfying Assumption~\ref{ass_Ah}
at the discrete level that we develop in Sections~\ref{sec_div_free_dec} and~\ref{sec_equil}
below. We would in particular like to identify a patchwise construction such that
\bse \label{eq_fra_prop} \ba
    \fr^\ver & \in \Hca, \\
    \fr \eq \suma \fr^\ver & \in \HcN, \\
    \Crl \fr & = \tj.
\ea \ese
At the continuous level, the solution is trivially
\be \label{eq_frv}
    \fr^\ver = \psia (\Crl \tA).
\ee

We now rewrite the above definition implicitly. The idea is to introduce
\be \label{eq_fra}
    \fr^\ver \eq \arg \min_{\substack{\tv \in \Hca\\ \Crl \tv= \tj^\ver}} \norm{\tv - \psia(\Crl \tA)}_\oma^2 \qquad \forall \ver \in \Vh
\ee
with a suitable curl constraint $\tj^\ver$. Since
\be \label{eq_psia_E}
\Crl (\psia (\Crl \tA))
=
\psia \underbrace{(\Crl (\Crl \tA))}_{\tj}
+
\underbrace{\Gr \psia \vp (\Crl \tA)}_{\bm{\theta}^\ver}
\ee
we have
\be \label{eq_ja}
\tj^\ver \eq \psia \tj + \bm{\theta}^\ver,
\qquad
\bm{\theta}^\ver \eq \Gr \psia \vp (\Crl \tA).
\ee
Importantly, it holds that
\bse \label{eq_tsea} \ba
    \bm{\theta}^\ver & \in \Hdva, \label{eq_tsea_sp}  \\
    \Dv \bm{\theta}^\ver & = \underbrace{\Crl \Gr \psia}_{0} \scp (\Crl \tA) - \Gr \psia \scp \underbrace{\Crl (\Crl \tA)}_{\tj} = - \Gr \psia \scp \tj, \label{eq_tsea_dv}  \\
    \suma \bm{\theta}^\ver & = \suma \Gr \psia \vp (\Crl \tA) = \nv, \label{eq_tsea_sum}
\ea \ese
where the last property follows by the partition of unity~\eqref{eq_PU}. Consequently,
\bse \label{eq_je} \ba
    \tj^\ver & = \psia \tj + \bm{\theta}^\ver \in \Hdva, \label{eq_je_a}\\
    \Dv \tj^\ver & = \Gr \psia \scp \tj + \psia \underbrace{\Dv \tj}_{0} + \Dv \bm{\theta}^\ver = 0, \label{eq_je_div}\\[-0.2cm]
    \suma \tj^\ver & = \tj, \label{eq_je_sum}
\ea \ese
\cye{which gives a decomposition of the divergence-free current density $\tj$ into divergence-free contributions $\tj^\ver$ defined over the vertex patch subdomains $\oma$.}
The above auxiliary fields $\bm{\theta}^\ver$ can also be defined implicitly as the solution
to the minimization problems:
\be \label{eq_sa}
    \bm{\theta}^\ver \eq \arg \min_{\substack{\tv \in \Hdva\\ \Dv \tv= - \Gr \psia \scp \tj}} \norm{\tv - \Gr \psia \vp (\Crl \tA)}_\oma^2 \qquad \forall \ver \in \Vh.
\ee
We shall now mimic~\eqref{eq_fra}, \eqref{eq_je}, and~\eqref{eq_sa} at the
discrete level.

\section{\cye{Main results}} \label{sec_main_res}

\cye{In this section, we summarize our main results.}

\subsection{Stable divergence-free patchwise decomposition of the \cye{given current density \texorpdfstring{$\tj$}{}}}
\label{sec_div_free_dec}

The central issue for our approach is a stable divergence-free patchwise decomposition of the current density $\tj$ in the spirit of~\eqref{eq_je}. \cye{For this purpose, we first} design an appropriate discrete variant of~\eqref{eq_sa}, where we crucially rely on the patchwise orthogonality stemming from Assumption~\ref{ass_Ah}.
We will initially be requested to work with the increased polynomial degree \ifSIAM$\pp \eq \min\{p,1\}$, \else
\be \label{eq_pp}
    \pp \eq \min\{p,1\},
\ee
\fi
recalling that $p \geq 0$ is fixed in Section~\ref{eq_appr}. \cye{We start with:}

\bd[Patchwise contributions $\tj_h^\ver$] \label{def_div_free_dec}
Let $\tj$ and $\tAh$ satisfy respectively Assumptions~\ref{ass_jh} and~\ref{ass_Ah}.
Carry out the three following steps:

\begin{enumerate} [parsep=0.2pt, itemsep=0.2pt, topsep=1pt, partopsep=1pt, leftmargin=13pt]

\item \label{it_1} For all vertices $\ver \in \Vh$, consider the $\pp$-degree Raviart--Thomas patchwise \cye{\cye{(seemingly} over-constrained\cye{)}} minimizations
\be \label{eq_tha}
\tha \eq \arg \hspace*{-1cm}
\min_{\substack{
\tv_h \in \RT_{\pp}(\Ta) \cap \Hdva
\\
\Dv \tv_h= \Pi_\pp(- \Gr \psia \scp \tj)
\\
(\tv_h, \tr_h)_\elm = (\Gr \psia \vp (\Crl \tAh), \tr_h)_\elm
\quad
\forall \tr_h \in [\PP_0(\elm)]^{3}, \, \forall \elm \in \Ta}} \hspace*{-1cm} \norm{\tv_h - \Gr \psia \vp (\Crl \tAh)}_\oma^2,
\ee
where \cye{in addition to the usual} normal trace and divergence\cye{, the constraints} additionally also \cye{concern} elementwise product with piecewise vector-valued constants.

\item \label{it_2} Extending $\tha$ by zero outside of the patch subdomains $\oma$, set
\be \label{eq_td}
    \bdl_h \eq \suma \tha.
\ee
For all tetrahedra $\elm \in \Th$, consider the $(p+1)$-degree Raviart--Thomas elementwise
minimizations:
\bse \label{eq_tda} \bat{3}
    \tda|_\elm & \! \eq \! \arg \hspace*{-0.7cm} \min_{\substack{\tv_h \in \RT_{1}(\elm)\\ \Dv \tv_h= 0\\ \tv_h \scp \tn_{\elm} = \RTproj{1}(\cye{(}\psia \bdl_h\cye{)|_\elm}) \scp \tn_{\elm} \, \text{ on } \pt \elm }} \hspace*{-0.7cm} \norm{\tv_h - \RTproj{1}(\cye{(}\psia \bdl_h\cye{)|_\elm})}_\elm^2 \quad & & \forall \ver \in \VK \, & & \text{ when } p=0, \label{eq_tda_0} \\
    \tda|_\elm & \! \eq \! \arg \hspace*{-0.4cm} \min_{\substack{\tv_h \in \RT_{p+1}(\elm)\\ \Dv \tv_h= 0\\ \tv_h \scp \tn_{\elm} = \psia \bdl_h \scp \tn_{\elm} \, \text{ on } \pt \elm }} \hspace*{-0.35cm} \norm{\tv_h - \psia \bdl_h}_\elm^2 & & \forall \ver \in \VK & & \text{ when } p \geq 1, \label{eq_tda_p}
\eat\ese
\cye{which yields the divergence-free decomposition
\[
    \bdl_h = \suma \tda.
\]}

\item \label{it_3} For all vertices $\ver \in \Vh$, define
\be \label{eq_tjha}
    \tj_h^\ver \eq \psia \tj + \tha - \tda.
\ee
\end{enumerate}
\ed

\cye{For a vertex $\ver \in \Vh$ and the extended \cye{(second-order)} patch $\tTa$, define the data oscillation term}
\be \label{eq_osc_j}
    \widetilde \eta_{\mathrm{osc},\tj}^\ver  \eq  \Bigg\{\sum_{\elm \in \tTa} \Big(\frac{h_\elm}{\pi} \norm{\tj - \bPi_\pp(\tj)}_\elm\Big)^2\Bigg\}^\ft.
\ee
\cye{We crucially have:}

\bt[\cye{Stable} divergence-free patchwise decomposition of $\tj$]
\label{thm_div_free_dec}
Let $\tj$ and $\tAh$ satisfy respectively Assumptions~\ref{ass_jh} and~\ref{ass_Ah}.
Let $\tj_h^\ver$ be given by Definition~\ref{def_div_free_dec} for all vertices $\ver \in \Vh$.
Then
\bse \label{eq_div_free_dec} \ba
    \tj_h^\ver & \in \Hdva, \label{eq_div_free_dec_contr} \\
    \Dv \tj_h^\ver & = \Gr \psia \scp (\tj - \bPi_\pp(\tj)), \label{eq_div_free_dec_dv} \\
    \suma \tj_h^\ver & = \tj, \label{eq_div_free_dec_dec}
\ea \ese
where the extension of $\tj_h^\ver$ by zero outside of the patch subdomains $\oma$ is
understood in the last two properties. Moreover, when $\tj \in \RT_p(\Th) \cap \HdvN$
is piecewise polynomial, \cye{then}\cye{, in strengthening of~\eqref{eq_div_free_dec_contr}--\eqref{eq_div_free_dec_dv},}
\bse \label{eq_div_free_dec_pol} \ba
    \tj_h^\ver & \in \RT_{p+1}(\Ta) \cap \Hdva, \label{eq_div_free_dec_contr_pol} \\
    \Dv \tj_h^\ver & = 0. \label{eq_div_free_dec_dv_pol}
\ea \ese

Let \cye{in addition} $\tA$ satisfying~\eqref{eq_maxwell_weak_II} be arbitrary
and let, as in~\eqref{eq_ja},
\[
    \tj^\ver \eq \psia \tj + \Gr \psia \vp (\Crl \tA).
\]
Then
\be \label{eq_dec_stab}
\norm{\tj^\ver - \tj_h^\ver}_\oma
\ls
h_\oma^{-1} \big[ \norm{\Crl (\tA - \tAh)}_\toma
+
\widetilde \eta_{\mathrm{osc},\tj}^\ver \big].
\ee
\et

\subsection*{Remarks}

Several remarks \cye{about Definition~\ref{def_div_free_dec} and Theorem~\ref{thm_div_free_dec}} are in order:

\begin{enumerate} [parsep=0.2pt, itemsep=0.2pt, topsep=1pt, partopsep=1pt, leftmargin=13pt]

\item At the discrete level, \cye{$\Gr \psia \vp (\Crl \tAh) \not \in \Hdva$, in contrast to $\Gr \psia \vp (\Crl \tA)$, see~\eqref{eq_ja}--\eqref{eq_tsea}}. The auxiliary field $\tha$ from~\eqref{eq_tha} \cye{is the projection of $\Gr \psia \vp (\Crl \tAh)$ to} $\RT_{\pp}(\Ta) \cap \Hdva$ satisfying $\Dv \tha= \Pi_\pp(- \Gr \psia \scp \tj)$. Step~\ref{it_1} of Definition~\ref{def_div_free_dec} thus mimics~\eqref{eq_sa} and achieves equivalents to~\eqref{eq_tsea_sp} and~\eqref{eq_tsea_dv}.
    Unfortunately, $\cye{\bdl_h = }\suma \tha$ \cye{given by~\eqref{eq_td}} \cye{typically does not equal $\nv$,} which would mimic~\eqref{eq_tsea_sum}.


\item
Step~\ref{it_2} of Definition~\ref{def_div_free_dec}  \cye{yields the corrected fields $\tha - \tda$ which mimic~\eqref{eq_tsea} entirely in that (see Lemma~\ref{lem_tdh} below for details)}
\ban
\tha - \tda & \in \RT_{p+1}(\Ta) \cap \Hdva, \\
\Dv(\tha - \tda) & = \Pi_\pp(- \Gr \psia \scp \tj), \\
\suma(\tha - \tda) & = \bdl_h - \bdl_h = \nv.
\ean

\item
Step~\ref{it_3} of Definition~\ref{def_div_free_dec}\cye{, in view of Theorem~\ref{thm_div_free_dec},} finally materializes~\eqref{eq_je}
at the discrete level.

\item
Property~\eqref{eq_dec_stab} from Theorem~\ref{thm_div_free_dec} shows that the
local discrete decomposition~\eqref{eq_tjha} compares in a $p$-robust way to the
continuous-level decomposition~\eqref{eq_je}, up to data oscillation \cye{given by~\eqref{eq_osc_j}}.


\item
\cye{Minimization}~\eqref{eq_tha} contains an \cye{additional} constraint on
the elementwise product with piecewise vector-valued constants. \cye{Existence, uniqueness, and $p$-robust stability theory for such problems is developed in Appendix~\ref{app_over_constr_min}. Section~\ref{sec_det_proofs} shows that this applies to our setting under the orthogonality in Assumption~\ref{ass_Ah}}.

\item \cye{The additional constraint in~\eqref{eq_tha} also ensures the existence, uniqueness, and $p$-robust stability of the elementwise problems~\eqref{eq_tda}, where $\tda$ form a divergence-free local decomposition of $\bdl_h$ following Appendix~\ref{app_div_free_dec}, see Lemma~\ref{lem_tdh} below.}





\end{enumerate}

\subsection{Equilibrated flux reconstruction based on local patchwise minimizations in \texorpdfstring{$\HC$}{H(curl)}} \label{sec_equil}

We now identify an appropriate discrete variant of~\eqref{eq_fra_prop}--\eqref{eq_fra}, giving a locally computable equilibrated flux $\frh$. \cye{Let \cye{$\tj_h^\ver$ be given by Definition~\ref{def_div_free_dec} and set}
\be \label{eq_bjha}
\bar \tj_h^\ver \eq \arg \hspace*{-1cm}
\min_{\substack{
\tv_h \in \RT_{p+1}(\Ta) \cap \Hdva
\\
\Dv \tv_h = 0
\\
(\tv_h, \tr_h)_\elm = (\tj_h^\ver, \tr_h)_\elm
\quad
\forall \tr_h \in [\PP_0(\elm)]^{3}, \, \forall \elm \in \Ta}} \hspace*{-1cm} \norm{\tv_h - \tj_h^\ver}_\oma^2.
\ee
When $\tj$ is piecewise polynomial, $\tj \in \RT_p(\Th) \cap \HdvN$, then~\eqref{eq_div_free_dec_pol} implies that $\bar \tj_h^\ver = \tj_h^\ver$, so that there is no need for~\eqref{eq_bjha}. In general, \cye{the role of~\eqref{eq_bjha} is to prepare a piecewise polynomial datum for a discrete version of~\eqref{eq_fra}: it} projects the non-polynomial and non-divergence-free $\tj_h^\ver$ to $\bar \tj_h^\ver \in \RT_{p+1}(\Ta) \cap \Hdva$ with $\Dv \bar \tj_h^\ver = 0$. \cye{Problem~\eqref{eq_bjha}} has the same form as problem~\eqref{eq_tha} and is well-posed following Appendix~\ref{app_over_constr_min} below. \cye{With the Raviart--Thomas divergence-free $\bar \tj_h^\ver$,} the following \cye{N\'ed\'elec} local equilibration problem is well posed by standard arguments, see, \eg, \cite{Bof_Brez_For_MFEs_13}:}

\bd[Equilibrated flux reconstruction based on local minimization in $\HC$] \label{def_EFR}
Let $\tj$ and $\tAh$ satisfy respectively Assumptions~\ref{ass_jh} and~\ref{ass_Ah} and let,
for all vertices $\ver \in \Vh$, $\tj_h^\ver$ be given by Definition~\ref{def_div_free_dec} \cye{and $\bar \tj_h^\ver$ by~\eqref{eq_bjha}}.
Consider the patchwise minimizations
\bse\be \label{eq_frha}
\frh^\ver
\eq
\arg \min_{\substack{
\tv_h \in \ND_{p+1}(\Ta) \cap \Hca
\\
\Crl \tv_h= \cye{\bar \tj_h^\ver}
}}
\norm{\tv_h - \psia(\Crl \tAh)}_\oma^2.
\ee
Extending $\frh^\ver$ by zero outside of $\oma$, define
\be \label{eq_frh}
    \frh \eq \suma \frh^\ver.
\ee \ese
\ed

\cye{Recall} $\widetilde \eta_{\mathrm{osc},\tj}^\ver$ \cye{from~\eqref{eq_osc_j}} and define
\be \label{eq_osc_jha}
\eta_{\mathrm{osc},\tj_h^\ver}^\ver
\eq
h_\oma \norm{\cye{\bar \tj_h^\ver} - \tj_h^\ver}_\oma.
\ee
\cye{Crucially, t}he construction of Definition~\ref{def_EFR} is a stable equilibration:

\bt[Equilibrium property and $p$-robust stability of the flux reconstruction] \label{thm_EFR}
Let $\tj$ and $\tAh$ satisfy respectively Assumptions~\ref{ass_jh} and~\ref{ass_Ah}. Then the equilibrated flux reconstruction $\frh$ from Definition~\ref{def_EFR} satisfies
\bse \label{eq_frh_prop} \ba
    \frh & \in \ND_{p+1}(\Th) \cap \HcN, \label{eq_frh_prop_sum}
\\
    \Crl \frh & = \tj \quad\text{ when } \tj \in \RT_p(\Th) \cap \HdvN. \label{eq_frh_prop_curl}
\ea \ese

Let \cye{in addition} $\tA$ satisfying~\eqref{eq_maxwell_weak_II} be arbitrary.
Then
\[
\norm{\frh^\ver  - \psia(\Crl \tAh)}_\oma
\ls
\norm{\Crl (\tA - \tAh)}_\toma
+
\eta_{\mathrm{osc},\tj_h^\ver}^\ver
+
\widetilde \eta_{\mathrm{osc},\tj}^\ver.
\]
\et

\subsection{Guaranteed, fully computable, constant-free, and $p$-robust a posteriori error estimates for the curl--curl problem}
\label{sec_a_post}

We apply here the results \cye{of Sections~\ref{sec_div_free_dec} and~\ref{sec_equil}} to a posteriori error analysis of the curl--curl problem~\eqref{eq_maxwell_strong}.

\bt[Guaranteed, fully computable, and constant-free upper bound]\label{thm_rel_eff}
Let $\tj$ satisfy Assumption~\ref{ass_jh}, let $\tA$ be the weak solution to
the curl--curl problem given by~\eqref{eq_maxwell_weak}, and let $\tAh$ be its
N\'ed\'elec finite element approximation given by~\eqref{eq_maxwell_discrete}.
Let $\tj_h^\ver$ be given by Definition~\ref{def_div_free_dec} for all vertices $\ver \in \Vh$, and let $\frh$ be given by
Definition~\ref{def_EFR}. Then
\[
\norm{\Crl (\tA - \tAh)}
\leq
\eta_{\mathrm{tot}}
\eq
\underbrace{\norm{\frh  - \Crl \tAh}}_{\eta}
+
\cye{\Clift} \underbrace{\Bigg\{\sum_{\elm \in \Th} \underbrace{\frac{h_\elm^2}{\pi^2} \norm{\tj - \Crl \frh}_\elm^2}_{\cye{\eta_{\mathrm{osc},\elm}^2}} \Bigg\}^\ft }_{\eta_{\mathrm{osc}}}
\]
and
\[
\norm{\frh  - \Crl \tAh}_\elm \cye{+\eta_{\mathrm{osc},\elm}}
\ls
\sum_{\ver \in \VK} \big[ \norm{\Crl (\tA - \tAh)}_\toma + \eta_{\mathrm{osc},\tj_h^\ver}^\ver + \widetilde \eta_{\mathrm{osc},\tj}^\ver \big].
\]
\et

\subsection*{Remarks} 

Several remarks are in order:

\begin{enumerate} [parsep=0.2pt, itemsep=0.2pt, topsep=1pt, partopsep=1pt, leftmargin=13pt]


\item \cye{On the discrete level, $\psia(\Crl \tAh) \not \in \Hca$, in contrast to $\psia (\Crl \tA)$ on the continuous level, see~\eqref{eq_fra_prop}--\eqref{eq_frv}. The equilibrated flux contribution $\frh^\ver$ from~\eqref{eq_frha} is its constrained projection to $\ND_{p+1}(\Ta) \cap \Hca$.
It mimics~\eqref{eq_fra} at the discrete level.}

\item \label{it_projs} When $\tj \in \RT_p(\Th) \cap \HdvN$\cye{, all the data oscillation estimators in Theorem~\ref{thm_rel_eff} vanish. Indeed, \eqref{eq_frh_prop_curl} implies that $\eta_{\mathrm{osc}} =0$, whereas $\bar \tj_h^\ver = \tj_h^\ver$} gives $\eta_{\mathrm{osc},\tj_h^\ver}^\ver = 0$, see~\eqref{eq_osc_jha}. Similarly, $\widetilde \eta_{\mathrm{osc},\tj}^\ver$ from~\eqref{eq_osc_j} vanishes as well (this is actually true up to $\tj \in \RT_\pp(\Th) \cap \HdvN$, since $\Dv \tj = 0$). Moreover, all these terms are higher-order with respect to $\norm{\Crl (\tA - \tAh)}$ if $\tj$ is piecewise smooth. \cye{We also note that} using~\eqref{eq_frh}, \eqref{eq_frha}, and~\eqref{eq_div_free_dec_dec}, the data oscillation term $\eta_{\mathrm{osc},\elm}$ can equivalently be rewritten with
\be \label{eq_osc_eq}
    \tj - \Crl \frh = \suma (\tj_h^\ver - \cye{\bar \tj_h^\ver}).
\ee
%



\item The equilibration of Definition~\ref{def_EFR} is performed in local N\'ed\'elec spaces of order $p+1$. This is in agreement with $p$-robust flux equilibrations from~\cite{Brae_Pill_Sch_p_rob_09, Ern_Voh_p_rob_15, Ern_Voh_p_rob_3D_20, Chaum_Voh_p_rob_3D_H_curl_22}. Similarly to~\cite{Brae_Pill_Sch_p_rob_09, Ern_Voh_adpt_IN_13}, it is also possible to design a downgrade of the orders of the local problems~\eqref{eq_frha} from $p + 1$ to $p$.

    \qquad Let us first discuss the case $p \geq 1$.
    The first step is to replace~\eqref{eq_tda_p} by~\eqref{eq_tda_0} with $\RT_1$ replaced by $\RT_p$. Then, according to Theorem~\ref{thm_div_free_dec_abs} with $\qq = p$, we obtain $\tda \in \RT_{p}(\Ta) \cap \Hdva$ in place of~\eqref{eq_tda_prop} \cye{below}.
    Second, we employ (elementwise) \cye{the projector} $\RTproj{p}(\psia \tj)$ in~\eqref{eq_tjha} \cye{in place of $\psia \tj$}. \cye{Let $\tj \in \RT_p(\Th) \cap \HdvN$. Then} $\suma \tj_h^\ver = \tj$ and $\Dv \tj_h^\ver = 0$ \cye{as in}~\eqref{eq_div_free_dec_dec}, \eqref{eq_div_free_dec_dv_pol}\cye{, but} $\tj_h^\ver \in \RT_p(\Ta) \cap \Hdva$ in place of~\eqref{eq_div_free_dec_contr_pol}. Consequently, \eqref{eq_frha} can be brought down to
    \be \label{eq_frha_p}
        \frh^\ver \eq \arg \min_{\substack{\tv_h \in \ND_p(\Ta) \cap \Hca\\ \Crl \tv_h= \cye{\tj_h^\ver}}} \norm{\tv_h - \NDproj{p} ( \psia(\Crl \tAh))}_\oma^2,
    \ee
    where $\NDproj{p}$ is the elementwise canonical $p$-degree N\'ed\'elec interpolate, analogue to~\eqref{eq_RT_proj}.
    This leads to a cheaper procedure where the guaranteed estimate of Theorem~\ref{thm_rel_eff} \cye{(with $\eta_{\mathrm{osc}} =0$)} still holds true \cye{when $\tj \in \RT_p(\Th) \cap \HdvN$; general $\tj$ can be covered by data oscillation terms.}
    Similarly, the local efficiency of Theorem~\ref{thm_rel_eff} is also preserved, with, however, the $p$-robustness theoretically lost. In particular, from~\eqref{eq_tsat_stab_2_A}, estimate~\eqref{eq_tsat_stab} below still holds true, up to a possibly $p$-dependent constant.

    \qquad Alternatively, for $p=0$ in particular, because of $\pp=p+1$ employed in~\eqref{eq_tha}, we need to replace~\eqref{eq_tjha} by
    \[
        \tj_h^\ver\cye{|_\elm} \eq \RTproj{0}(\cye{(}\psia \tj + \tha - \tda\cye{)|_\elm}) \qquad \cye{\forall \elm \in \Th}.
    \]
    Let $\tj \in \RT_0(\Th) \cap \HdvN$. Then clearly $\tj_h^\ver \in \RT_0(\Ta) \cap \Hdva$. Moreover, from~\eqref{eq_com_prop}, $\Dv \RTproj{0}(\tda\cye{|_\elm}) = \PP_0 (\Dv \cye{(}\tda\cye{)|_\elm}) = 0$, whereas $\Dv \RTproj{0}(\cye{(}\psia \tj\cye{)|_\elm}) = \PP_0 (\Dv (\psia \tj)\cye{|_\elm}) = \cye{(}\Gr \psia \scp \tj\cye{)|_\elm}$, also using that $\tj \in [\PP_0(\Ta)]^3$ as above in point~\ref{it_projs}, and similarly $\Dv \RTproj{0}(\tha\cye{|_\elm}) = \PP_0 (\Dv \cye{(}\tha\cye{)|_\elm}) = \cye{(}- \Gr \psia \scp \tj\cye{)|_\elm}$. Thus $\Dv \tj_h^\ver = 0$. Finally, $\cye{\big(}\suma \tj_h^\ver\cye{\big)|_\elm} = \RTproj{0}(\suma(\psia \tj + \tha - \tda)\cye{|_\elm}) = \RTproj{0}(\tj\cye{|_\elm}) = \tj\cye{|_\elm}$ by the linearity of the Raviart--Thomas projector $\RT_0$. Then the above discussion for $p \geq 1$ applies.

\item The approach of~\cite{Ged_Gee_Per_a_post_Maxw_20,Ged_Gee_Per_Sch_post_Maxw_21} includes solutions of local, a priori over-determined, problems on vertex patches in a multi-stage procedure. The present \cye{(again a priori} over-constrained\cye{)} problems~\eqref{eq_tha} and consecutive steps in Definitions~\ref{def_div_free_dec} and~\ref{def_EFR} share this spirit, though the minimizations directly determine the best-possible local energy error estimator contributions.



\end{enumerate}

\section{Numerical illustration}
\label{sec_num}

This section presents some numerical examples illustrating the key features
of the estimator of Theorem~\ref{thm_rel_eff}. We impose the Dirichlet boundary
condition on the whole boundary, \ie, $\GD \eq \pt \Om$. We consider both
structured meshes and unstructured meshes.
When we speak about a ``structured'' mesh, we
mean a Cartesian partition of $\Om$ into $N \times N \times N$ cubes where each
cube is first subdivided into 6 pyramids (with the basis a face and the apex the barycenter
of the cube) and then each pyramid into 4 tetrahedra. The corresponding mesh size is
$h = \sqrt{3}/(2N)$.
\tcf{On the other hand, the ``unstructured'' meshes are generated with the software pacakge
{\tt MMG3D}~\cite{Dobr_mmg3d}, where we simply require a maximum element size. These are typically
quasi-uniform, but do not have any particular reapeting structure (every vertex patch is different).
For both types of meshes, w}e consider the N\'ed\'elec finite element
approximation~\eqref{eq_maxwell_discrete} with varying degree $p \geq 1$.

\subsection{\texorpdfstring{$\tH^3(\Om)$ s}{S}olution with a polynomial right-hand side}
\label{sec_num_1}

We first consider the unit cube $\Om \eq (0,1)^3$ and a polynomial right-hand side
$\tj \eq (0,0,1)$, so that the data oscillation estimator $\eta_{\mathrm{osc}}$ vanishes.
We can show that the solution is given by $\tA = (0,0,A_3)$ with
\be
\label{eq_solution_1}
A_3(\tx)
\eq
\frac{16}{\pi^4}
\sum_{n,m \geq 1}
\frac{1}{nm(n^2+m^2)} \sin(n\pi\tx_1)\sin(m\pi\tx_2).
\ee
This function belongs to $H^3(\Om)$ but not to $H^4(\Om)$.
In practice, we cut the series at $n=m=100$, and obtain $\Crl \tA$
by analytically differentiating~\eqref{eq_solution_1}.

We first fix the polynomial degree and consider a sequence of meshes.
We use $p=1$ and structured meshes with $N=1,2,4,8$, and then $p=2$ and a
sequence of unstructured meshes. Figure~\ref{figure_hconv_pol} presents the
corresponding errors, estimates, and effectivity indices. We observe the expected
convergence rate $h^2$ (recall that $\tA \in \tH^{3}(\Om)$ merely). The estimator
$\eta = \eta_{\mathrm{tot}}$ \cye{($\eta_{\mathrm{osc}}=0$ here)} closely follows the error $\norm{\Crl(\tA-\tAh)}$,
and the effectivity index \cye{given by the ratio $\eta / \norm{\Crl (\tA - \tAh)}$} is \cye{above but} close to the optimal value $1$; we actually numerically
observe asymptotic exactness.

We then fix a mesh and increase the polynomial degree $p$ from $1$ to $6$.
We consider two configurations: a structured mesh where the unit cube is
split into $24$ tetrahedra as described above and an unstructured mesh consisting
of $\cye{20}$ tetrahedra. Figure~\ref{figure_pconv_pol} reports the results. The convergence
is not exponential, which is expected because of the solution's finite regularity.
Also in this setting, the estimator closely follows the actual error, and the effectivity
index always remains \cye{above but} close to $1$. In particular, the effectivity index does not increase
with $p$, which illustrates the $p$-robustness of the estimator.

Although this is not reported in the figures, we also numerically check that the
reconstructed flux $\frh$ is indeed equilibrated, \ie, $\norm{\tj-\Crl\frh} = 0$.
Because of finite precision arithmetics, this value is not exactly zero, but ranges
between $10^{-15}$ and $10^{-11}$, which is perfectly reasonable compared to
the actual error levels.

\begin{figure}
\begin{minipage}{.45\linewidth}
\begin{tikzpicture}
\begin{axis}%
[
	width=0.96\linewidth,
	xlabel={$h$ ($p=1$\cye{, structured meshes})},
	ylabel={Error and estimate},
	ymode = log,
	xmode = log,
	x dir = reverse
]

\addplot[thick,solid,black,mark=*, mark options={solid}]
table[x expr=sqrt(3)/(2*\thisrow{N}),y=err]{figures/data/hcurve1_pol.txt}
node[pos=.27,pin={[pin distance=1.0cm, pin edge=solid]-90:{$\norm{\Crl (\tA - \tAh)}$}}] {};
\addplot[thick,red,dashed,mark=square,mark options={solid}]
table[x expr=sqrt(3)/(2*\thisrow{N}),y=est_err] {figures/data/hcurve1_pol.txt}
node[pos=.8,pin={[pin distance=1.0cm, pin edge=solid]90:{$\eta$}}] {};

\SlopeTriangle{.65}{-.1}{.15}{-2}{$h^2$}{}

\end{axis}
\end{tikzpicture}
\end{minipage}
\qquad
\begin{minipage}{.45\linewidth}
\begin{tikzpicture}
\begin{axis}%
[
	width=0.96\linewidth,
	xlabel={\hspace{-0.35cm}$h$ ($p=2$\cye{, unstructured meshes})},
	ylabel={Error and estimate},
	ymode = log,
	xmode = log,
	x dir = reverse
]

\addplot[thick,solid,black,mark=*, mark options={solid}] table[x=h,y=err]     {figures/data/hcurve2_uns_pol.txt}
node[pos=.3,pin={[pin distance=1.0cm, pin edge=solid]-90:{$\norm{\Crl (\tA - \tAh)}$}}] {};
\addplot[thick,red,dashed,mark=square,mark options={solid}] table[x=h,y=est_err] {figures/data/hcurve2_uns_pol.txt}
node[pos=.8,pin={[pin distance=1.0cm, pin edge=solid]90:{$\eta$}}] {};

\SlopeTriangle{.65}{-.1}{.15}{-2}{$h^{2}$}{}

\end{axis}
\end{tikzpicture}
\end{minipage}

\medskip

\begin{minipage}{.45\linewidth}
\begin{tikzpicture}
\begin{axis}%
[
	width=0.96\linewidth,
	xlabel={$h$ ($p=1$\cye{, structured meshes})},
	ylabel={Effectivity index},
	xmode = log,
	ymin  = 1.04,
	ymax  = 1.062,
	x dir = reverse
]

\addplot[thick,red,dashed,mark=square,mark options={solid}]
table[x expr=sqrt(3)/(2*\thisrow{N}),y expr=\thisrow{est_err}/\thisrow{err}]{figures/data/hcurve1_pol.txt}
node[pos=.5,pin={[pin distance=0.5cm, pin edge=solid]-90:{$\frac{\eta}{\norm{\Crl (\tA - \tAh)}}$}}] {};

\end{axis}
\end{tikzpicture}
\end{minipage}
\qquad
\begin{minipage}{.45\linewidth}
\begin{tikzpicture}
\begin{axis}%
[
	width=0.96\linewidth,
	xlabel={\hspace{-0.35cm}$h$ ($p=2$\cye{, unstructured meshes})},
	ylabel={Effectivity index},
	xmode = log,
	ymin  = 1.047,
	ymax  = 1.08,
	x dir = reverse
]

\addplot[thick,red,dashed,mark=square,mark options={solid}] table[x=h,y expr=\thisrow{est_err}/\thisrow{err}]
{figures/data/hcurve2_uns_pol.txt}
node[pos=.5,pin={[pin distance=0.5cm, pin edge=solid]-90:{$\frac{\eta}{\norm{\Crl (\tA - \tAh)}}$}}] {};

\end{axis}
\end{tikzpicture}
\end{minipage}
\caption{[Smooth solution \cye{with limited regularity}~\eqref{eq_solution_1}] Uniform mesh refinement.}
\label{figure_hconv_pol}
\end{figure}
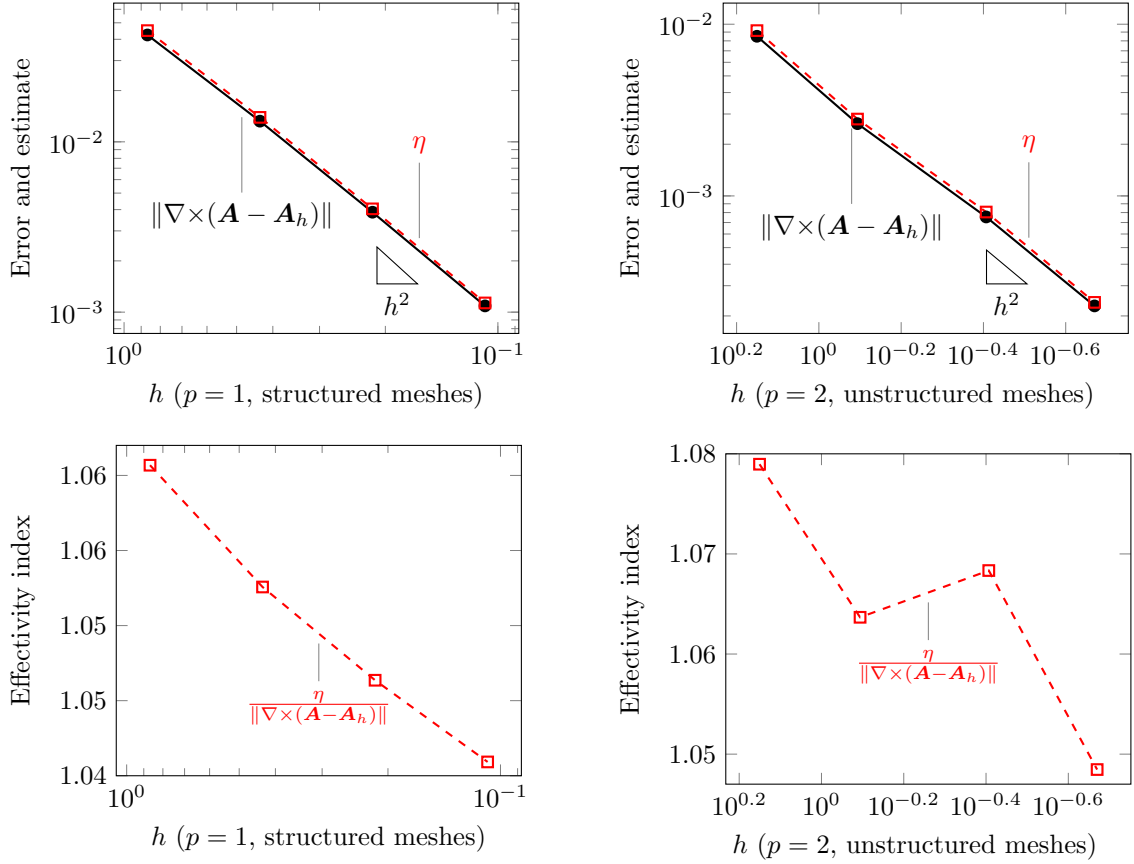

\begin{figure}
\begin{minipage}{.45\linewidth}
\begin{tikzpicture}
\begin{axis}%
[
	width=0.96\linewidth,
	xlabel={$p$ (structured mesh)},
	ylabel={Error and estimate},
	ymode = log,
	xtick={1,2,3,4,5,6}
]

\addplot[thick,solid,black,mark=*, mark options={solid}] table[x=p,y=err]     {figures/data/pcurve_pol.txt}
node[pos=.33,pin={[pin distance=1.2cm, pin edge=solid]-90:{$\norm{\Crl (\tA - \tAh)}$}}] {};
\addplot[thick,red,dashed,mark=square,mark options={solid}] table[x=p,y=est_err] {figures/data/pcurve_pol.txt}
node[pos=.8,pin={[pin distance=1cm, pin edge=solid]90:{$\eta$}}] {};

\end{axis}
\end{tikzpicture}
\end{minipage}
\qquad
\begin{minipage}{.45\linewidth}
\begin{tikzpicture}
\begin{axis}%
[
	width=0.96\linewidth,
	xlabel={$p$ (unstructured mesh)},
	ylabel={Error and estimate},
	ymode = log,
	xtick={1,2,3,4,5,6}
]

\addplot[thick,solid,black,mark=*, mark options={solid}] table[x=p,y=err]     {figures/data/pcurve_uns_pol.txt}
node[pos=.33,pin={[pin distance=1.2cm, pin edge=solid]-90:{$\norm{\Crl (\tA - \tAh)}$}}] {};
\addplot[thick,red,dashed,mark=square,mark options={solid}] table[x=p,y=est_err] {figures/data/pcurve_uns_pol.txt}
node[pos=.8,pin={[pin distance=1cm, pin edge=solid]90:{$\eta$}}] {};

\end{axis}
\end{tikzpicture}
\end{minipage}

\medskip

\begin{minipage}{.45\linewidth}
\begin{tikzpicture}
\begin{axis}%
[
	width=0.96\linewidth,
	xlabel={$p$ (structured mesh)},
	ylabel={Effectivity index},
	ymin=1.05,
	ymax=1.20,
	xtick={1,2,3,4,5,6}
]

\addplot[thick,red,dashed,mark=square,mark options={solid}] table[x=p,y expr=\thisrow{est_err}/\thisrow{err}]
{figures/data/pcurve_pol.txt}
node[pos=.555,pin={[pin distance=.52cm, pin edge=solid]-90:{$\frac{\eta}{\norm{\Crl (\tA - \tAh)}}$}}] {};

\end{axis}
\end{tikzpicture}
\end{minipage}
\qquad
\begin{minipage}{.45\linewidth}
\begin{tikzpicture}
\begin{axis}%
[
	width=0.96\linewidth,
	xlabel={$p$ (unstructured mesh)},
	ylabel={Effectivity index},
	ymin=1.05,
	ymax=1.20,
	xtick={1,2,3,4,5,6}
]

\addplot[thick,red,dashed,mark=square,mark options={solid}] table[x=p,y expr=\thisrow{est_err}/\thisrow{err}]
{figures/data/pcurve_uns_pol.txt}
node[pos=.655,pin={[pin distance=.52cm, pin edge=solid]-90:{$\frac{\eta}{\norm{\Crl (\tA - \tAh)}}$}}] {};

\end{axis}
\end{tikzpicture}
\end{minipage}
\caption{[Smooth solution \cye{with limited regularity}~\eqref{eq_solution_1}] Uniform polynomial degree refinement.}
\label{figure_pconv_pol}
\end{figure}
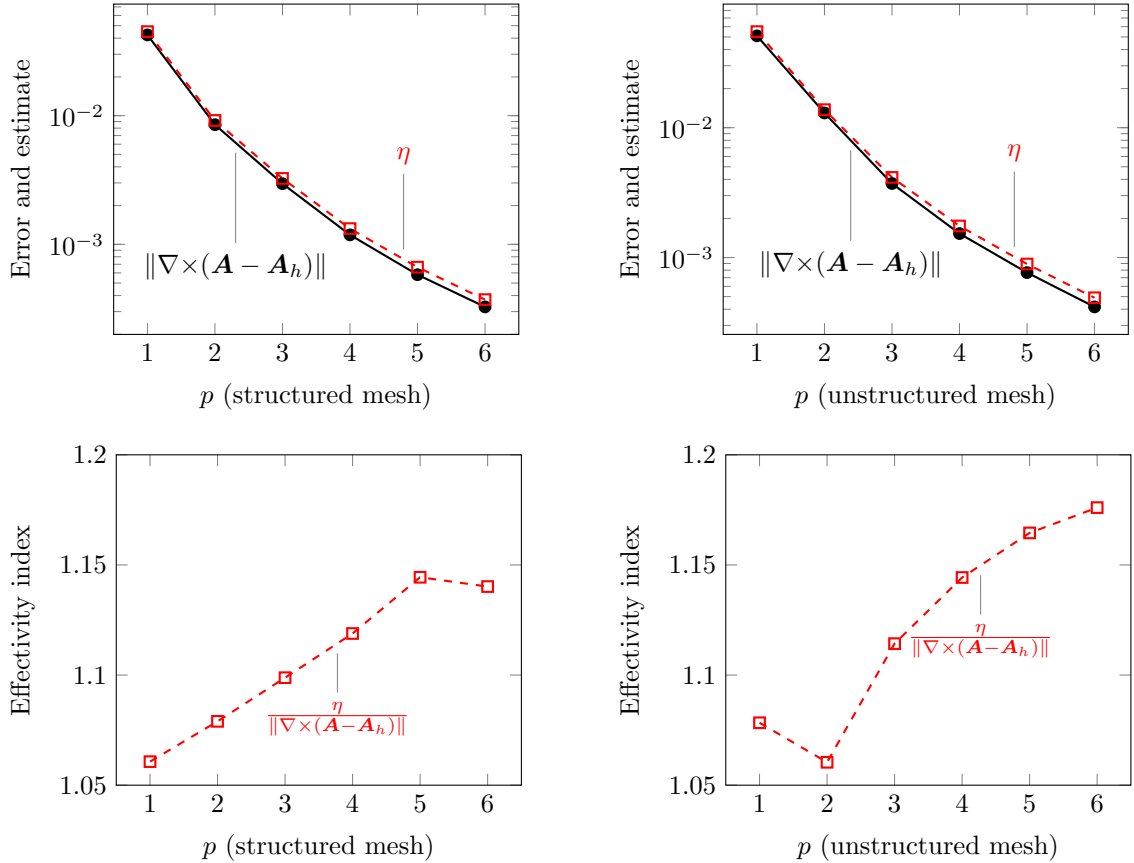

\subsection{Analytical solution with a general right-hand side}\label{sec_num_2}

We consider again the unit cube $\Om \eq (0,1)^3$, this time with a non-polynomial
right-hand side $\tj \eq 8\pi^2(\sin(2\pi \tx_2)\sin(2\pi \tx_3),0,0)$.
The associated solution is analytic,
\be
\label{eq_solution_2}
\tA \eq (\sin(2\pi\tx_2)\sin(2\pi\tx_3),0,0).
\ee

\begin{figure}
\begin{minipage}{.45\linewidth}
\begin{tikzpicture}
\begin{axis}%
[
	width=0.96\linewidth,
	xlabel={$h$ ($p=1$\cye{, structured meshes})},
	ylabel={Error and estimate},
	ymode = log,
	xmode = log,
	x dir = reverse
]

\addplot[thick,solid,black,mark=*,mark options={solid}]
table[x expr=.5*sqrt(3)/\thisrow{N},y=err]{figures/data/hcurve1_sin.txt}
node[pos=.7,pin={[pin distance=1cm, pin edge=solid]90:{$\norm{\Crl (\tA - \tAh)}$}}] {};
\addplot[thick,red,dashed,mark=square,mark options={solid}]
table[x expr=.5*sqrt(3)/\thisrow{N},y=est_err]{figures/data/hcurve1_sin.txt}
node[pos=.05,pin={[pin distance=1cm, pin edge=solid]-90:{$\eta$}}] {};
\addplot[thick,blue,dashdotted,mark=o,mark options={solid}]
table[x expr=.5*sqrt(3)/\thisrow{N},y expr=\thisrow{est_err}+\thisrow{osc}]%
{figures/data/hcurve1_sin.txt}
node[pos=.05,pin={[pin distance=.5cm, pin edge=solid]0:{$\eta_{\mathrm{tot}}$}}] {};

\SlopeTriangle{.65}{-.1}{.15}{-2}{$h^2$}{}

\end{axis}
\end{tikzpicture}
\end{minipage}
\qquad
\begin{minipage}{.45\linewidth}
\begin{tikzpicture}
\begin{axis}%
[
	width=0.96\linewidth,
	xlabel={\hspace{-0.35cm}$h$ ($p=2$\cye{, unstructured meshes})},
	ylabel={Error and estimate},
	ymode = log,
	xmode = log,
	x dir = reverse
]

\addplot[thick,solid,black,mark=*,mark options={solid}]
table[x=h,y=err]
{figures/data/hcurve2_uns_sin.txt}
node[pos=.7,pin={[pin distance=1cm, pin edge=solid]90:{$\norm{\Crl (\tA - \tAh)}$}}] {};
\addplot[thick,dashed,red,mark=square,mark options={solid}]
table[x=h,y=est_err]
{figures/data/hcurve2_uns_sin.txt}
node[pos=.1,pin={[pin distance=1cm, pin edge=solid]-90:{$\eta$}}] {};
\addplot[thick,dashdotted,blue,mark=o,mark options={solid}]
table[x=h,y expr=\thisrow{est_err}+\thisrow{osc}]
{figures/data/hcurve2_uns_sin.txt}
node[pos=.05,pin={[pin distance=.5cm, pin edge=solid]0:{$\eta_{\mathrm{tot}}$}}] {};

\SlopeTriangle{.65}{-.1}{.15}{-3}{$h^3$}{}

\end{axis}
\end{tikzpicture}
\end{minipage}

\medskip

\begin{minipage}{.45\linewidth}
\begin{tikzpicture}
\begin{axis}%
[
	width=0.96\linewidth,
	xlabel={$h$ ($p=1$\cye{, structured meshes})},
	ylabel={Effectivity index},
	xmode = log,
	x dir = reverse,
	ymin=.89,
	ymax=2.15
]

\addplot[thick,red,dashed,mark=square,mark options={solid}]
table[x expr=.5*sqrt(3)/\thisrow{N},y expr=\thisrow{est_err}/\thisrow{err}]
{figures/data/hcurve1_sin.txt}
node[pos=.6,pin={[pin distance=0.25cm, pin edge=solid]160:{$\frac{\eta}{\norm{\Crl (\tA - \tAh)}}$}}] {};

\addplot[thick,dashdotted,blue,mark=o,mark options={solid}]
table[x expr=.5*sqrt(3)/\thisrow{N},y expr=(\thisrow{est_err}+\thisrow{osc})/\thisrow{err}]
{figures/data/hcurve1_sin.txt}
node[pos=.5,pin={[pin distance=0.4cm, pin edge=solid]75:{$\frac{\eta_{\mathrm{tot}}}{\norm{\Crl (\tA - \tAh)}}$}}] {};


\end{axis}
\end{tikzpicture}
\end{minipage}
\qquad
\begin{minipage}{.45\linewidth}
\begin{tikzpicture}
\begin{axis}%
[
	width=0.96\linewidth,
	xlabel={\hspace{-0.35cm}$h$ ($p=2$\cye{, unstructured meshes})},
	ylabel={Effectivity index},
	xmode = log,
	x dir = reverse,
	ymin=.89,
	ymax=2.18
]

\addplot[thick,red,dashed,mark=square,mark options={solid}]
table[x=h,y expr=\thisrow{est_err}/\thisrow{err}]
{figures/data/hcurve2_uns_sin.txt}
node[pos=.7,pin={[pin distance=0.25cm, pin edge=solid]155:{$\frac{\eta}{\norm{\Crl (\tA - \tAh)}}$}}] {};

\addplot[thick,dashdotted,blue,mark=o,mark options={solid}]
table[x=h,y expr=(\thisrow{est_err}+\thisrow{osc})/\thisrow{err}]
{figures/data/hcurve2_uns_sin.txt}
node[pos=.5,pin={[pin distance=0.4cm, pin edge=solid]75:{$\frac{\eta_{\mathrm{tot}}}{\norm{\Crl (\tA - \tAh)}}$}}] {};


\end{axis}
\end{tikzpicture}
\end{minipage}
\caption{[Analytical solution~\eqref{eq_solution_2}] Uniform mesh refinement.}
\label{figure_hconv_sin}
\end{figure}
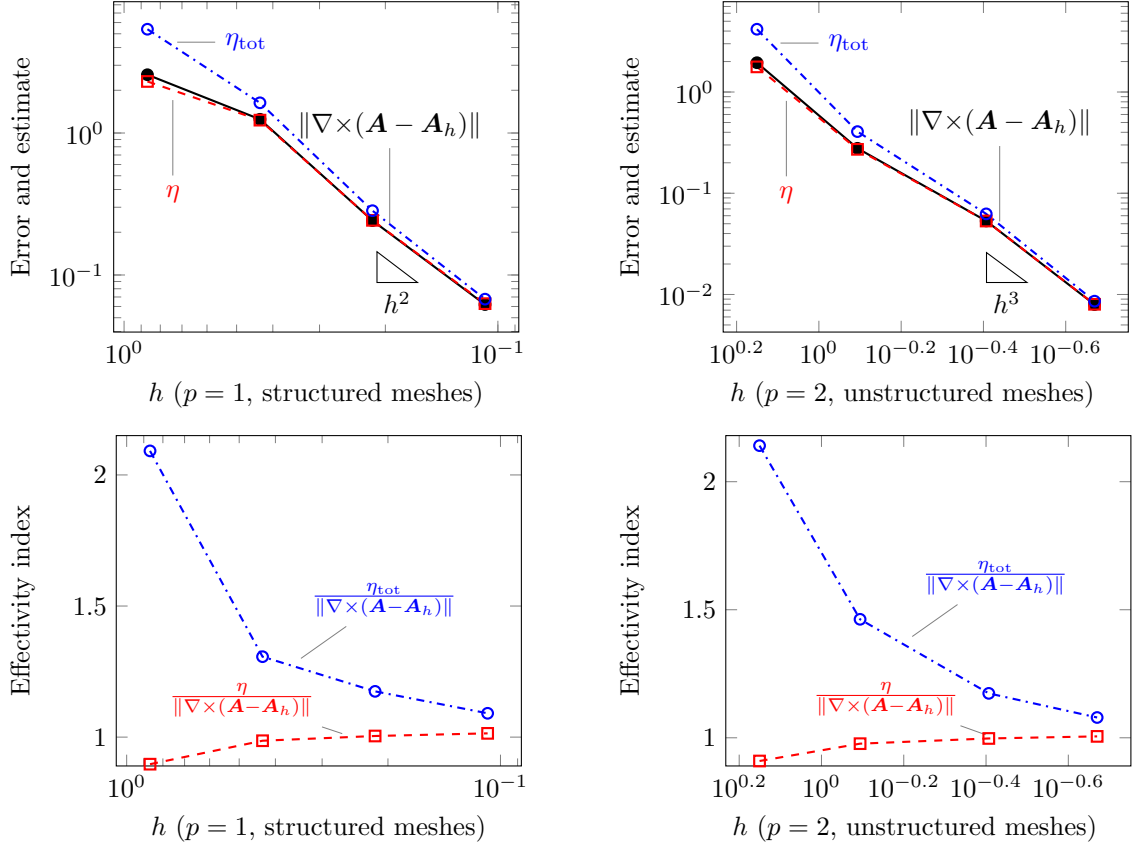

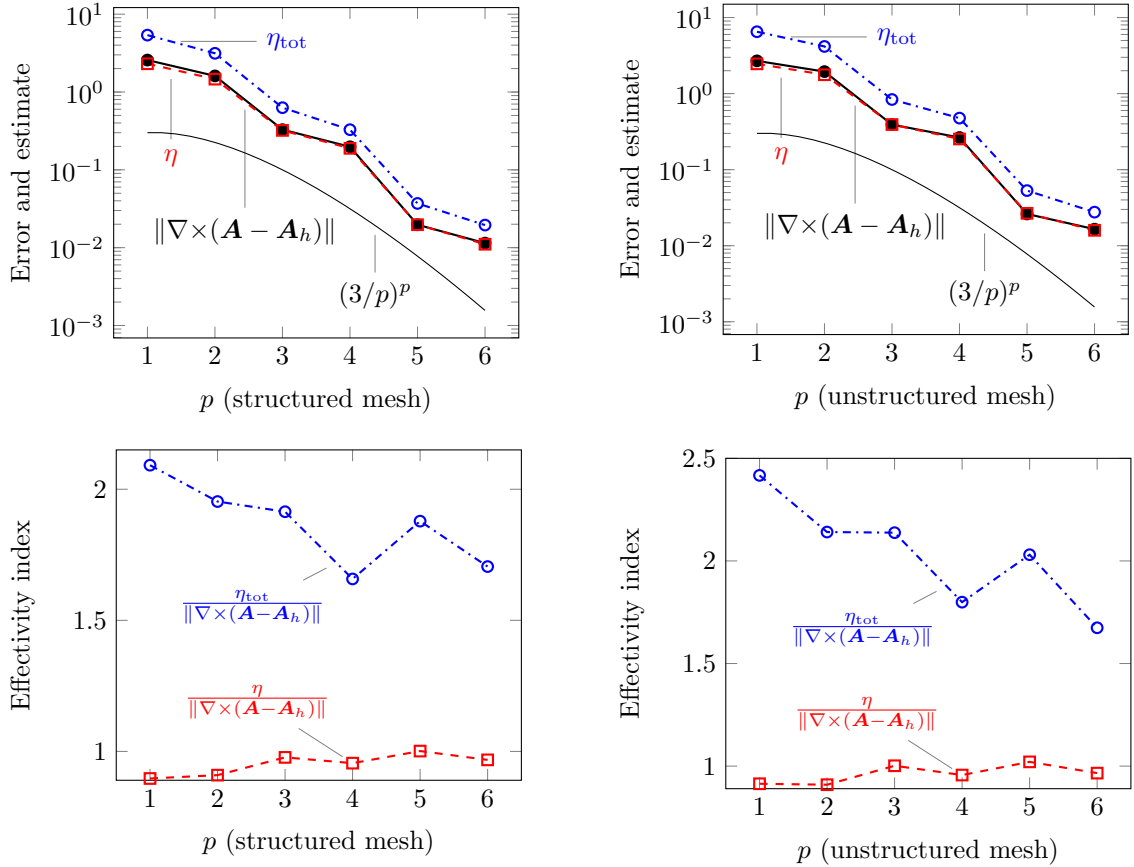
\begin{figure}
\begin{minipage}{.45\linewidth}
\begin{tikzpicture}
\begin{axis}%
[
	width=0.96\linewidth,
	xlabel={$p$ (structured mesh)},
	ylabel={Error and estimate},
	ymode = log,
	xtick={1,2,3,4,5,6}
]

\addplot[thick,solid,black,mark=*,mark options={solid}]
table[x=p,y=err]
{figures/data/pcurve_sin.txt}
node[pos=.25,pin={[pin distance=1.3cm, pin edge=solid]-90:{$\norm{\Crl (\tA - \tAh)}$}}] {};
\addplot[thick,dashed,red,mark=square,mark options={solid}]
table[x=p,y=est_err]
{figures/data/pcurve_sin.txt}
node[pos=.05,pin={[pin distance=.8cm, pin edge=solid]-90:{$\eta$}}] {};
\addplot[thick,dashdotted,blue,mark=o,mark options={solid}]
table[x=p,y expr=\thisrow{est_err}+\thisrow{osc}]
{figures/data/pcurve_sin.txt}
node[pos=.05,pin={[pin distance=1cm, pin edge=solid]0:{$\eta_{\mathrm{tot}}$}}] {};

\plot[domain=1:6] {.1*(3/x)^(x)}
node[pos=0.6,pin={[pin distance=0.5cm,pin edge=solid]-90:{$(3/p)^p$}}] {};


\end{axis}
\end{tikzpicture}
\end{minipage}
\qquad
\begin{minipage}{.45\linewidth}
\begin{tikzpicture}
\begin{axis}%
[
	width=0.96\linewidth,
	xlabel={$p$ (unstructured mesh)},
	ylabel={Error and estimate},
	ymode = log,
	xtick={1,2,3,4,5,6}
]

\addplot[thick,black,mark=*,mark options={solid}]
table[x=p,y=err]
{figures/data/pcurve_uns_sin.txt}
node[pos=.25,pin={[pin distance=1.3cm, pin edge=solid]-90:{$\norm{\Crl (\tA - \tAh)}$}}] {};
\addplot[thick,red,dashed,mark=square,mark options={solid}]
table[x=p,y=est_err]
{figures/data/pcurve_uns_sin.txt}
node[pos=.05,pin={[pin distance=.8cm, pin edge=solid]-90:{$\eta$}}] {};
\addplot[thick,dashdotted,blue,mark=o,mark options={solid}]
table[x=p,y expr=\thisrow{est_err}+\thisrow{osc}]
{figures/data/pcurve_uns_sin.txt}
node[pos=.05,pin={[pin distance=1cm, pin edge=solid]0:{$\eta_{\mathrm{tot}}$}}] {};

\plot[domain=1:6] {.1*(3/x)^(x)}
node[pos=0.6,pin={[pin distance=0.5cm,pin edge=solid]-90:{$(3/p)^p$}}] {};


\end{axis}
\end{tikzpicture}
\end{minipage}

\medskip

\begin{minipage}{.45\linewidth}
\begin{tikzpicture}
\begin{axis}%
[
	width=0.96\linewidth,
	xlabel={$p$ (structured mesh)},
	ylabel={Effectivity index},
	ymin=.89,
	ymax=2.15,
	xtick={1,2,3,4,5,6}
]

\addplot[thick,red,dashed,mark=square,mark options={solid}]
table[x=p,y expr=\thisrow{est_err}/\thisrow{err}]
{figures/data/pcurve_sin.txt}
node[pos=.6,pin={[pin distance=0.35cm, pin edge=solid]110:{$\frac{\eta}{\norm{\Crl (\tA - \tAh)}}$}}] {};

\addplot[thick,dashdotted,blue,mark=o,mark options={solid}]
table[x=p,y expr=(\thisrow{est_err}+\thisrow{osc})/\thisrow{err}]
{figures/data/pcurve_sin.txt}
node[pos=.55,pin={[pin distance=0.1cm, pin edge=solid]-95:{$\frac{\eta_{\mathrm{tot}}}{\norm{\Crl (\tA - \tAh)}}$}}] {};


\end{axis}
\end{tikzpicture}
\end{minipage}
\qquad
\begin{minipage}{.45\linewidth}
\begin{tikzpicture}
\begin{axis}%
[
	width=0.96\linewidth,
	xlabel={$p$ (unstructured mesh)},
	ylabel={Effectivity index},
	ymin=.89,
	ymax=2.5,
	xtick={1,2,3,4,5,6}
]

\addplot[thick,red,dashed,mark=square,mark options={solid}]
table[x=p,y expr=\thisrow{est_err}/\thisrow{err}]
{figures/data/pcurve_uns_sin.txt}
node[pos=.6,pin={[pin distance=0.35cm, pin edge=solid]110:{$\frac{\eta}{\norm{\Crl (\tA - \tAh)}}$}}] {};

\addplot[thick,dashdotted,blue,mark=o,mark options={solid}]
table[x=p,y expr=(\thisrow{est_err}+\thisrow{osc})/\thisrow{err}]
{figures/data/pcurve_uns_sin.txt}
node[pos=.55,pin={[pin distance=0.1cm, pin edge=solid]-95:{$\frac{\eta_{\mathrm{tot}}}{\norm{\Crl (\tA - \tAh)}}$}}] {};


\end{axis}
\end{tikzpicture}
\end{minipage}
\caption{[Analytical solution~\eqref{eq_solution_2}] Uniform polynomial degree refinement.}
\label{figure_pconv_sin}
\end{figure}

Figure~\ref{figure_hconv_sin} presents an $h$ convergence experiment
with the same settings as above. The optimal convergence rate
$h^{p+1}$ is observed for $\norm{\Crl(\tA-\tAh)}$. The oscillation-free
estimator $\eta$ closely follows the actual error,
with a possible slight underestimation\cye{, whereas the total estimator including data oscillation $\eta_{\mathrm{tot}} = \eta + \cye{\Clift} \eta_{\mathrm{osc}}$ from Theorem~\ref{thm_rel_eff} gives a guaranteed upper bound with a slight overestimation; as discussed in Section~\ref{sec_ineq}, we can take here $\Clift = 1$. In agreement with the theory, the influence of $\eta_{\mathrm{osc}}$ diminishes with mesh refinement, and we again numerically observe asymptotic exactness. }
We then consider a $p$ convergence test. In Figure~\ref{figure_pconv_sin},
we now observe the expected exponential convergence rate of the error and
a perfect behavior of the effectivity indices. \cye{More precisely, as the mesh is not refined here, $\eta_{\mathrm{osc}}$ does not necessarily go faster to zero than the error; this would be the case if the $hp$-version of~\eqref{eq_PF}, with $\bPi_p$ and $C h_\elm/(p+1)$ in place of respectively $\bPi_0$ and $h_\elm/\pi$, was used. We put forward here~\eqref{eq_PF}, where the is no unknown constant $C$, leading to a fully computable $\eta_{\mathrm{osc}}$.}

\subsection{Adaptivity with a singular solution}\label{sec_num_3}

Our last experiment features a singular solution
in a nonconvex domain, following~\cite{Chaum_Ern_Voh_Maxw_22,Ged_Gee_Per_a_post_Maxw_20}.
Specifically, we consider an L-shape example where $\Om \eq L \times (0,1)$, with
\[
L \eq \left \{
\tx = (r\cos\theta,r\sin\theta); \,
|\tx_1|,|\tx_2| \leq 1, \quad 0 \leq \theta \leq 3\pi/2
\right \}.
\]
The right-hand side $\tj$ is non-polynomial and chosen
such that
\be
\label{eq_solution_3}
\tA(\tx) = \big(0,0,\chi(r) r^\alpha \sin(\alpha \theta)\big),
\ee
where $\alpha \eq 3/2$, $r^2 \eq |\tx_1|^2 + |\tx_2|^2$,
$(\tx_1,\tx_2) = r(\cos\theta,\sin\theta)$, and $\chi:(0,1) \to \mathbb R$ is
a smooth cutoff function such that $\chi = 0$ in a neighborhood of $1$.
One easily checks that $\Dv \tA = 0$. Besides, since
$\Delta \left (r^\alpha \sin(\alpha \theta)\right ) = 0$ near the origin,
the right-hand side is non-singular (\ie,  $\tj \in \tLt$), and
the singularity appearing in the solution is solely due to the re-entrant
edge.

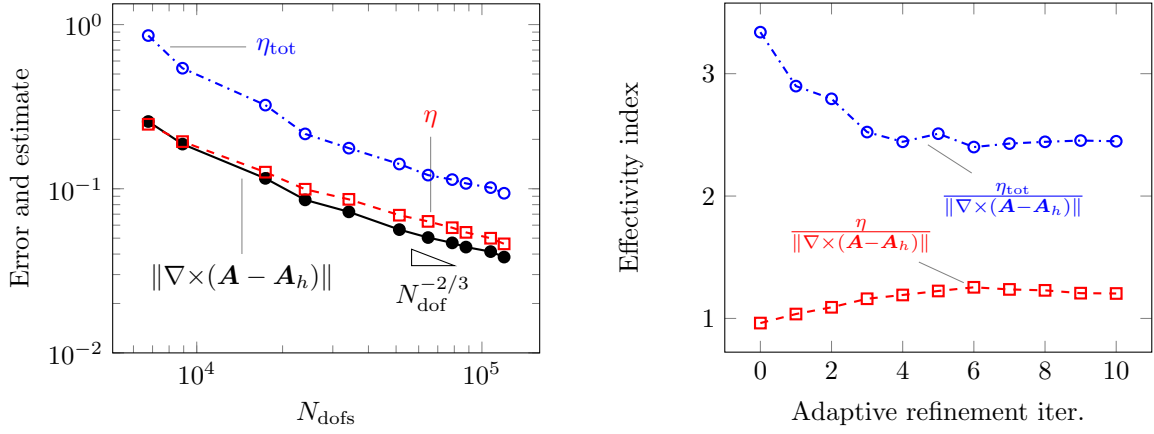
\begin{figure}
\begin{minipage}{.45\linewidth}
\begin{tikzpicture}
\begin{axis}%
[
	width=\linewidth,
	xlabel={$N_{\rm dofs}$},
	ylabel={Error and estimate},
	ymode = log,
	xmode = log,
	ymin = 0.01
]

\addplot[thick,black,solid,mark=*,mark options={solid}]
table[x=dof,y=err]{figures/data/adaptivity.txt}
node[pos=.29,pin={[pin distance=1cm, pin edge=solid]-90:{$\norm{\Crl (\tA - \tAh)}$}}] {};

\addplot[thick,red,dashed,mark=square,mark options={solid}]
table[x=dof,y=est]{figures/data/adaptivity.txt}
node[pos=.8,pin={[pin distance=1cm, pin edge=solid]90:{$\eta$}}] {};

\addplot[thick,dashdotted,blue,mark=o,mark options={solid}]
table[x=dof,y expr=\thisrow{est}+\thisrow{osc}]{figures/data/adaptivity.txt}
node[pos=.05,pin={[pin distance=1cm, pin edge=solid]0:{$\eta_{\mathrm{tot}}$}}] {};


\SlopeTriangle{.7}{-.1}{.25}{-2./3.}{$N_{\rm dof}^{-2/3}$}{}

\end{axis}
\end{tikzpicture}
\end{minipage}
\qquad
\begin{minipage}{.45\linewidth}
\begin{tikzpicture}
\begin{axis}%
[
	width=\linewidth,
	xlabel={Adaptive refinement iter.},
	ylabel={Effectivity index}
]

\addplot[thick,red,dashed,mark=square,mark options={solid}]
table[x=i,y expr=\thisrow{est}/\thisrow{err}]
{figures/data/adaptivity.txt}
node[pos=.6,pin={[pin distance=0.35cm, pin edge=solid]140:{$\frac{\eta}{\norm{\Crl (\tA - \tAh)}}$}}] {};

\addplot[thick,dashdotted,blue,mark=o,mark options={solid}]
table[x=i,y expr=(\thisrow{est}+\thisrow{osc})/\thisrow{err}]
{figures/data/adaptivity.txt}
node[pos=.45,pin={[pin distance=0.35cm, pin edge=solid]-70:{$\frac{\eta_{\mathrm{tot}}}{\norm{\Crl (\tA - \tAh)}}$}}] {};

\end{axis}
\end{tikzpicture}
\end{minipage}
\caption{[Singular solution~\eqref{eq_solution_3}] Adaptive mesh refinement.}
\label{figure_adaptivity}
\end{figure}

We couple our estimator with an adaptive strategy based on
D\"orfler's marking~\cite{Dorf_cvg_FE_96} for $\eta_\elm \eq \norm{\frh  - \Crl \tAh}_\elm$
and {\tt MMG3D}~\cite{Dobr_mmg3d} to build a series of adaptively refined meshes. We select $p=2$ and
an initial mesh made of $415$ elements.

The behaviors of the error and of the estimators $\eta$ \cye{and $\eta_{\mathrm{tot}}$} with respect to the number of degrees of freedom $N_{\rm dofs}$ are presented in Figure~\ref{figure_adaptivity}. \cye{Here we still take $\Clift = 1$ in front of $\eta_{\mathrm{osc}}$, though we do not anymore have a theoretical support for this.}
The effectivity ind\cye{ices} stay close to one, even on unstructured and locally refined meshes\cye{, with $\eta_{\mathrm{tot}} / \norm{\Crl (\tA - \tAh)}$ always above one}. Besides,
the optimal convergence rate is observed (it is limited to $-2/3$ when using isotropic
elements in the presence of an edge singularity\cye{, see~\cite[Section 4.2.3]{Apel_asnis_FEs_99}}).
This seems to indicate that the estimator is perfectly suited to drive adaptive \cye{mesh refinement}, and illustrates the local efficiency \cye{of Theorem~\ref{thm_rel_eff}}.

\begin{figure}
\includegraphics[width=\linewidth]{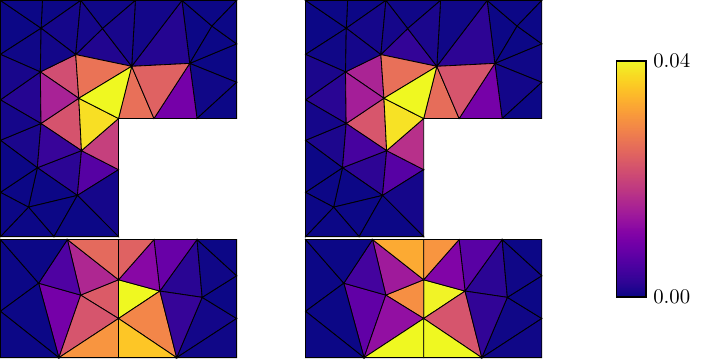}
\caption{[Singular solution~\eqref{eq_solution_3}].
Estimated (left) and actual (right) error distributions on the initial mesh.
Top view (top) and side view (bottom).}
\label{figure_adaptivity_it0}
\end{figure}

\begin{figure}
\includegraphics[width=\linewidth]{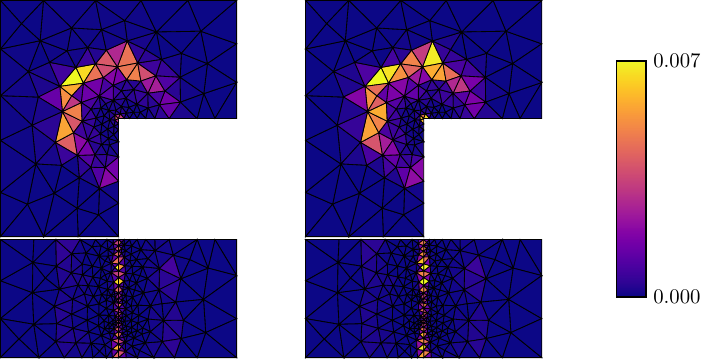}
\caption{[Singular solution~\eqref{eq_solution_3}]
Estimated (left) and actual (right) error distributions at adaptive mesh refinement iteration \#10.
Top view (top) and side view (bottom).}
\label{figure_adaptivity_it10}
\end{figure}

Finally, Figures~\ref{figure_adaptivity_it0}--\ref{figure_adaptivity_it10} present the meshes
generated by the adaptive algorithm, the estimators $\eta_\elm = \norm{\frh  - \Crl \tAh}_\elm$, and the elementwise errors $\norm{\Crl (\tA - \tAh)}_\elm$ (the top face and the faces sharing the re-entrant edge).
The meshes are refined close to the re-entrant edge, as expected. The
estimated error distribution closely matches the actual one, illustrating again
the local efficiency of the estimator.

\section{\cye{Technical details and proofs}} \label{sec_det_proofs}

\cye{This section collects some technical details and the proofs of all the claims above.}

\subsection{\cye{Equivalent form of~Assumption~\ref{ass_Ah}}}
\label{sec_orth}

Recall \cye{from Section~\ref{sec_mesh}} the piecewise affine ``hat'' function $\psia$ associated with the vertex $\ver \in \Vh$, as well as the notation $\Hsa$ from Section~\ref{sec_cont_sp_oma}. The following technical result holds true:

\bl[Equivalence of images by the curl operator]\label{lem_im_curl} There holds
\be \label{eq_im_curl}
    \Crl\left[\underset{\ver \in \Vh}{\operatorname{span}} \; \left (\psia|_\oma \Gr(\PP_1(\Ta)\cap\Hsa)  \right ) \right] = \Crl [\ND_0(\Th) \cap \HcD].
\ee
\el

\bp
Let $\ver \in \Vh$. For any $q_h \in \PP_1(\Ta) \cap \Hsa$, clearly $\psia|_\oma \Gr q_h$,
extended by zero outside of the patch subdomain $\oma$, lies in $\HcD$ (though in general
not in $\ND_0(\Th)$). Moreover, $\Crl(\psia|_\oma \Gr q_h) = \Gr \psia|_\oma \vp \Gr q_h$,
which is a piecewise constant vector-valued polynomial on the patch $\Ta$ whose extension
by zero outside of the patch subdomain $\oma$ has a continuous normal trace on interfaces
and zero normal trace on $\GD$. Thus, this extension belongs to the lowest-order divergence-free
Raviart--Thomas space, which implies $\Crl(\psia|_{\oma}\Gr q_h) = \Crl \tw_h$ on $\oma$
for $\tw_h$ which belongs to $\ND_0(\Th) \cap \HcD$. Thus, in~\eqref{eq_im_curl}, there holds
the inclusion $\subseteq$.

Conversely, following, \eg, Monk~\cite[Section 5.5.1]{Monk_FEs_Maxwell_03} or
Ern and Guermond~\cite[Section~15.1]{Ern_Guermond_FEs_I_21}, the space
$\ND_0(\Th) \cap \HcD$ is spanned by the set of the ``edge functions''
$\{\bpsi^\edg\}_{\edg \in \EhD}$, where $\EhD$ denotes the mesh edges
not lying in $\overline{\GD}$. If $\edg$ is the edge between vertices
$\ver,\vertt \in \Vh$, then $\bpsi^\edg = \psia \Gr \psib-\psib\Gr \psia$.
Moreover, if one of the vertices of $\edg$ lies in $\overline{\GD}$, we chose
the convention that $\ver \in \overline{\GD}$, so that we have
$(\psib-c_{\vertt})|_\oma \in \Hsa$ for some constant $c_{\vertt}$ in all cases.
Now, since
$
\Crl \bpsi^\edg
=
2 \Gr \psia \times \Gr \psib
=
2 \Crl(\psia \Gr \psib)
=
2 \Crl(\psia \Gr (\psib - c_{\vertt}))
$,
we have found $q_h \eq (\psib-c_{\vertt})|_\oma/2 \in \PP_1(\Ta) \cap \Hsa$ such that,
after zero extension, $\Crl(\psia|_\oma \Gr q_h) = \Crl \bpsi^\edg$, and the inclusion
$\supseteq$ in~\eqref{eq_im_curl} holds.
\ep

The following alternative \cye{formulation of} Assumption~\ref{ass_Ah} is crucial:

\bl[Patchwise orthogonality]
\label{lem_orth}
Let $\tj$ satisfy Assumption~\ref{ass_jh}. Then $\tAh$ satisfies Assumption~\ref{ass_Ah}
if and only if $\tAh \in \ND_p(\Th) \cap \HcD$ and
\be
\label{eq_orth}
(\psia \tj, \Gr q_h)_\oma \! + \! (\Gr \psia \vp (\Crl \tAh), \Gr q_h)_\oma \!
=
\! 0 \quad \forall q_h \! \in \! \PP_1(\Ta) \cap \Hsa, \, \forall \ver \in \Vh.\!\!\!
\ee
\el

\bp Since $\Gr \psia|_\oma \vp \Gr q_h = \Crl(\psia|_{\oma}\Gr q_h)$,
\[
    (\Gr \psia \vp (\Crl \tAh), \Gr q_h)_\oma = - (\Crl \tAh, \Gr \psia \vp \Gr q_h)_\oma = - (\Crl \tAh, \Crl(\psia \Gr q_h))_\oma.
\]
For any $\tv \in \HcN$ such that $\tj = \Crl \tv$, the Green theorem in turn gives
\[
    (\psia\tj,\Gr q_h)_{\oma} = (\tj,\psia\Gr q_h)_{\oma} = (\tv,\Crl(\psia\Gr q_h))_{\oma}.
\]
Finally, again by the Green theorem, for any $\tv_h \in \ND_0(\Th) \cap \HcD$,
\[
    (\tj,\tv_h) = (\Crl \tv, \tv_h) = (\tv, \Crl \tv_h).
\]
Applying these identities respectively in~\eqref{eq_orth} and~\eqref{eq_Ah_orth},
the assertion follows from Lemma~\ref{lem_im_curl}.
\ep

\subsection{\cye{Properties of the auxiliary fields \texorpdfstring{$\tha$, $\bdl_h$, and $\tda$}{} from Definition~\ref{def_div_free_dec}}}
\label{sec_dec_prop}

We collect here some important results on $\tha$, $\bdl_h$, and $\tda$
from~\eqref{eq_tha}--\eqref{eq_tda}. We start with the following application
of the self-standing result on over-constrained minimization in the
Raviart--Thomas spaces that we present in Appendix~\ref{app_over_constr_min} below. \cye{Let
$\eta_{\mathrm{osc},\tj}^\ver$ be defined as $\widetilde \eta_{\mathrm{osc},\tj}^\ver$ in~\eqref{eq_osc_j} but on the patch $\Ta$ only.}

\bl[Existence, uniqueness, and stability of $\tha$ from~\eqref{eq_tha}]
\label{lem_tsa}
There exists a unique solution $\tha$ to problem~\eqref{eq_tha} for all $\ver \in \Vh$. Moreover, it satisfies the stability estimate
\[
    \norm{\tha - \Gr \psia \vp (\Crl \tAh)}_\oma
    \ls \min_{\substack{\tv \in \Hdva\\ \Dv \tv = - \Gr \psia \scp \tj}}
    \norm{\tv - \Gr \psia \vp (\Crl \tAh)}_\oma + h_\oma^{-1} \eta_{\mathrm{osc},\tj}^\ver.
\]
\el

\bp
We choose $g^\ver \eq (- \Gr \psia \scp \tj)|_\oma$,
$\btau_{h}^\ver \eq (\Gr \psia \vp (\Crl \tAh))|_\oma$,
$q \eq p$ and verify the assumptions of Theorem~\ref{thm_over_constr_min}
in three steps. Note that $\Pi_\pp(\Gr \psia \scp \tj) = \Gr \psia \scp \bPi_\pp (\tj)$
and that $\norm{\Gr \psia \scp (\tj - \bPi_\pp(\tj))}_\elm
\leq \norm{\Gr \psia}_{\infty, \oma} \norm{\tj - \bPi_\pp(\tj)}_\elm
\ls h_\oma^{-1} \norm{\tj - \bPi_\pp(\tj)}_\elm$, where $\norm{\Gr \psia}_{\infty,\oma} \ls h_\oma^{-1}$ follows from the shape regularity
of the mesh, which gives rise to $h_\oma^{-1} \eta_{\mathrm{osc},\tj}^\ver$ from the data oscillation term in Theorem~\ref{thm_over_constr_min}.

{\em Step~1. Assumption~\eqref{eq_g_tau}.}
From~\eqref{eq_j_Hdv}, $g^\ver \in \Lti{\oma}$, so that the first condition
in~\eqref{eq_g_tau} is satisfied. From~\eqref{eq_Ah_conf}, in turn, on $\oma$,
it follows that $\Crl \tAh \in [\PP_p(\Ta)]^3$, see, \eg,
\cite[Corollary~2.3.2]{Bof_Brez_For_MFEs_13}, so that
$\btau_{h}^\ver = \Gr \psia \vp (\Crl \tAh) \in [\PP_p(\Ta)]^3 \subset \RT_p(\Ta) \subset \RT_\pp(\Ta)$.
Thus the second (polynomial) condition in~\eqref{eq_g_tau} is also satisfied.

{\em Step~2. Assumption~\eqref{eq_ga_comp}.}
For vertices $\ver \in \Vh$ such that $\ver \not \in \overline \GD$,
the Green theorem and $\tj \in \HdvN$ from~\eqref{eq_j_Hdv} together with
$\Dv \tj = 0$ from~\eqref{eq_j_div} imply
\[
    - (\Gr \psia \scp \tj,1)_\oma = - (\Gr \psia, \tj)_\oma = (\psia, \Dv \tj)_\oma = 0.
\]

{\em Step~3. Assumption~\eqref{eq_orth_A}.}
For any $q_h \in \PP_1(\Ta)\cap\Hsa$, again the Green theorem yields
\[
-(\Gr \psia \scp \tj, q_h)_\oma
\reff{eq_j_div}{=}
-(\Dv (\psia \tj), q_h)_\oma
=
(\psia \tj, \Gr q_h)_\oma,
\]
so that the patchwise orthogonality property~\eqref{eq_orth} implies
\be \label{eq_orth_cons}
   (\Gr \psia \vp (\Crl \tAh), \Gr q_h)_\oma - (\Gr \psia \scp \tj, q_h)_\oma = 0.
\ee
\ep

Similarly, an important part of the results of the following lemma are consequences
of Appendix~\ref{app_div_free_dec} below:

\bl[Auxiliary correction fields $\bdl_h$ and $\tda$]
\label{lem_tdh}
For $\bdl_h$ given by~\eqref{eq_td}, there holds
\be \label{eq_td_prop}
    \bdl_h \in \RT_{\pp}(\Th) \cap \HdvN \quad \text{ and } \quad \Dv \bdl_h = 0.
\ee
In addition, there exists a unique solution $\tda|_\elm$ to problems~\eqref{eq_tda} for all
tetrahedra $\elm \in \Th$ and all vertices $\ver \in \VK$, yielding the local divergence-free decomposition
\bse \ba
    \tda & \in \RT_{p+1}(\Ta) \cap \Hdva \quad \text{ and } \quad \Dv \tda= 0 \qquad \forall \ver \in \Vh, \label{eq_tda_prop}\\
    \bdl_h & = \suma \tda. \label{eq_tsat_dec}
\ea \ese
Moreover, for all tetrahedra $\elm \in \Th$ and all vertices $\ver \in \VK$, there holds
the local stability estimate
\be \label{eq_tsat_stab}
    \norm{\tda}_\elm \ls \norm{\bdl_h}_\elm.
\ee
\el

\bp
The patchwise contributions $\tha$ extended by zero outside of the patch subdomains $\oma$
belong to $\RT_{\pp}(\Th) \cap \HdvN$, so that the first property in~\eqref{eq_td_prop}
is immediate. The second property in~\eqref{eq_td_prop} then follows by the divergence
constraint in~\eqref{eq_tha}, the linearity of the projector $\Pi_\pp$, and the partition
of unity~\eqref{eq_PU}\cye{, since}
\[
\Dv \bdl_h
=
\suma \Dv \tha
=
\suma \Pi_\pp(- \Gr \psia \scp \tj)
=
\Pi_\pp\Bigg[ \suma - \Gr \psia \scp \tj \Bigg]
=
\Pi_\pp(0)
=
0.
\]

Let $\elm \in \Th$ and $\tr_h \in [\PP_0(\elm)]^{3}$ be fixed. Then definition~\eqref{eq_td},
which gives $\bdl_h|_\elm = \sum_{\vertt \in \VK} \thatt$, the partition of unity~\eqref{eq_PU},
which implies $\sum_{\vertt \in \VK}(\Gr \psib \vp (\Crl \tAh))|_\elm = \nv$, and the elementwise
orthogonality constraint in~\eqref{eq_tha} lead to
\[
    (\bdl_h, \tr_h)_{\elm} = \sum_{\vertt \in \VK} ( \thatt - \Gr \psib \vp (\Crl \tAh), \tr_h)_{\elm} = 0.
\]
This is condition~\eqref{eq_elm_orth_A}. Thus, Theorem~\ref{thm_div_free_dec_abs} can be
employed, where we choose $q \eq \pp$ together with $\qq \eq \pp$ for $p = 0$ and $\qq \eq \pp+1$
for $p \geq 1$. This implies the existence and uniqueness of solutions $\tda|_\elm$ to
problems~\eqref{eq_tda}, the properties~\eqref{eq_tda_prop}, the decomposition~\eqref{eq_tsat_dec},
and the stability bound~\eqref{eq_tsat_stab}. Note in particular that we only
employ~\eqref{eq_tsat_stab_2_A} with $\qq = q$ in the lowest-order case with $q=1$, \cye{whereas in other cases, we employ~\eqref{eq_tsat_stab_2_A} with $\qq = q+1$,}
so there is indeed no polynomial degree dependence in~\eqref{eq_tsat_stab}.
\ep

\subsection{\cye{Decomposition of the current density \texorpdfstring{$\tj$}{} and its stability \cye{from Theorem~\ref{thm_div_free_dec}}}} \cye{We are now ready to prove Theorem~\ref{thm_div_free_dec}.}

\bp[\cye{Proof of Theorem~\ref{thm_div_free_dec} (decomposition)}]
Property~\eqref{eq_div_free_dec_contr} is immediate since
$\psia \tj \in \Hdva$ in view of assumption~\eqref{eq_j_Hdv},
from~\eqref{eq_tha} which gives $\tha \in \RT_{\pp}(\Ta) \cap \Hdva$,
and from the first property in~\eqref{eq_tda_prop}.
Property~\eqref{eq_div_free_dec_dv} follows since
$\Dv(\psia \tj) = \Gr \psia \scp \tj $ in view of assumption~\eqref{eq_j_div} and
using $\Dv \tha= \Pi_\pp(- \Gr \psia \scp \tj) = - \Gr \psia \scp \bPi_\pp (\tj)$
from~\eqref{eq_tha} and $\Dv \tda= 0$, which
is the second property in~\eqref{eq_tda_prop}. Finally, \eqref{eq_div_free_dec_dec}
follows from the partition of unity~\eqref{eq_PU} which gives $\suma \psia \tj = \tj$
together with~\eqref{eq_td} and~\eqref{eq_tsat_dec}. When $\tj \in \RT_p(\Th) \cap \HdvN$,
\eqref{eq_div_free_dec_pol} immediately follows from~\eqref{eq_div_free_dec} and the fact
that $\psia \tj \in \RT_{p+1}(\Ta) \cap \Hdva$.
\ep

\bp[\cye{Proof of Theorem~\ref{thm_div_free_dec} (stability)}]
We develop
\[
    \tj^\ver - \tj_h^\ver = \Gr \psia \vp (\Crl \tA) - \tha + \tda = \Gr \psia \vp (\Crl (\tA - \tAh)) - (\tha - \Gr \psia \vp (\Crl \tAh)) + \tda.
\]
For the first term above, we immediately see
\[
    \norm{\Gr \psia \vp (\Crl (\tA - \tAh))}_\oma \leq \norm{\Gr \psia}_{\infty,\oma} \norm{\Crl (\tA - \tAh)}_\oma.
\]
For the second term above, we employ Lemma~\ref{lem_tsa} with $\tv = \Gr \psia \vp (\Crl \tA)$, which lies in $\Hdva$ with divergence equal to $- \Gr \psia \scp \tj$ by virtue of~\eqref{eq_tsea}, which leads to
\be \label{eq_b1}
\norm{\tha - \Gr \psia \vp (\Crl \tAh)}_\oma
\ls
\norm{\Gr \psia \vp (\Crl (\tA - \tAh))}_\oma
+
h_\oma^{-1} \eta_{\mathrm{osc},\tj}^\ver.
\ee
For the last term, we first recall~\eqref{eq_tsat_stab},
\ie, $\norm{\tda}_\elm \ls \norm{\bdl_h}_\elm$ for every $\elm \in \Ta$.
Now definition~\eqref{eq_td}, the partition of unity~\eqref{eq_PU},
and the triangle inequality imply
\[
\norm{\bdl_h}_\elm
=
\norm[\Bigg]{\sum_{\vertt \in \VK}(\thatt - \Gr \psib \vp (\Crl \tAh))}_\elm
\leq
\sum_{\vertt \in \VK} \norm{\thatt - \Gr \psib \vp (\Crl \tAh)}_\omat,
\]
which extends by one layer beyond the patch $\oma$ \cye{and can be bounded by~\eqref{eq_b1}.} The shape regularity of the mesh ensures that
$\norm{\Gr \psia}_{\infty,\oma} \ls h_\oma^{-1}$
and $\norm{\Gr \psib}_{\infty,\omat} \simeq \norm{\Gr \psia}_{\infty,\oma}$
for all vertices $\vertt$ in the patch $\Ta$. Hence, \eqref{eq_dec_stab} follows
upon combining the above developments.
\ep

\subsection{\cye{Equilibrated flux reconstruction from Section~\ref{sec_equil} and its stability}} \label{sec_equil_proofs}

\cye{To prove Theorem~\ref{thm_EFR}, we rely on the following crucial result:}

\begin{theoremclr}[$p$-robust \cye{$\Hci{\oma}$} stability]
\label{thm_p_rob}
For a vertex $\ver \in \Vh$, let $\tAh \in \ND_{p}(\Ta) \cap \Hci{\oma}$ and
$\cye{\bar \tj_h^\ver} \in \RT_{p+1}(\Ta) \cap \Hdva$ with $\Dv \cye{\bar \tj_h^\ver} = 0$ be given.
Then
\be \label{eq_p_rob}
    \min_{\substack{\tv_h \in \ND_{p+1}(\Ta) \cap \Hca\\ \Crl \tv_h= \cye{\bar \tj_h^\ver}}} \norm{\tv_h - \psia(\Crl \tAh)}_\oma \ls \min_{\substack{\tv \in \Hca\\ \Crl \tv= \cye{\bar \tj_h^\ver}}} \norm{\tv - \psia(\Crl \tAh)}_\oma.
\ee
\end{theoremclr}

On a single tetrahedron $\elm$ in place of the vertex patch $\Ta$, \cye{Theorem~\ref{thm_p_rob}}
follows by the seminal contributions of Costabel and
McIntosh~\cite[Proposition~4.2]{Cost_McInt_Bog_Poinc_10} and Demkowicz
\eal~\cite[Theorem~7.2]{Demk_Gop_Sch_ext_II_09}, see~\cite[Theorem~2]{Chaum_Ern_Voh_curl_elm_20}.
On an edge patch, such a result has been established
in~\cite[Theorem~3.1]{Chaum_Ern_Voh_Maxw_22}. \cye{The further extension to a vertex
patch has recently been established in~\cite[\tcf{Theorem 3.3\cye{, see also Corollary~4.3}}]{Chaum_Voh_p_rob_3D_H_curl_22}.}

\bp[\cye{Proof of Theorem~\ref{thm_EFR} (equilibration)}]
Property~\eqref{eq_frh_prop_sum} follows immediately
from $\cye{\frh^\ver} \in \ND_{p+1}(\Ta) \cap \Hca$ of~\eqref{eq_frha} \cye{and~\eqref{eq_frh}}. \cye{For piecewise polynomial $\tj \in \RT_p(\Th) \cap \HdvN$, $\Crl \frh^\ver = \bar \tj_h^\ver = \tj_h^\ver$ from~\eqref{eq_frha} and~\eqref{eq_bjha}.} Property~\eqref{eq_frh_prop_curl} is then an direct consequence of~\eqref{eq_div_free_dec_dec} and~\eqref{eq_frh}. \ep

\bp[\cye{Proof of Theorem~\ref{thm_EFR} (stability)}]
Fix a vertex $\ver \in \Vh$ and use
$\tj^\ver = \psia \tj + \Gr \psia \vp (\Crl \tA) = \Crl (\psia (\Crl \tA))$
as in property~\eqref{eq_psia_E}. This implies
$ (\tj^\ver, \tv)_\oma = (\psia(\Crl \tA), \Crl \tv)_\oma$ for any $\tv \in \Hcdual{\oma}$.
Then \cye{Theorem~\ref{thm_p_rob}} and a primal--dual equivalence as in,
\eg, \cite[Lemma~5.5]{Chaum_Ern_Voh_Maxw_22} imply
\ban
\norm{\frh^\ver  - \psia(\Crl \tAh)}_\oma
\ls \!\! {} &
\min_{\substack{
\tv \in \Hca
\\
\Crl \tv = \cye{\bar \tj_h^\ver}
}}
\norm{\tv - \psia(\Crl \tAh)}_\oma
\\
= \!\! {} &
\sup_{\substack{\tv \in \Hcdual{\oma}\\ \norm{\Crl \tv}_\oma = 1}}
\big\{(\cye{\bar \tj_h^\ver}, \tv)_\oma - (\psia(\Crl \tAh), \Crl \tv)_\oma\big\}
\\
\leq \!\! {} &
\sup_{\substack{\tv \in \Hcdual{\oma}\\ \norm{\Crl \tv}_\oma = 1}}
(\cye{\bar \tj_h^\ver} - \tj^\ver, \tv)_\oma  + \norm{\psia (\Crl (\tA  - \tAh))}_\oma
\\
\leq \!\! {} &
\sup_{\substack{\tv \in \Hcdual{\oma}\\ \norm{\Crl \tv}_\oma = 1}}
(\cye{\bar \tj_h^\ver} - \tj^\ver, \tv)_\oma + \norm{\Crl (\tA  - \tAh)}_\oma.
\ean
We are thus left to treat the first term above.

Fix $\tv \in \Hcdual{\oma}$ with $\norm{\Crl \tv}_\oma = 1$. Consider $q \in \Hsa$ such that
\[
    (\Gr q, \Gr w)_\oma = (\tv, \Gr w)_\oma \qquad \forall w \in \Hsa.
\]
Then $\tilde \tv \eq \tv - \Gr q$ lies in both $\Hcdual{\oma}$ and $\Hdva$ and is divergence-free,
$\Dv \tilde \tv =0$. Thus, the Poincar\'e--Friedrichs--Weber inequality~\eqref{eq_weber_patch}
implies
\be \label{eq_Web}
\norm{\tilde \tv}_\oma \ls h_\oma \norm{\Crl \tilde \tv}_\oma
=
h_\oma \norm{\Crl \tv}_\oma
=
h_\oma.
\ee
Note that $\cye{\bar \tj_h^\ver} - \tj^\ver \in \Hdva$ with
$\Dv (\cye{\bar \tj_h^\ver} - \tj^\ver) = 0$; indeed, this follows from~\eqref{eq_je_a}--\eqref{eq_je_div} together with~\eqref{eq_bjha}. Thus, the Green theorem gives
\be \label{eq_orth_grad}
    (\cye{\bar \tj_h^\ver} - \tj^\ver, \Gr q)_\oma = 0.
\ee
\cye{Thanks to this crucial property, we can play in $\tilde \tv$ and use~\eqref{eq_Web}: employing additionally} the Cauchy--Schwarz inequality and the triangle inequality, we have
\be \label{eq_jj} \bs
(\cye{\bar \tj_h^\ver} - \tj^\ver, \tv)_\oma
&=
(\cye{\bar \tj_h^\ver} - \tj^\ver, \tilde \tv)_\oma
\leq
\norm{\cye{\bar \tj_h^\ver} - \tj^\ver}_\oma \norm{\tilde \tv}_\oma
\\
&
\cye{\ls}
h_\oma \big[\norm{\cye{\bar \tj_h^\ver} - \tj_h^\ver}_\oma
+
\norm{\tj_h^\ver - \tj^\ver}_\oma \big],
\es \ee
and we conclude by~\eqref{eq_dec_stab} from Theorem~\ref{thm_div_free_dec}.
\ep

\subsection{\cye{A posteriori error estimates from Section~\ref{sec_a_post}}}

\cye{We can finally prove Theorem~\ref{thm_rel_eff}.}

\bp[\cye{Proof of Theorem~\ref{thm_rel_eff} (reliability)}]
For a piecewise polynomial current density, $\tj \in \RT_p(\Th) \cap \HdvN$,
Theorem~\ref{thm_EFR} implies $\frh \in \HcN$ with $\Crl \frh  = \tj$. Thus,
in this case the claim follows with $\eta_{\mathrm{osc}} = 0$ by the Prager--Synge
theorem~\cite{Prag_Syng_47} in the $\HC$-context, see, \eg,
\cite[Theorem~10]{Braess_Scho_a_post_edge_08} or~\cite[Theorem~3.1]{Ged_Gee_Per_a_post_Maxw_20}.

In general, we proceed as follows. Since $\tA, \tAh \in \HcD$,
\[
    \norm{\Crl (\tA - \tAh)} = \max_{\substack{\tv \in \HcD\\ \norm{\Crl \tv} = 1}} (\Crl (\tA  - \tAh), \Crl \tv).
\]
Fix $\tv \in \HcD$ with $\norm{\Crl \tv} = 1$ and consider $\tw$ from~\eqref{eq_estimate_lift}. Note that since $\frh \in \HcN$ from Theorem~\ref{thm_EFR}, the Green theorem and $\Crl \tw = \Crl \tv$ give
\[
    (\Crl \frh, \tw) = (\frh, \Crl \tw) = (\frh, \Crl \tv).
\]
Similarly, $\Crl \tw = \Crl \tv$ and the weak solution characterization~\eqref{eq_maxwell_weak_II} lead to
\[
    (\Crl \tA, \Crl \tv) = (\Crl \tA, \Crl \tw) = (\tj, \tw).
\]
Thus
\[
    (\Crl (\tA  - \tAh), \Crl \tv) = (\tj - \Crl \frh, \tw) + (\frh - \Crl \tAh, \Crl \tv).
\]
The second term is trivially bounded by the estimator $\eta$ via the Cauchy--Schwarz inequality, so that we are left with bounding the first one.

\cye{Property~\eqref{eq_osc_eq} and the additional orthogonality constraint in~\eqref{eq_bjha} lead to
\[
    (\tj - \Crl \frh, \tw) = \sum_{\elm \in \Th} \Bigg(\sum_{\ver \in \VK} (\tj_h^\ver - \cye{\bar \tj_h^\ver}), \tw\Bigg)_\elm = \sum_{\elm \in \Th} \big(\tj - \Crl \frh, \tw - \bPi_0(\tw)\big)_\elm.
\]
Consequently, the Poincar\'e inequality~\eqref{eq_PF}, \eqref{eq_estimate_lift}, and $\norm{\Crl \tv} = 1$ give}
\[\cye{
    (\tj - \Crl \frh, \tw) \leq \sum_{\elm \in \Th} \eta_{\mathrm{osc},\elm} \norm{\Gr \tw}_\elm \leq \eta_{\mathrm{osc}} \norm{\Gr \tw} \leq \Clift \eta_{\mathrm{osc}} \norm{\Crl \tv} = \Clift \eta_{\mathrm{osc}}.}
\]
\ep

\cye{\br[Comparison with~\eqref{eq_jj}] Above, we could also write
\[
    (\tj - \Crl \frh, \tw) = \suma (\tj_h^\ver - \cye{\bar \tj_h^\ver}, \tw)_\oma,
\]
where the terms in the sum are similar to~\eqref{eq_jj} from Section~\ref{sec_equil_proofs}. In contrast to~\eqref{eq_jj}, it seems that we cannot pass through the Poincar\'e--Friedrichs--Weber inequality~\eqref{eq_weber_patch} in absence of an exactly divergence-free field \cye{(recall from~\eqref{eq_div_free_dec_dv} that $\Dv \tj_h^\ver = \Gr \psia \scp (\tj - \bPi_\pp(\tj))$ only in general)}, and rather need to resort to the switch from $\tv \in \HcD$ to $\tw$ of~\eqref{eq_estimate_lift} and to make use of the Poincar\'e inequality~\eqref{eq_PF}.\er}

\bp[\cye{Proof of Theorem~\ref{thm_rel_eff} (efficiency)}]
\cye{Property~\eqref{eq_osc_eq}, the triangle inequality, and definition~\eqref{eq_osc_jha} immediately lead to $\eta_{\mathrm{osc},\elm} \leq \sum_{\ver \in \VK} \eta_{\mathrm{osc},\tj_h^\ver}^\ver$. Moreover,} Definition~\eqref{eq_frh} together with the partition of unity~\eqref{eq_PU} imply
\[
\norm{\frh  - \Crl \tAh}_\elm
=
\norm[\Bigg]{\sum_{\ver \in \VK}(\frh^\ver  - \psia(\Crl \tAh))}_\elm
\leq
\sum_{\ver \in \VK} \norm{\frh^\ver  - \psia(\Crl \tAh)}_\oma.
\]
Thus employing Theorem~\ref{thm_EFR} concludes the proof.
\ep

\appendix

\section{Over-constrained minimization in Raviart--Thomas' spaces}
\label{app_over_constr_min}

\subsection{\cgr{Assumption and statement of the over-constrained minimization}}

In this appendix, we consider a fixed mesh vertex $\ver \in \Vh$.
Let an integer $q \geq 0$ be fixed and set $\qq \eq \min\{q,1\}$.
We employ the notation of Section~\ref{sec_not} and in particular
recall that $\ls$ means smaller or equal to up to a constant only
depending on the mesh shape-regularity parameter $\kappa_\Th$.
We also assume a polynomial form, mean value zero, and patchwise orthogonality
conditions \cye{on the two data $g^\ver$ and $\btau_h^\ver$}:

\bas[Data $g^\ver$ and $\btau_h^\ver$]
\label{ass_data}
The data $g^\ver$ and $\btau_h^\ver$ satisfy
\bse\ba
g^\ver & \in \Lti{\oma} \quad \text{ and } \quad \btau_h^\ver \in \RT_{\qq}(\Ta), \label{eq_g_tau}
\\
(g^\ver, 1)_\oma & = 0 \quad \text{ when } \ver \not \in \overline \GD, \label{eq_ga_comp}
\\
(\btau_h^\ver, \Gr q_h)_\oma + (g^\ver, q_h)_\oma & = 0
\qquad
\forall q_h \in \PP_1(\Ta)\cap\Hsa. \label{eq_orth_A}
\ea \ese
\eas

We consider the following \cye{(seemingly} over-constrained\cye{)} minimization problem in the Raviart--Thomas space
$\RT_{\qq}(\Ta) \cap \Hdva$:
\be \label{eq_tha_A}
\tha
\eq
\arg \min_{\substack{
\tv_h \in \RT_{\qq}(\Ta) \cap \Hdva
\\
\Dv \tv_h= \Pi_\qq(g^\ver)
\\
(\tv_h, \tr_h)_\elm = (\btau_h^\ver, \tr_h)_\elm
\quad
\forall \tr_h \in [\PP_0(\elm)]^3, \, \forall \elm \in \Ta}}
\norm{\tv_h - \btau_h^\ver}_\oma^2.
\ee

The following result is of independent interest:

\bt[Over-constrained minimization in the Raviart--Thomas spaces]
\label{thm_over_constr_min}
Let Assumption~\ref{ass_data} hold.
Then there exists a unique solution $\tha$ to problem~\eqref{eq_tha_A},
satisfying the stability estimate
\[
\norm{\tha - \btau_h^\ver}_\oma
\ls
\min_{\substack{\tv \in \Hdva\\ \Dv \tv = g^\ver}}
\norm{\tv - \btau_h^\ver}_\oma
+
\Bigg\{\sum_{\elm \in \Ta}\Big(
\frac{h_\elm}{\pi} \norm{\Pi_\qq(g^\ver) - g^\ver}_\elm
\Big)^2\Bigg\}^\ft.
\]
\et

\subsection{\cgr{Auxiliary conventional minimization}}

\cgr{In addition to~\eqref{eq_tha_A}, i}t will be useful to \cgr{also} consider
\be \label{eq_tsab}
\thab
\eq
\arg \min_{\substack{
\tv_h \in \RT_{\qq}(\Ta) \cap \Hdva
\\
\Dv \tv_h= \Pi_\qq(g^\ver)
}}
\norm{\tv_h - \btau_h^\ver}_\oma^2.
\ee
Minimizations~\eqref{eq_tsab} are in a conventional format
in that the constraints only concern normal trace and divergence.
Moreover, they fulfill the following important property:

\bl[Existence, uniqueness, and stability of $\thab$ from~\eqref{eq_tsab}] \label{lem_tsab}
Let Assumption~\ref{ass_data} hold.
Then there exists a unique solution $\thab$ to problem~\eqref{eq_tsab},
satisfying the stability estimate
\be \label{eq_tsab_stab}
\norm{\thab - \btau_h^\ver}_\oma
\ls
\min_{\substack{\tv \in \Hdva \\ \Dv \tv = g^\ver}}
\norm{\tv - \btau_h^\ver}_\oma
+
\Bigg\{\sum_{\elm \in \Ta} \Big(
\frac{h_\elm}{\pi} \norm{\Pi_\qq(g^\ver) - g^\ver}_\elm
\Big)^2\Bigg\}^\ft.
\ee
\el

\bp
Existence and uniqueness of $\thab$ from~\eqref{eq_tsab} are classical following, \eg,
\cite{Bof_Brez_For_MFEs_13}, thanks to the Neumann boundary compatibility
condition~\eqref{eq_ga_comp}; note that this implies $(\Pi_\qq(g^\ver), 1)_\oma = 0$ when
$\ver \not \in \overline \GD$.
Moreover, since $\Pi_\qq(g^\ver) \in \PP_{\qq}(\Ta)$ and $\btau_h^\ver \in \RT_{\qq}(\Ta)$,
taking $p=\qq$, $\btau_p = - \btau_h^\ver$, and $r_\elm = (\Pi_\qq(g^\ver) - \Dv \btau_h^\ver)|_\elm$
in~\cite[Corollaries~3.3 and 3.6]{Ern_Voh_p_rob_3D_20} for an interior vertex $\ver$
and~\cite[Corollary~3.8]{Ern_Voh_p_rob_3D_20}
\tcf{\cye{and}~\cite[Proposition 3.1]{Chaum_Voh_p_rob_3D_H_curl_22}} for a boundary vertex $\ver$ lead to
\[
\norm{\thab - \btau_h^\ver}_\oma
\ls
\min_{\substack{\tv \in \Hdva\\ \Dv \tv = \Pi_\qq(g^\ver)}}
\norm{\tv - \btau_h^\ver}_\oma = \norm{\Gr \tilde r^\ver}_\oma.
\]
The equality above is a classical primal--dual equivalence, with $\tilde r^\ver \in \Hsa$ given by
\[
(\Gr \tilde r^\ver, \Gr v)_\oma
=
(\btau_h^\ver, \Gr v)_\oma + (\Pi_\qq(g^\ver), v)_\oma
\qquad
\forall v \in \Hsa.
\]
Thus, as in, \eg, \cite[Theorem~3.17]{Ern_Voh_p_rob_15},
\ban
\norm{\Gr \tilde r^\ver}_\oma
&=
\max_{\substack{v \in \Hsa\\ \norm{\Gr v}_\oma=1}}
\{(\btau_h^\ver, \Gr v)_\oma + (\Pi_\qq(g^\ver), v)_\oma \}
\\
&=
\max_{\substack{v \in \Hsa\\ \norm{\Gr v}_\oma=1}}
\{ (\btau_h^\ver, \Gr v)_\oma + (g^\ver, v)_\oma + (\Pi_\qq(g^\ver) - g^\ver, v)_\oma \}.
\ean
The projection orthogonality and the elementwise Poincar\'e inequality then lead to
\ban
|(\Pi_\qq(g^\ver) - g^\ver, v)_\oma|
&=
\Bigg|\sum_{\elm \in \Ta} (\Pi_\qq(g^\ver) - g^\ver, v - \Pi_0 v)_\elm \Bigg|
\\
&\leq
\Bigg\{\sum_{\elm \in \Ta} \Big(
\frac{h_\elm}{\pi} \norm{\Pi_\qq(g^\ver) - g^\ver}_\elm
\Big)^2\Bigg\}^\ft
\norm{\Gr v}_\oma,
\ean
and~\eqref{eq_tsab_stab} follows, since
\[
\max_{\substack{v \in \Hsa\\ \norm{\Gr v}_\oma=1}}
\{ (\btau_h^\ver, \Gr v)_\oma + (g^\ver, v)_\oma \}
=
\min_{\substack{\tv \in \Hdva\\ \Dv \tv = g^\ver}}
\norm{\tv - \btau_h^\ver}_\oma
\]
by the same primal--dual equivalence argument.
\ep

\subsection{\cgr{Auxiliary first-order over-constrained minimization and proof of Theorem~\ref{thm_over_constr_min}}}

Let, \cgr{in addition to~\eqref{eq_tha_A} and~\eqref{eq_tsab},} the {\em first-order} Raviart--Thomas piecewise polynomial $\tzab$ be given by
\be \label{eq_tzab}
\tzab
\eq
\arg \min_{\substack{
\tv_h \in \RT_1(\Ta) \cap \Hdva
\\
\Dv \tv_h= 0
\\
(\tv_h, \tr_h)_\elm = (\btau_h^\ver - \thab, \tr_h)_\elm
\quad \forall \tr_h \in [\PP_0(\elm)]^3, \, \forall \elm \in \Ta}}
\norm{\tv_h - \btau_h^\ver + \thab}_\oma^2.
\ee
The field $\tzab$ can be seen as the correction of $\thab$ \cgr{from~\eqref{eq_tsab}} necessary to fulfill the constraints
on the elementwise product with piecewise vector-valued constants in~\eqref{eq_tha_A}.
As one might expect, the patchwise orthogonality assumption~\eqref{eq_orth_A}, turns to be the
key for the following crucial technical result:

\bl[Existence, uniqueness, and stability of $\tzab$ from~\eqref{eq_tzab}]
\label{lem_tzab}
Let Assumption~\ref{ass_data} hold.
Then there exists a unique solution $\tzab$ to problem~\eqref{eq_tzab},
and the following stability estimate holds true:
\be \label{eq_tzab_stab}
    \norm{\tzab}_\oma \ls \norm{\btau_h^\ver - \thab}_\oma.
\ee
\el

\cgr{We postpone the proof Lemma~\ref{lem_tzab} to the sections below; let us now first show that Lemma~\ref{lem_tzab} implies the results announced in Theorem~\ref{thm_over_constr_min}.}

\bp[Proof of Theorem~\ref{thm_over_constr_min}]
It follows straightforwardly from~\eqref{eq_tsab} and~\eqref{eq_tzab} that $\thab + \tzab$
lies in the minimization set of~\eqref{eq_tha_A}. Consequently, the existence and uniqueness
of~\eqref{eq_tha_A} follows since the minimized functional in~\eqref{eq_tha_A} is convex. Moreover,
the triangle inequality together with Lemma~\ref{lem_tzab} implies
\[
\norm{\tha - \btau_h^\ver}_\oma
\leq
\norm{\thab + \tzab - \btau_h^\ver}_\oma
\leq
\norm{\tzab}_\oma + \norm{\thab - \btau_h^\ver}_\oma
\ls \norm{\thab - \btau_h^\ver}_\oma,
\]
and we conclude by Lemma~\ref{lem_tsab}.
\ep

\subsection{\cgr{Piola mappings}}

\cgr{In order to prove the technical results below, we will rely on}
\tcf{Piola mappings.
An extensive description may be found in, \eg,
\cite[Chapters~7.2 and~9.2]{Ern_Guermond_FEs_I_21},
and we \cgr{only} list here the essential properties we need.

If $U,V \subset \mathbb R^3$ are open sets with Lipschitz
boundaries and $\phi: U \to V$ is a bilipschitz
mapping, the gradient-, curl-, diverence-preserving and broken Piola mappings
are the applications $\phi^{\rm g}: L^2(U) \to L^2(V)$,
$\phi^{\rm c}, \phi^{\rm d}: \tL^2(U) \to \tL^2(V)$, and
$\phi^{\rm b}: L^2(U) \to L^2(V)$ respectively defined by
\begin{equation}
\label{eq_definition_piola}
\phi^{\rm g}(v)
=
v \circ \phi^{-1},
\, \,
\phi^{\rm c}(\tw)
=
(\JAC^{-T} \tw) \circ \phi^{-1},
\, \,
\phi^{\rm d}(\tw)
=
\left (\frac{\JAC }{|\JAC|} \tw \right ) \circ \phi^{-1},
\, \,
\phi^{\rm b}(v)
=
\left (|\JAC| v \right ) \circ \phi^{-1},
\end{equation}
for} \tcf{all $v \in L^2(U)$ and $\tw \in \tL^2(V)$, where $\JAC$
is the Jacobian matrix of $\phi$ \cgr{and $|\JAC|$ its determinant}. These mappings are invertible. In addition,
if $\gamma_U \subset \partial U$ and $\gamma_V \eq \phi(V)$, \cgr{with a similar notation as in Section~\ref{sec_cont_sp_BC}, $\phi^{\rm g}$, $\phi^{\rm c}$, and $\phi^{\rm d}$}
\begin{equation}
\label{eq_piola_bijection}
H^1_{0,\gamma_U}\!(U) \! \to \! H^1_{0,\gamma_V}\!(V),
\tH_{0,\gamma_U}\!(\crl,U) \! \to \! \tH_{0,\gamma_V}\!(\crl,V),
\tH_{0,\gamma_U}\!(\dv,U) \! \to \! \tH_{0,\gamma_V}\!(\dv,V)
\end{equation}
are} \tcf{bijections, and more generally, the full, tangential\cgr{,} and normal traces
on $\gamma_U$ are respectively transported by $\phi^{\rm g}$, $\phi^{\rm c}$\cgr{,}
and $\phi^{\rm d}$ on $\gamma_V$. The commutativity properties
\begin{equation}
\label{eq_piola_commute}
\phi^{\rm c}(\Gr v) = \Gr (\phi^{\rm g} (v))
\quad
\Dv (\phi^{\rm d}(\tw)) = \phi^{\rm b}(\Dv \tw)
\end{equation}
for $v \in H^1(U)$ and $\tw \in \tH(\dv,U)$
will be useful. We will also need the formula
\begin{equation}
\label{eq_piola_green}
((\phi^{\rm d})^{-1}(\tv),\tw)_U = \varepsilon (\tv,\phi^{\rm c}(\tw))_V
\quad
\forall \tv \in \tL^2(V), \; \forall \tw \in \tL^2(U),
\end{equation}
where} \tcf{$\varepsilon$ is the (constant) sign of the determinant of $\JAC$.
Finally, if $U$ is \cgr{a polyhedron} covered by a tetrahedral mesh $\T_U$ and $\phi|_\elm$
is affine for each $\elm \in \T_U$, denoting by $\T_V$ the tetrahedral mesh
of $V$ induced by $\phi$, we have the bijections
\begin{equation}
\label{eq_piola_discrete_bijection}
\phi^{\rm g}: \PP_q(\T_U) \to \PP_q(\T_V), \, \,
\phi^{\rm c}: \ND_q(\T_U) \to \ND_q(\T_V), \, \,
\phi^{\rm d}: \RT_q(\T_U) \to \RT_q(\T_V)
\end{equation}
for all integer $q \geq 0$. In addition, in this case, we have
\begin{equation}
\label{eq_bound_piola}
\norm{\phi^{\rm d}}\norm{(\phi^{\rm d})^{-1}}
\leq
C(\kappa_{\T_U},\kappa_{\T_V})
\end{equation}
with $\norm{\phi^{\rm d}}$ denoting the norm operator of $\phi^{\rm d}: \tL^2(U) \to \tL^2(V)$
(and vice-versa for $\norm{(\phi^{\rm d})^{-1}}$), and $\kappa_{\T_U}$, $\kappa_{\T_V}$
the shape-regularity constants of $\T_U$ and $\T_V$ \cgr{as in Section~\ref{sec_mesh}}.}

\subsection{\cgr{A preliminary result}}

\tcf{Before proving Lemma~\ref{lem_tzab}, we establish the following preliminary result:}

\tcf{\bl[Orthogonalities]
\label{lem_W}
Let $\mean \in \tL^2(\oma)$ satisfy $(\mean,\Gr q_h)_{\oma} = 0$ for all $q_h \in \PP_1(\Ta) \cap \Hsa$.
Then, the following set is non-empty:
\[
W_{h}(\Ta,\mean)
\eq
\left \{
\tv_h \in \RT_1(\Ta) \cap \Hdva
\; \left | \;
\begin{array}{l}
\Dv \tv_h = 0
\\
(\tv_h, \tr_h)_{\elm} = (\mean, \tr_h)_{\elm}\\
\quad
\forall \tr_h \in [\PP_0(\elm)]^{3}, \, \forall \elm \in \Ta
\end{array}
\right .
\right \}.
\]
\el}

\begin{figure}
\centerline{\includegraphics[width=0.35\textwidth]{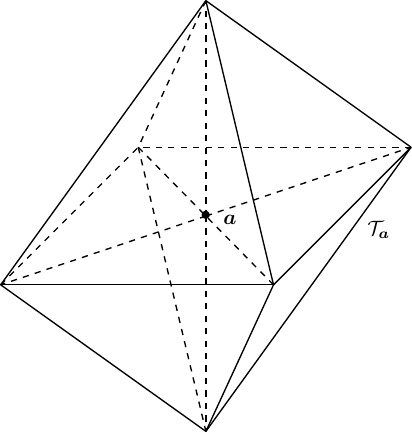} \qquad \includegraphics[width=0.25\textwidth]{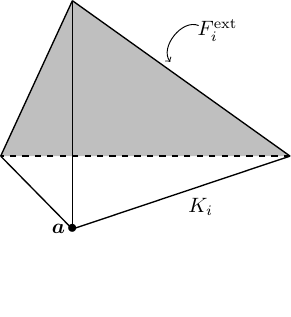}}
\caption{Interior vertex patch (left) and the element $\elm_1$ with the face $\sd_1^{\mathrm{ext}}$ (right)}
\label{fig_patch}
\end{figure}

\bp
{\em Step~1: interior patches.}
We start with the case where the vertex $\ver \in \Vh$
does not lie on the boundary $\pt \Om$, \cf~Figure~\ref{fig_patch}, left.
\tcf{We will construct a particular element $\tw_h \in W_h(\Ta,\mean)$} by an explicit
run through the patch $\Ta$ of tetrahedra sharing the vertex $\ver \in \Vh$,
similarly as in~\cite{Brae_Pill_Sch_p_rob_09, Ern_Voh_p_rob_3D_20}.
Specifically, following the concept of shelling of a polytopal complex,
see~\cite[Theorem~8.12]{Zieg_poly_95} and~\cite[Lemma~B.1]{Ern_Voh_p_rob_3D_20},
there exists an enumeration $\elmi$,
$1 \leq i \leq |\Ta|$, of the tetrahedra in the patch $\Ta$ such that, except for
the first tetrahedron in the enumeration $\elmo$: (i) if there are at least two faces
corresponding to the neighbors of $\elmi$ which have been already enumerated, then all
the tetrahedra of $\Ta$ sharing this edge come sooner in the enumeration; (ii) except
for the last element $\elm_{|\Ta|}$, there are one or two neighbors of $\elmi$ which have
been already enumerated and correspondingly two or one neighbors of $\elmi$ which have not
been enumerated yet.

Consider a pass through the patch $\Ta$ in the sense of the above enumeration.
For the tetrahedron $\elmi$, $1 \leq i \leq |\Ta|$, let us denote by $\F_i^\sharp$ the faces
of $\elmi$ corresponding to the neighbors of $\elmi$ which have been already passed through
and $\sd^j = \pt \elmi \cap \pt \elm_j \in \F_i^\sharp$ the face corresponding to the neighbor
$\elm_j$. Also, let $\sd_i^{\mathrm{ext}}$ be the face of $\elmi$ lying on the patch boundary
$\pt \oma$. Consider the problem
\be \label{eq_tzi}
\tw_h^i
\eq
\arg \min_{\substack{
\tv_h \in \RT_1(\elmi)
\\
\Dv \tv_h= 0
\\
\tv_h \scp \tn_\elmi = 0 \text{ on } \sd_i^{\mathrm{ext}}
\\
\tv_h \scp \tn_\elmi = \tw_h^j \scp \tn_\elmi \text{ on all } \sd^j \in \F_i^\sharp
\\
(\tv_h, \tr_h)_\elmi = (\tcf{\mean}, \tr_h)_\elmi
\quad
\forall \tr_h \in [\PP_0(\elmi)]^3
}}
\norm{\tv_h\tcf{-\mean}}_\elmi^2,
\ee
similar to~\eqref{eq_tzab} but reduced to the single tetrahedron $\elmi$. If a solution
to~\eqref{eq_tzi} exists, $\tw_h$ defined as $\tw_h|_\elmi \eq \tw_h^i$
\tcf{belongs to $W_h(\Ta,\mean)$}. We are thus left to establish the existence and uniqueness
of~\eqref{eq_tzi}.

{\em Step~1a: the first element $\elmo$.}
Let us start with the first element $\elmo$, \cf~Figure~\ref{fig_patch}, right.
Then the set $\F_1^\sharp$ is empty, and \cye{the constraints in~\eqref{eq_tzi} lead us to} ask whether in the first-order Raviart--Thomas
space $\RT_1(\elmo)$, one can impose simultaneously
the divergence, the normal flux through one face, and moments against constant
functions. We will reason by the canonical degrees of freedom, see, \eg,
\cite[Proposition~2.3.4 and Figure~2.14.c]{Bof_Brez_For_MFEs_13} \cye{or~\eqref{eq_RT_proj}}, and find a suitable
$\tv_h \in \RT_1(\elmo)$ \cye{satisfying the constraints in~\eqref{eq_tzi}}. First, we see that in $\RT_1(\elmo)$, the normal flux
$\tv_h \scp \tn_\elmo$ on $\sd_1^{\mathrm{ext}}$ can be fixed to zero and the moments
against constants $(\tv_h, \tr_h)_\elmo$ can be fixed as in~\eqref{eq_tzi}. We still have
the freedom to choose the normal fluxes $\tv_h \scp \tn_\elmo$ on the faces of $\elmo$ different
from $\sd_1^{\mathrm{ext}}$, and the question is whether this can be done so as to fix the
divergence of $\tv_h$ to zero. By~\cite[Proposition~2.3.3]{Bof_Brez_For_MFEs_13}, there holds
\[
\Dv \tv_h = 0
\quad \Leftrightarrow \quad
(\Dv \tv_h, q_h)_\elmo = 0 \qquad \forall q_h \in \PP_1(\elmo).
\]
Employing the Green theorem and the fact that $\tv_h \scp \tn_\elmo = 0$ on $\sd_1^{\mathrm{ext}}$,
\[
(\Dv \tv_h, q_h)_\elmo
=
\<\tv_h \scp \tn_\elmo, q_h\>_{\pt \elmo \setminus \sd_1^{\mathrm{ext}}}
-
(\tv_h, \Gr q_h)_\elmo.
\]
Now, since $\Gr q_h \in [\PP_0(\elmo)]^3$, the last term above is fixed from the last constraint
in~\eqref{eq_tzi}, so the question becomes: can one choose $\tv_h \scp \tn_\elmo$ on
$\pt \elmo \setminus \sd_1^{\mathrm{ext}}$ such that
\be \label{eq_4_9}
\<\tv_h \scp \tn_\elmo, q_h\>_{\pt \elmo \setminus \sd_1^{\mathrm{ext}}}
=
(\tcf{\mean}, \Gr q_h)_\elmo \qquad \forall q_h \in \PP_1(\elmo),
\ee
which gives $4$ conditions for the $9$ remaining degrees of freedom (there are \cye{$4$ degrees of freedom in $\PP_1(\elmo)$ and} $3$ degrees of freedom per face in $\RT_1(\elmo)$ following~\cite[Proposition~2.3.4]{Bof_Brez_For_MFEs_13}).

We proceed as follows. Out of the three faces of $\elmo$ different from
$\sd_1^{\mathrm{ext}}$, choose one and impose $\tv_h \scp \tn_\elmo = 0$ therein.
Then we are left to set $\tv_h \scp \tn_\elmo$ on two faces, say $\sd$ and $\widetilde \sd$.
For $\sd$, consider the three hat basis functions $\psi_\sd^k$, $1 \leq k \leq 3$, as in
Section~\ref{sec_mesh}, corresponding to its three vertices. Restricted to $\widetilde \sd$,
which is necessary a face neighboring with $\sd$, they belong to $\PP_1(\widetilde \sd)$, and
one of the restrictions, say $\psi_\sd^3$, is zero on $\widetilde \sd$. Thus, there holds \ifSIAM$\<\tv_h \scp \tn_\elmo, \psi_\sd^3\>_{\widetilde \sd} = 0,$ \else
\[
\<\tv_h \scp \tn_\elmo, \psi_\sd^3\>_{\widetilde \sd} = 0,
\]\fi
and, following~\cite[Proposition~2.3.4]{Bof_Brez_For_MFEs_13}, we can prescribe
\[
\<\tv_h \scp \tn_\elmo, \psi_\sd^k\>_{\widetilde \sd} \eq 0 \qquad 1 \leq k \leq 2.
\]
Moreover, restricted to $\sd$, $\psi_\sd^k$ create a basis of $\PP_1(\sd)$, whereas restricted
to $\elmo$, they belong to $\PP_1(\elmo)$. Thus, we can also set
\[
\<\tv_h \scp \tn_\elmo, \psi_\sd^k\>_\sd
\eq
(\tcf{\mean}, \Gr \psi_\sd^k)_\elmo \qquad 1 \leq k \leq 3.
\]
With the choices made so far, we see that~\eqref{eq_4_9} holds for the three hat functions
$\psi_\sd^k$, $1 \leq k \leq 3$. Finally, consider $\psi_\sd^4$, the hat basis function
corresponding to the vertex opposite to the face $\sd$. Restricted to $\sd$, it is zero,
so that
\ifSIAM$\<\tv_h \scp \tn_\elmo, \psi_\sd^4\>_\sd = 0.$ \else
\[
\<\tv_h \scp \tn_\elmo, \psi_\sd^4\>_\sd = 0.
\]\fi
Moreover, restricted to $\widetilde \sd$, it completes $\psi_\sd^1$ and $\psi_\sd^2$
(restricted to $\widetilde \sd$) to create a basis of $\PP_1(\widetilde \sd)$, and restricted
to $\elmo$, it belongs to $\PP_1(\elmo)$, so that we can prescribe
\[
\<\tv_h \scp \tn_\elmo, \psi_\sd^4\>_{\widetilde \sd} \eq (\tcf{\mean}, \Gr \psi_\sd^4)_\elmo.
\]
Thus, \eqref{eq_4_9} also holds for $\psi_\sd^4$, and since $\psi_\sd^k$, $1 \leq k \leq 4$,
restricted to $\elmo$ create a basis of $\PP_1(\elmo)$, \eqref{eq_4_9} holds true, and a
unique $\tw_h^1$ from~\eqref{eq_tzi} exists.

{\em Step~1b: any element $\elmi$ with $|\F_i^\sharp|=1$.}
We now investigate those consecutive elements $\elmi$ which are such that two neighbors of
$\elmi$ have not been passed through yet. This means that exactly one neighbor of $\elmi$,
say $\elmj$, has already been passed through, so there is one face $\sd^j$ in the set
$\F_i^\sharp$. Since $\tv_h \scp \tn_\elmi = \tw_h^j \scp \tn_\elmi$ on $\sd^j$ is requested
in~\eqref{eq_tzi}, \eqref{eq_4_9} asks if can one choose $\tv_h \scp \tn_\elmi$
on $\pt \elmi \setminus \{\sd_i^{\mathrm{ext}}, \sd^j\}$ such that
\be \label{eq_4_6}
\<\tv_h \scp \tn_\elmi, q_h\>_{\pt \elmi \setminus \{\sd_i^{\mathrm{ext}}, \sd^j\}}
=
(\tcf{\mean}, \Gr q_h)_\elmi
-
\<\tw_h^j \scp \tn_\elmi, q_h\>_{\sd^j}
\ee
for all $q_h \in \PP_1(\elmi)$, which is still undetermined, giving $4$ conditions for the $6$ remaining degrees of freedom.
The reasoning is similar as for $\elmo$. Still denoting $\sd$ and $\widetilde \sd$ the two
remaining faces and $\psi_\sd^k$, $1 \leq k \leq 4$ the hat basis functions, we again have
\[
\<\tv_h \scp \tn_\elmi, \psi_\sd^3\>_{\widetilde \sd} = 0,
\qquad
\<\tv_h \scp \tn_\elmi, \psi_\sd^4\>_\sd = 0.
\]
Moreover, imposing
\ban
\<\tv_h \scp \tn_\elmi, \psi_\sd^k\>_{\widetilde \sd}
&\eq
0 \qquad 1 \leq k \leq 2,
\\
\<\tv_h \scp \tn_\elmi, \psi_\sd^4\>_{\widetilde \sd}
&\eq
(\tcf{\mean}, \Gr \psi_\sd^4)_\elmi
-
\<\tw_h^j \scp \tn_\elmi, \psi_\sd^4\>_{\sd^j},
\\
\<\tv_h \scp \tn_\elmi, \psi_\sd^k\>_\sd
&\eq
(\tcf{\mean}, \Gr \psi_\sd^k)_\elmi
-
\<\tw_h^j \scp \tn_\elmi, \psi_\sd^k\>_{\sd^j} \qquad 1 \leq k \leq 3
\ean
yields~\eqref{eq_4_6}, and $\tw_h^i$ exists.

{\em Step~1c: any element $\elmi$ with $|\F_i^\sharp|=2$.}
We now investigate those consecutive elements $\elmi$ which are such that only
one neighbor of $\elmi$ has not been passed through yet, with $\elm_{j_1}$ and
$\elm_{j_2}$ already passed through and faces $\sd^{j_1}$, $\sd^{j_2}$ in the set
$\F_i^\sharp$. Denote $\sd$ the only remaining face, so that $\sd_i^{\mathrm{ext}}$,
$\sd^{j_1}$, $\sd^{j_2}$, and $\sd$ are the four faces of the tetrahedron $\elmi$.
As in~\eqref{eq_4_9} and~\eqref{eq_4_6}, we need to ensure that
\be \label{eq_4_3}
\<\tv_h \scp \tn_\elmi, q_h\>_\sd
=
(\tcf{\mean}, \Gr q_h)_\elmi
-
\<\tw_h^{j_1} \scp \tn_\elmi, q_h\>_{\sd^{j_1}}
-
\<\tw_h^{j_2} \scp \tn_\elmi, q_h\>_{\sd^{j_2}}
\ee
for all $q_h \in \PP_1(\elmo)$.
This time, the system is over-determined in that we request $4$ conditions for the $3$
remaining degrees of freedom of the normal components $\tv_h \scp \tn_\elmi$ on the face
$\sd$. As above, we can impose%
\[
\<\tv_h \scp \tn_\elmi, \psi_\sd^k\>_\sd
\eq
(\tcf{\mean}, \Gr \psi_\sd^k)_\elmi
-
\<\tw_h^{j_1} \scp \tn_\elmi, \psi_\sd^k\>_{\sd^{j_1}}
-
\<\tw_h^{j_2} \scp \tn_\elmi, \psi_\sd^k\>_{\sd^{j_2}}
\quad 1 \leq k \leq 3,
\]
which fixes $\tv_h \scp \tn_\elmi$ on the face $\sd$. Now, noting that $\<\tv_h \scp \tn_\elmi, \psi_\sd^4\>_\sd = 0$, it follows that to prove~\eqref{eq_4_3}, we need to show that
\be
\label{eq_Ki_2}
(\tcf{\mean}, \Gr \psi_\sd^4)_\elmi
-
\<\tw_h^{j_1} \scp \tn_\elmi, \psi_\sd^4\>_{\sd^{j_1}}
-
\<\tw_h^{j_2} \scp \tn_\elmi, \psi_\sd^4\>_{\sd^{j_2}}
=
0.
\ee

To prove~\eqref{eq_Ki_2}, recall from property (i) of the enumeration (giving that all other
elements sharing the edge $\edg$ common to $\sd^{j_1}$ and $\sd^{j_2}$ have been already passed
through) and the previous steps, see~\eqref{eq_4_9} and~\eqref{eq_4_6}, that
\be \label{eq_Kj}
\<\tw_h^j \scp \tn_{\elmj}, \psi_\sd^4\>_{\pt \elmj}
=
(\tcf{\mean}, \Gr \psi_\sd^4)_{\elmj}
\ee
for all the tetrahedra $\elmj$ sharing the edge $\edg$, different from $\elmi$.
Recalling \cye{from assumptions of Lemma~\ref{lem_W}}, we have
\be \label{eq_orth_ome}
0
=
(\tcf{\mean}, \Gr \psi_\sd^4)_\oma
=
(\tcf{\mean}, \Gr \psi_\sd^4)_\ome,
\ee
where $\ome$ is the part of $\oma$ corresponding to the elements sharing the edge $\edg$;
the last equality holds since in the vertex patch subdomain $\oma$,
$\psi_\sd^4$ is only supported on the edge patch subdomain $\ome$. Denote by $\tome$
the part of $\ome$ without the element $\elmi$. Then the normal traces
orientation, the zero normal trace boundary conditions
$\tw_h^j \scp \tn_{\elmj} = 0$ on the faces $\sd_j^{\mathrm{ext}}$ together with the
zero values of $\psi_\sd^4$ \cye{on $\pt \ome \setminus \pt \oma$}, the Green theorem first applied on $\tome$ and later individually on $\elmj$,
\cye{and the notation $\tw_h|_{\elmj} = \tw_h^j$ for the previous $\elmj^\circ \subset \tome$} give
\be\label{eq_IPPs}
\bs
{}&-\<\tw_h^{j_1} \scp \tn_\elmi, \psi_\sd^4\>_{\sd^{j_1}}
-
\<\tw_h^{j_2} \scp \tn_\elmi, \psi_\sd^4\>_{\sd^{j_2}}
=
\<\tw_h^{j_1} \scp \tn_\tome, \psi_\sd^4\>_{\sd^{j_1}}
+
\<\tw_h^{j_2} \scp \tn_\tome, \psi_\sd^4\>_{\sd^{j_2}}
\\
= {} & \cye{\<\tw_h \scp \tn_\tome, \psi_\sd^4\>_{\pt \tome} = (\tw_h, \Gr \psi_\sd^4)_{\tome}
+ (\Dv \tw_h, \psi_\sd^4)_{\tome}} \\
= {} &
\sum_{j; \, \elmj^\circ \subset \tome} \big\{(\tw_h^j, \Gr \psi_\sd^4)_{\elmj}
+
(\Dv \tw_h^j, \psi_\sd^4)_{\elmj}\big\}
=
\sum_{j; \, \elmj^\circ \subset \tome}  \<\tw_h^j \scp \tn_{\elmj}, \psi_\sd^4\>_{\pt \elmj}
\\
\reff{eq_Kj}{=} \!\!\!\! {} &
\sum_{j; \, \elmj^\circ \subset \tome} (\tcf{\mean}, \Gr \psi_\sd^4)_{\elmj} \cye{=(\tcf{\mean}, \Gr \psi_\sd^4)_\ome -(\tcf{\mean}, \Gr \psi_\sd^4)_\elmi}
\reff{eq_orth_ome}{=} -(\tcf{\mean}, \Gr \psi_\sd^4)_\elmi,
\es \ee
which is~\eqref{eq_Ki_2}. Thus, there exists a unique $\tw_h^i$ from~\eqref{eq_tzi}
also on this $\elmi$.

{\em Step~1d: the last element $\elml$.}
According to property (ii) of the enumeration, the last element $\elml$ is such that
$|\F_{|\Ta|}^\sharp|=3$, so that all the neighbors have been already passed through.
Consequently, all the degrees of freedom of $\tv_h$ are fixed from the last three constraints
in~\eqref{eq_tzi}, and we need to show that $\Dv \tv_h = 0$, \ie, that
\[
(\Dv \tv_h, q_h)_\elml = 0 \qquad \forall q_h \in \PP_1(\elml),
\]
since $\Dv \tv_h \in \PP_1(\elml)$.
\cye{Using the Green theorem and the constraints in~\eqref{eq_tzi} as above, t}his is equivalent to verifying that
\be \label{eq_4_0} \bs
0 = {} &
(\tcf{\mean}, \Gr \psi_\sd^k)_\elml
-
\<\tw_h^{j_1} \scp \tn_\elml, \psi_\sd^k\>_{\sd^{j_1}} \\
& -
\<\tw_h^{j_2} \scp \tn_\elml, \psi_\sd^k\>_{\sd^{j_2}}
-
\<\tw_h^{j_3} \scp \tn_\elml, \psi_\sd^k\>_{\sd^{j_3}}
\es\ee
for all $1 \leq k \leq 4$, where $\sd^{j_1}$, $\sd^{j_2}$, $\sd^{j_3}$ are
the three faces in $\F_{|\Ta|}^\sharp$ and $\psi_\sd^k$ are the hat basis
functions associated with the four vertices of $\elml$. As in~\eqref{eq_orth_ome},
\cye{assumptions of Lemma~\ref{lem_W}} imply
\be \label{eq_orth_oma}
    0 = (\tcf{\mean}, \Gr \psi_\sd^k)_\oma \qquad 1 \leq k \leq 4.
\ee
Moreover, as in~\eqref{eq_Kj}, \cye{it follows from~\eqref{eq_4_9}, \eqref{eq_4_6}, and~\eqref{eq_4_3}}
\be \label{eq_Kjt}
\<\tw_h^j \scp \tn_{\elmj}, \psi_\sd^k\>_{\pt \elmj}
=
(\tcf{\mean}, \Gr \psi_\sd^k)_{\elmj}
\qquad 1 \leq k \leq 4
\ee
is satisfied on all elements $\elmj$ of the patch $\Ta$ other than $\elml$.
Let $\toma$ correspond to the patch subdomain $\oma$ without the element $\elml$.
Then, as in~\eqref{eq_IPPs},
\ban
{}&
-\<\tw_h^{j_1} \scp \tn_\elml, \psi_\sd^k\>_{\sd^{j_1}}
-\<\tw_h^{j_2} \scp \tn_\elml, \psi_\sd^k\>_{\sd^{j_2}}
-\<\tw_h^{j_3} \scp \tn_\elml, \psi_\sd^k\>_{\sd^{j_3}}
\\
={}&
\<\tw_h^{j_1} \scp \tn_\toma, \psi_\sd^k\>_{\sd^{j_1}}
+
\<\tw_h^{j_2} \scp \tn_\toma, \psi_\sd^k\>_{\sd^{j_2}}
+
\<\tw_h^{j_2} \scp \tn_\toma, \psi_\sd^k\>_{\sd^{j_3}}
\\
= {} & \cye{\<\tw_h \scp \tn_\toma, \psi_\sd^k\>_{\pt \toma} = (\tw_h, \Gr \psi_\sd^k)_{\toma}
+ (\Dv \tw_h, \psi_\sd^k)_{\toma}} \\
={}&
\sum_{j; \, \elmj^\circ \subset \toma} \big\{(\tw_h^j, \Gr \psi_\sd^k)_{\elmj}
+
(\Dv \tw_h^j, \psi_\sd^k)_{\elmj}\big\}
=
\sum_{j; \, \elmj^\circ \subset \toma}  \<\tw_h^j \scp \tn_{\elmj}, \psi_\sd^k\>_{\pt \elmj}
\\
\reff{eq_Kjt}{=} \!\!\!\! {} &
\sum_{j; \, \elmj^\circ \subset \toma} (\tcf{\mean}, \Gr \psi_\sd^k)_{\elmj}
\cye{=(\tcf{\mean}, \Gr \psi_\sd^k)_\oma -(\tcf{\mean}, \Gr \psi_\sd^k)_\elml}
\reff{eq_orth_oma}{=} - (\tcf{\mean}, \Gr \psi_\sd^k)_\elml
\ean
for all $1 \leq k \leq 4$, \ie, \eqref{eq_4_0}.
Thus, there exists a minimizer $\tw_h^{|\Ta|}$ of~\eqref{eq_tzi} on $\elml$.

\tcf{{\em Step~2. Boundary patches with flat boundaries.}
We now investigate the case where the vertex $\ver \in \Vh$ lies on the boundary $\pt \Om$.
We present in this step in detail the case of a boundary patch $\Ta$
for which $\Gamma_\ver$, the part of $\partial \oma$ that contains the faces sharing the vertex $\ver$, is contained in a plane $H$, which we call a ``flat boundary'' case.
For the sake of simplicity, assume that either
$\Gamma_\ver \subset \GD$ or $\Gamma_\ver \subset \GN$ and, without loss generality,
that $H = \{ \tx \in \mathbb R^3; \; \tx_3 = 0 \}$. The symmetrization operator around
the plane $H$, $\phi: \tx =(\tx_1,\tx_2,\tx_3) \to (\tx_1,\tx_2,-\tx_3)$, \cgr{as in~\cite[Section~7]{Ern_Voh_p_rob_3D_20} and~\cite[Section~7]{Chaum_Voh_p_rob_3D_H_curl_22},}
will be instrumental in the proof.
Specifically, we introduce the symmetrized patch $\tTa \eq \Ta \cup \phi(\Ta)$,
with the associated domain $\toma$, obtained by mapping the elements of $\Ta$ by $\phi$. We will employ the Piola mappings \cgr{from~\eqref{eq_definition_piola}} to relate the set $W_h(\Ta,\mean)$ from the announcement of Lemma~\ref{lem_W} to a set $W_h(\tTa,\tmean)$ with an extended datum $\tmean$. Then, the result will follow by {\em Step~1}, since $\tTa$ is an interior patch.}

\tcf{{\em Step~2a. The case $\ver \in \GD$}. 
We start by defining the extended datum $\tmean \in \tL^2(\toma)$ from $\mean$:
we simply set $\tmean \eq \mean$ in $\oma$ and $\tmean \eq \phi^{\rm d}(\mean)$
on $\toma \setminus \oma$. Let $\widetilde q_h \in \PP_1(\tTa) \cap \Hsaw$.
Recalling that $\phi$ is a symmetrization operator, its \cgr{(constant)} Jacobian matrix has
a negative determinant. As a result, using~\eqref{eq_piola_green} and~\eqref{eq_piola_commute},}
\tcf{\begin{align*}
(\tmean,\Gr \widetilde q_h)_{\toma}
\! &= \!
(\mean,\Gr \widetilde q_h)_{\oma}
\!\! + \!
(\phi^{\rm d}(\mean),\Gr \widetilde q_h)_{\toma \setminus \oma}
\!\! = \!
(\mean,\Gr \widetilde q_h)_{\oma}
\!\! - \!
(\mean,(\phi^{\rm c})^{-1}(\Gr \widetilde q_h))_{\oma}
\\
\! &= \!
(\mean,\Gr \widetilde q_h)_{\oma}
\!\! - \!
(\mean,\Gr \cgr{(}(\phi^{\rm g})^{-1}(\widetilde q_h))_{\oma}
\!\! = \!
(\mean,\Gr q_h)_{\oma}
\end{align*}
with} \tcf{$q_h \eq \widetilde q_h-(\phi^{\rm g})^{-1}(\widetilde q_h)$. Because $(\phi^{\rm g})^{-1}$
preserves the trace on $H$, we see that $q_h = 0$ on $H$
(see~\eqref{eq_piola_bijection}), and since it also maps piecewise
polynomials to piecewise polynomials (see~\eqref{eq_piola_discrete_bijection}),
$q_h \in \PP_1(\Ta) \cap \Hsa$ (recall from Section~\ref{sec_cont_sp_oma} that $\Hsa = \{v \in \Hoi{\oma}; \, v = 0$ on $\gD = \Gamma_\ver\}$ here). Hence, $(\tmean,\Gr \widetilde q_h)_{\toma} = 0$
by our assumption $(\mean,\Gr q_h)_{\oma} = 0$. Thus}
\tcf{$\tmean$ satisfies the assumption of
Lemma~\ref{lem_W} on the interior patch $\tTa$, and therefore {\em Step~1} ensures that the set $W_h(\tTa,\tmean)$ is non-empty.}

\tcf{We now consider an arbitrary element $\widetilde \tw_h \in W_h(\tTa,\tmean)$ and set $\tw_h \eq \widetilde \tw_h|_{\oma}$. Since $\widetilde \tw_h \in \RT_1(\tTa) \cap \Hdvat$, it is clear that $\tw_h \in \RT_1(\Ta) \cap \Hdva$, namely as no normal trace boundary conditions are required on $\Gamma_\ver \subset H$. Indeed, in this case, $\Gamma_\ver = \gD$ in the notation of Section~\ref{sec_cont_sp_oma}, so that
$\Hdva = \{ \tv \in \tH(\dv,\oma); \; \tv {\cdot} \tn_{\oma} = 0 \text{ on }
\partial \oma \setminus \Gamma_\ver\}$.
Moreover, $\Dv \tw_h = \Dv \widetilde \tw_h = 0$
on $\oma$. Finally, $(\tw_{h}, \tr_h)_{\elm} = (\mean, \tr_h)_{\elm}$ for all $\tr_h \in [\PP_0(\elm)]^{3}$ and all $\elm \in \Ta$ since $\Ta \subset \tTa$ and simply $(\tw_{h}, \tr_h)_{\elm} = (\widetilde \tw_{h}, \tr_h)_{\elm}$, $(\mean, \tr_h)_{\elm} = (\tmean, \tr_h)_{\elm}$, and $\widetilde \tw_h \in W_h(\tTa,\tmean)$, so that $(\widetilde \tw_{h}, \tr_h)_{\elm} = (\tmean, \tr_h)_{\elm}$. This concludes the proof that the set $W_h(\Ta,\mean)$ is non-empty in this case.}

\tcf{{\em Step~2b. The case $\ver \in \GN$}. 
In this case, we extend the datum $\mean$ by setting
$\tmean \eq \mean$ on $\oma$ and $\tmean \eq \boldsymbol 0$ on
$\toma \setminus \oma$. If $\widetilde q_h \in \PP_1(\tTa) \cap
\Hsaw$, we have
\[
(\tmean,\Gr\widetilde q_h)_{\toma}
=
(\mean,\Gr\widetilde q_h)_{\oma}
=
0
\]
since $\widetilde q_h|_{\oma} \in \PP_1(\Ta) \cap H^1(\oma)$, whose gradients have the same span as those of $\PP_1(\Ta) \cap \Hsa$, the zero mean value subspace of $\PP_1(\Ta) \cap H^1(\oma)$ following Section~\ref{sec_cont_sp_oma} in this case. It thus follows
from {\em Step~1} that $W_h(\tTa,\tmean)$ is non-empty.

Consider an element $\widetilde \tw_h \in W_h(\tTa,\tmean)$ and set $\tw_h \eq \widetilde \tw_h|_\oma -(\phi^{\rm d})^{-1}(\widetilde \tw_h|_{\toma \setminus \oma})$. We need to show that $\tw_h \in W_h(\Ta,\mean)$. Recall that here the functions in $\Hdva$ need to satisfy the no-flow boundary condition on the whole patch boundary $\partial \oma$ and in particular on $H$: $\Hdva = \{\tv \in \Hdvi{\oma}; \, \tv \scp \tn_{\oma}=0$ on $\pt \oma\}$ from Section~\ref{sec_cont_sp_oma} in this case.
Since the Piola mapping $(\phi^{\rm d})^{-1}$ maps piecewise
Raviart--Thomas polynomials to piecewise Raviart--Thomas polynomials
(\cf~\eqref{eq_piola_discrete_bijection}) and
preserves the divergence (\cf~\eqref{eq_piola_commute}) and the normal trace
(\cf~\eqref{eq_piola_bijection}),
it is clear that $\tw_h \in \RT_1(\Ta) \cap \Hdva$ and $\Dv \tw_h = 0$.
It remains to show that $(\tw_{h}, \tr_h)_{\elm} = (\mean, \tr_h)_{\elm}$ for all
$\tr_h \in [\PP_0(\elm)]^{3}$ and all $\elm \in \Ta$. Let $\elm \in \Ta$ and
$\tr_h \in [\PP_0(\elm)]^{3}$ and let $\widetilde \elm$ be the tetrahedron corresponding to
$\elm$ by the symmetry map $\phi$. Then}
\tcf{
\begin{align*}
(\tw_h,\tr_h)_\elm
&=
(\widetilde \tw_h,\tr_h)_\elm
-
((\phi^{\rm d})^{-1}(\widetilde \tw_h),\tr_h)_\elm
=
(\widetilde \tw_h,\tr_h)_\elm
+
(\widetilde \tw_h,\phi^{\rm c}(\tr_h))_{\widetilde \elm}\\
&=(\tmean,\tr_h)_\elm
+
(\tmean,\phi^{\rm c}(\tr_h))_{\widetilde \elm}
=
(\tmean,\tr_h)_\elm = (\mean,\tr_h)_\elm,
\end{align*}
where} \tcf{we have used~\cred{\eqref{eq_piola_green},} that the Piola mapping $\phi^{\rm c}$ maps piecewise
constant vectors onto \cred{piecewise} constant vectors (this can be seen from the definition
\eqref{eq_definition_piola} of $\phi^{\rm c}$ since its Jacobian matrix is \cred{constant here}),
that $\widetilde \tw_h \in W_h(\tTa,\tmean)$,
and finally that $\tmean$ is the extension of $\mean$ by zero
to the symmetrized patch.
This concludes the proof that $W_h(\Ta,\mean)$ is non-empty in this case.}

\tcf{{\em Step~3: general boundary patches.} For general boundary patches, the proof follows the lines of {\em Step~2} while employing the extension and restriction operators introduced in~\cite[Section 7]{Chaum_Voh_p_rob_3D_H_curl_22} instead of the (simpler) symmetrization
operator $\phi$ of {\em Step~2}. We do not give details here.}
\ep

\subsection{\cgr{Proof of Lemma~\ref{lem_tzab}}} \cgr{We can now finally establish a proof of Lemma~\ref{lem_tzab}.}

\enlargethispage{\baselineskip}

\bp[Proof of Lemma~\ref{lem_tzab}]
\tcf{
{\em Step 1. Existence and uniqueness.}
The minimization set in~\eqref{eq_tzab} is the set $W_h(\Ta,\mean)$ of Lemma~\ref{lem_W}
with $\mean \eq \btau_h^\ver-\thab$. Since the minimization functional in~\eqref{eq_tzab} is convex, it is sufficient
to show that $W_h(\Ta,\mean)$ is non-empty to ensure the existence and uniqueness of $\tzab$ from~\eqref{eq_tzab}.
From Lemma~\ref{lem_W}, we need to show that $(\mean,\Gr q_h)_{\oma} = 0$ for all
$q_h \in \PP_1(\Ta) \cap \Hsa$. This is actually a direct consequence
of assumption~\eqref{eq_orth_A}. Indeed, from the divergence constraint in~\eqref{eq_tsab} and since $\qq \geq 1$ and $q_h \in \PP_1(\Ta) \cred{\cap \Hsa}$, we have
\ban
(\mean,\Gr q_h)_{\oma}
& =
(\btau_h^\ver,\Gr q_h)_{\oma}
-
(\thab,\Gr q_h)_{\oma}
=
(\btau_h^\ver,\Gr q_h)_{\oma}
+
(\Dv \thab,q_h)_{\oma}
\\
& =
(\btau_h^\ver,\Gr q_h)_{\oma}
+
(\Pi_\qq(g^\ver),q_h)_{\oma}
=
(\btau_h^\ver,\Gr q_h)_{\oma}
+
(g^\ver,q_h)_{\oma}
\reff{eq_orth_A}{=}
0.
\ean
}

{\em Step~2. Stability bound.}
We now proceed with the proof of the stability~\eqref{eq_tzab_stab}.

{\em Step~2a: generic stability bound.}
Set again $\tcf{\mean} \eq \btau_h^\ver-\thab$, and denote by $\tcf{\mean_h}$
the $\tLti{\oma}$-orthogonal projection of $\tcf{\mean}$ onto $[\PP_1(\Ta)]^3$.
Considering the Euler(--Lagrange) equations associated with~\eqref{eq_tzab}, it is clear that
we can equivalently replace $\tcf{\mean}$ by $\tcf{\mean_h}$ in the definition~\eqref{eq_tzab} of
$\tzab$. Furthermore, because~\eqref{eq_tzab} is a quadratic minimization problem with
linear constraints, the operator $T: [\PP_1(\Ta)]^3 \ni \tcf{\mean_h} \to \tzab \in \RT_1(\Ta) \cap \Hdva$ (well-defined from {\em Step~1}) is linear. Since both $[\PP_1(\Ta)]^3$ and
$\RT_1(\Ta) \cap \Hdva$ are finite-dimensional spaces, the operator $T$ is continuous,
and there exists a constant $C(\Ta)$ such that
\be \label{tmp_stab_tzab}
\norm{\tzab}_{\oma}
\leq
C(\Ta)\norm{\tcf{\mean_h}}_{\oma}
\leq
C(\Ta)\norm{\tcf{\mean}}_{\oma},
\ee
where we used the fact that $\tcf{\mean_h}$ is defined from $\tcf{\mean}$
by projection in the last inequality. The constant $C(\Ta)$ is independent of the
polynomial degree $q$ \cye{(recall that~\eqref{eq_tzab} works with $\RT_1$ elements only)}
but depends on the patch $\Ta$ in an unspecified way.
To make the dependence explicit, we resort in the next step to a reference patch.

{\em Step~2b: explicit stability bound.}
For a fixed shape-regularity parameter $\kappa_{\Th}$ from Section~\ref{sec_mesh},
there exists a maximal number of elements $N(\kappa_{\Th})$ allowed in any patch $\Ta$. In turn,
for any $N(\kappa_{\Th})$, there exists a finite set of reference patches
$\{\widehat \T\}$ such that for all vertex patches $\Ta$, there exists a reference
patch $\widehat \T$ and a bilipschitz mapping $\phi: \oma \to \widehat \omega$
($\widehat \omega$ being the open domain associated with $\widehat \T$) such that
$\phi|_\elm$ is an affine mapping between the tetrahedron $\elm \in \Ta$ and a tetrahedron
$\widehat \elm \in \widehat \T$. The associated Piola mapping $\phi^{\rm d}$
from~\eqref{eq_definition_piola} will be useful.

Crucially, we observe that for all $\widehat \elm \in \widehat \T$,
$\tv \in \tLti{\elm}$, and $\widehat \tr_h \in [\PP_0(\widehat \elm)]^3$,
there exists $\tr_h \in [\PP_0(\elm)]^3$ such that
$(\phi^{\rm d}(\tv),\widehat \tr_h)_{\widehat \elm} = (\tv,\tr_h)_\elm$,
since, elementwise, the Piola transform amounts to a multiplication
by a constant matrix and a change of coordinates.
It follows that $\phi^{\rm d}$ maps the minimization set of~\eqref{eq_tzab}
on $\Ta$ into the minimization set of the equivalent problem set on $\widehat \T$ with
constraints $\phi^{\rm d}(\tcf{\mean})$.

Now, on the reference patch $\widehat \T$, if $\htzab_h$ is the minimizer
of~\eqref{eq_tzab} with the datum $\phi^{\rm d}(\tcf{\mean})$, we conclude from
{\em Step~2a} that
\[
\norm{\htzab_h}_{\widehat \omega}
\leq
C(\kappa_{\Th})
\norm{\phi^{\rm d}(\tcf{\mean})}_{\widehat \omega}
\leq
C(\kappa_{\Th})
\norm{\phi^{\rm d}}\norm{\tcf{\mean}}_{\oma}.
\]
On the other hand, since
$(\phi^{\rm d})^{-1}(\htzab_h)$ belongs the minimization set on $\Ta$, we have
\[
\norm{\tzab-\tcf{\mean}}_{\oma}
\leq
\norm{(\phi^{\rm d})^{-1}(\htzab_h) - \tcf{\mean}}_{\oma}
\leq
\norm{(\phi^{\rm d})^{-1}}\norm{\htzab_h}_{\widehat \omega} + \norm{\tcf{\mean}}_{\oma},
\]
so that
\[
\norm{\tzab-\tcf{\mean}}_{\oma}
\leq
(1+ C(\kappa_{\Th})\norm{(\phi^{\rm d})^{-1}}\norm{\phi^{\rm d}})\norm{\tcf{\mean}}_{\oma}.
\]
At \cgr{this} point, we conclude the proof since
$\norm{(\phi^{\rm d})^{-1}}\norm{\phi^{\rm d}}$
only depends on $\kappa_{\Th}$ due to~\eqref{eq_bound_piola}
and $\norm{\tzab}_{\oma} \leq \norm{\tcf{\mean}}_{\oma} + \norm{\tzab-\tcf{\mean}}_{\oma}$.
\ep

\section{Decomposition of a divergence-free pie\-ce\-wise polynomial
with an elementwise orthogonality into local divergence-free contributions}
\label{app_div_free_dec}

Let $q \geq 0$ be a fixed integer and recall the notation of Section~\ref{sec_not}; namely,
$\RTproj{q}$ is the canonical $q$-degree Raviart--Thomas interpolate \cye{on the given mesh element $\elm \in \Th$} from~\eqref{eq_RT_proj}\cye{,} and $\ls$ means smaller or equal to up to a constant only
depending on the mesh shape-regularity parameter $\kappa_\Th$. The following result is
of independent interest:

\bt[Decomposition of a divergence-free Raviart--Thomas piecewise polynomial with an elementwise
orthogonality constraint into local divergence-free contributions]
\label{thm_div_free_dec_abs}
Let
\be \label{eq_td_prop_A}
    \bdl_h \in \RT_q(\Th) \cap \HdvN \quad \text{ with } \quad \Dv \bdl_h = 0
\ee
be a divergence-free $q$-degree Raviart--Thomas piecewise polynomial that is elementwise
orthogonal to vector-valued constants,
\be \label{eq_elm_orth_A}
    (\bdl_h, \tr_h)_\elm = 0 \quad \forall \tr_h \in [\PP_0(\elm)]^3, \, \forall \elm \in \Th.
\ee
Then there exists a unique solution to the $\qq$-degree Raviart--Thomas elementwise
minimizations, $\qq = q$ or $\qq = q + 1$,
\be \label{eq_tda_A}
\tda|_\elm
\eq
\arg \min_{\substack{
\tv_h \in \RT_{\qq}(\elm)
\\
\Dv \tv_h= 0
\\
\tv_h \scp \tn_{\elm} = \RTproj{\qq}(\cye{(}\psia \bdl_h\cye{)|_\elm}) \scp \tn_{\elm} \, \text{ on } \pt \elm }}
\norm{\tv_h - \RTproj{\qq}(\cye{(}\psia \bdl_h\cye{)|_\elm})}_\elm^2
\ee
for all tetrahedra $\elm \in \Th$ and all vertices $\ver \in \VK$.
This yields patchwise divergence-free contributions
\be \label{eq_tda_prop_A}
\tda \in \RT_{\qq}(\Ta) \cap \Hdva
\quad \text{ with } \quad
\Dv \tda= 0
\qquad \forall \ver \in \Vh,
\ee
decomposing $\bdl_h$ as
\be \label{eq_tsat_dec_A}
    \bdl_h = \suma \tda.
\ee
Moreover, for all tetrahedra $\elm \in \Th$ and all vertices $\ver \in \VK$, there hold the
local stability estimates
\bse
\label{eq_tsat_stab_A}
\ba
\label{eq_tsat_stab_1_A}
\norm{\tda - \RTproj{\qq}(\cye{(}\psia \bdl_h\cye{)|_\elm})}_\elm & \ls \norm{\bdl_h}_\elm,
\\
\label{eq_tsat_stab_2_A}
\norm{\tda}_\elm & \ls_{\qq} \norm{\bdl_h}_\elm,
\ea
\ese
where $\ls_{\qq}$ means $\ls$ for $\qq = q + 1$ and up to a constant only depending on the
mesh shape-regularity parameter $\kappa_\Th$ and the degree $q$ when $\qq = q$.
\et

\br[The two settings $\qq = q$ or $\qq = q + 1$ in Theorem~\ref{thm_div_free_dec_abs}]
With the choice $\qq = q$, the contributions $\tda$ in Theorem~\ref{thm_div_free_dec_abs}
stay in the same degree Raviart--Thomas space as the datum $\bdl_h$, but, unfortunately,
the stability~\eqref{eq_tsat_stab_2_A} is not necessarily $q$-robust. For $q$-robustness,
the choice $\qq = q + 1$, increasing the degree of $\tda$ by one, is to be used. Note that
in this case, the Raviart--Thomas interpolator $\RTproj{\qq}$ can be disregarded, since then
$\RTproj{\qq}(\cye{(}\psia \bdl_h\cye{)|_\elm}) = \cye{(}\psia \bdl_h\cye{)|_\elm}$.
\er

\bp
Let $\bdl_h$ satisfy~\eqref{eq_td_prop_A} and~\eqref{eq_elm_orth_A}.
We address~\eqref{eq_tda_A}--\eqref{eq_tsat_stab_A} in four steps.

{\em Step~1. Proof of the well-posedness of~\eqref{eq_tda_A}.}
Fix $\elm \in \Th$ and $\ver \in \VK$. The existence and uniqueness of $\tda|_\elm$
from~\eqref{eq_tda_A} are classical following, \eg, \cite{Bof_Brez_For_MFEs_13},
when the Neumann compatibility condition $\<\RTproj{\qq}(\cye{(}\psia \bdl_h\cye{)|_\elm}) \scp \tn_{\elm}, 1\>_{\pt \elm} = 0$ is satisfied. This can be shown via~\eqref{eq_RT_proj_face}, the Green theorem, the assumption $\Dv \bdl_h =0$
in~\eqref{eq_td_prop_A}, and the elementwise orthogonality assumption~\eqref{eq_elm_orth_A}
(note that $(\Gr \psia)|_\elm \in [\PP_0(\elm)]^3$) as
   \[
       \<\RTproj{\qq}(\cye{(}\psia \bdl_h\cye{)|_\elm}) \scp \tn_{\elm}, 1\>_{\pt \elm} = \ifSIAM \else \<\psia \bdl_h \scp \tn_{\elm}, 1\>_{\pt \elm} = \fi \<\bdl_h \scp \tn_{\elm}, \psia\>_{\pt \elm} =
       (\Dv \bdl_h, \psia)_\elm + (\bdl_h, \Gr \psia)_\elm = 0.
   \]

{\em Step~2. Proof of the stability estimates~\eqref{eq_tsat_stab_A}.}
Still for a fixed $\elm \in \Th$ and $\ver \in \VK$, consider the problem
\be\label{eq_tsah_A}
\tdah|_\elm
\eq
\arg \min_{\substack{
\tv_h \in \RT_{\qq}(\elm)
\\
\Dv \tv_h=\cye{(}-\Gr \psia \scp \bdl_h\cye{)|_\elm}
\\
\tv_h \scp \tn_{\elm} = 0 \, \text{ on } \pt \elm
}}
\norm{\tv_h}_\elm^2.
\ee
This problem is again well-posed since, from~\eqref{eq_elm_orth_A}, $(\Gr \psia \scp \bdl_h,1)_\elm = (\bdl_h,\Gr \psia )_\elm = 0$; moreover, $\cye{(}\Gr \psia \scp \bdl_h\cye{)|_\elm} \in \PP_q(\elm) \subset \PP_\qq(\elm)$, since from $\Dv \bdl_h = 0$, it follows that
$\bdl_h|_\elm \in [\PP_q(\elm)]^3$ (see, \eg, \cite[Corollary~2.3.1]{Bof_Brez_For_MFEs_13}).
It follows that $\tdah|_\elm = \tda|_\elm - \RTproj{\qq}\cye{(}(\psia \bdl_h)|_\elm\cye{)}$; indeed, the commuting property~\eqref{eq_com_prop} yields, on the simplex $\elm$,
$
\Dv (\RTproj{\qq}(\psia \bdl_h))
=
\PP_{\qq} (\Dv (\psia \bdl_h))
=
\PP_{\qq} (\Gr \psia \scp \bdl_h)
=
\Gr \psia \scp \bdl_h
$.
Problem~\eqref{eq_tsah_A} fits the framework of~\cite[Lemma~A.3]{Ern_Voh_p_rob_3D_20}
with $r_\sd = 0$, $r_\elm = \cye{(}- \Gr \psia \scp \bdl_h\cye{)|_\elm}$, and $p=\qq$, so that
\[
\norm{\tda - \RTproj{\qq}(\cye{(}\psia \bdl_h\cye{)|_\elm})}_\elm
=
\norm{\tdah}_\elm
=
\!\!\!\!\!\min_{\substack{
\tv_h \in \RT_{\qq}(\elm)
\\
\Dv \tv_h= - \Gr \psia \scp \bdl_h
\\
\tv_h \scp \tn_{\elm} = 0 \, \text{ on } \pt \elm
}}
\!\!\!\!\!\norm{\tv_h}_\elm
\ls \!\!\!\!
\min_{\substack{
\tv \in \Hdvi{\elm}
\\
\Dv \tv = -\Gr \psia \scp \bdl_h
\\
\tv \scp \tn_{\elm} = 0 \, \text{ on } \pt \elm
}}
\!\!\!\!\norm{\tv}_\elm
=
\norm{\Gr \zeta_\elm}_\elm.
\]
Here, by primal--dual equivalence, $\zeta_\elm \in \Hois{\elm}$ is such that
\[
(\Gr \zeta_\elm, \Gr v)_\elm
=
-(\Gr \psia \scp \bdl_h, v)_\elm
\qquad
\forall v \in \Hois{\elm}
\]
with $\Hois{\elm} \eq \{v \in \Hoi{\elm}; \, (v,1)_\elm=0\}$,
where the Poincar\'e inequality gives $\norm{v}_\elm \ls h_\elm \norm{\Gr v}_\elm$.
Then the Cauchy--Schwarz inequality and shape regularity yield
\[\norm{\Gr \zeta_\elm}_\elm
\! = \!\!\!\!
\max_{\substack{ v \in \Hois{\elm}\\ \norm{\Gr v}_\elm =1}}
\!\!\!\!(\Gr \zeta_\elm, \Gr v)_\elm
\! = \!\!
\!\!\max_{\substack{v \in \Hois{\elm}\\ \norm{\Gr v}_\elm =1}}
\!\!\!\!-(\Gr \psia \scp \bdl_h, v)_\elm
\!
\ls \!
\norm{\Gr \psia}_{\infty, \elm} \norm{\bdl_h}_\elm h_\elm
\! \ls \!
\norm{\bdl_h}_\elm.
\]
Combining the two above estimates gives the desired stability
result~\eqref{eq_tsat_stab_1_A}. The other stability result~\eqref{eq_tsat_stab_2_A}
follows from~\eqref{eq_tsat_stab_1_A} by the triangle inequality
together with the non-$q$-robust stability bound
$
\norm{\RTproj{\qq}(\cye{(}\psia \bdl_h\cye{)|_\elm})}_\elm
\ls_{\qq}
\norm{\psia \bdl_h}_\elm
\leq
\norm{\bdl_h}_\elm
$
when $\qq = q$, whereas
$\norm{\RTproj{\qq}(\cye{(}\psia \bdl_h\cye{)|_\elm})}_\elm = \norm{\psia \bdl_h}_\elm \leq \norm{\bdl_h}_\elm$
when $\qq = q + 1$.

{\em Step~3. Proof of the patchwise properties~\eqref{eq_tda_prop_A}.}
The first property in~\eqref{eq_tda_prop_A} follows from the prescription of the normal
components in~\eqref{eq_tda_A}, whereas the second one is the divergence prescription
in~\eqref{eq_tda_A}.

\enlargethispage{0.9cm}

{\em Step~4. Proof of the decomposition~\eqref{eq_tsat_dec_A}.}
Finally, in order to prove~\eqref{eq_tsat_dec_A}, set $\tilde \bdl_h \eq \suma \tda$.
Now fix an element $\elm \in \Th$ and remark that from the normal trace constraint
in~\eqref{eq_tda_A} and the linearity of the interpolator $\RTproj{\qq}$,
on $\pt \elm$,
\ban
\tilde \bdl_h\cye{|_\elm} \scp \tn_{\elm}
& =
\sum_{\ver \in \VK}  \tda\cye{|_\elm} \scp \tn_\elm
=
\sum_{\ver \in \VK}  \RTproj{\qq}(\cye{(}\psia \bdl_h\cye{)|_\elm}) \scp \tn_{\elm}
=
\RTproj{\qq}\Bigg[\sum_{\ver \in \VK}  (\psia \bdl_h)\cye{|_\elm}\Bigg] \scp \tn_{\elm} \\
& =
\RTproj{\qq}(\bdl_h\cye{|_\elm}) \scp \tn_{\elm}
=
\bdl_h\cye{|_\elm} \scp \tn_{\elm}
\ean
also using the partition of unity~\eqref{eq_PU}. Similarly, by the divergence constraint
in~\eqref{eq_tda_A} and $\Dv \bdl_h = 0$ from~\eqref{eq_td_prop_A}, on $\elm$,
\[
\Dv \tilde \bdl_h = \sum_{\ver \in \VK} \Dv \tda = 0 = \Dv \bdl_h.
\]
Consequently, $(\tilde \bdl_h - \bdl_h)|_\elm \in \RT_{\qq}(\elm)$ has zero normal trace and
divergence. Moreover, the Euler conditions of problem~\eqref{eq_tda_A} state
\[
(\tda - \RTproj{\qq}(\cye{(}\psia \bdl_h\cye{)|_\elm}), \tv_h)_\elm = 0
\quad
\forall \tv_h \in \RT_{\qq}(\elm)
\text{ with }
\Dv \tv_h= 0
\text{ and }
\tv_h \scp \tn_{\elm} = 0
\]
on $\pt \elm$. Summing this over all vertices $\ver \in \VK$ and using again the linearity of $\RTproj{\qq}$,
\[
(\tilde \bdl_h - \bdl_h, \tv_h)_\elm = 0
\qquad
\forall \tv_h \in \RT_{\qq}(\elm)
\text{ with }
\Dv \tv_h= 0
\text{ and }
\tv_h \scp \tn_{\elm} = 0 \, \text{ on } \pt \elm,
\]
so that indeed $\tilde \bdl_h = \bdl_h$ on any mesh element $\elm \in \Th$.
\ep

\enlargethispage{0.8cm}

\ifSIAM\bibliographystyle{siamplain}\else\bibliographystyle{acm_mod}\fi
\bibliography{biblio}

\end{document}